\RequirePackage{rotating} 
\documentclass[opre,nonblindrev]{informs3} 

\OneAndAHalfSpacedXI 

\usepackage{endnotes}
\let\footnote=\endnote
\let\enotesize=\normalsize
\def\notesname{Endnotes}%
\def\makeenmark{$^{\theenmark}$}
\def\enoteformat{\rightskip0pt\leftskip0pt\parindent=1.75em
  \leavevmode\llap{\theenmark.\enskip}}

\usepackage{natbib}
 \bibpunct[, ]{(}{)}{,}{a}{}{,}%
 \def\bibfont{\footnotesize}%
 \def\BIBand{and}%

\usepackage[utf8]{inputenc}
\usepackage[english]{babel}
\usepackage[margin=1in]{geometry}
\usepackage{bbding}
\usepackage[shortlabels]{enumitem}
\usepackage[super]{nth}
\usepackage{bm, bbm}
\usepackage{amsmath, amssymb}
\usepackage{mathtools}
\usepackage{subcaption, placeins}
\usepackage{graphicx}
\usepackage[hidelinks]{hyperref}
\usepackage{pgfplots}

\usepackage{algorithm}
\usepackage{algorithmic}
\usepackage{subcaption}
\usepackage{booktabs}
\usepackage[flushleft]{threeparttable}  

\usepackage{multirow} 
\usepackage{placeins}
\usepackage{booktabs}
\usepackage{accents}
\usepackage{graphics}

\usepackage{tikz}
\usepackage{placeins}
\usepackage{float}
\usepackage{graphicx}
\usepackage{accents}
\usepackage{rotating}
\usepackage{graphics}

\usepackage{duckuments}

\usepackage{color}
\newcommand{\blue}{\color{black}}

\definecolor{blue}{rgb}{0,0,1}

\setlist[itemize]{leftmargin=*, nosep}

\usepackage{etoolbox}
\newcommand{\zerodisplayskips}{%
  \setlength{\abovedisplayskip}{7pt}%
  \setlength{\belowdisplayskip}{5pt}%
  \setlength{\abovedisplayshortskip}{0pt}%
  \setlength{\belowdisplayshortskip}{0pt}}
\appto{\normalsize}{\zerodisplayskips}
\appto{\small}{\zerodisplayskips}
\appto{\footnotesize}{\zerodisplayskips}
\usepackage{xfrac}

\TheoremsNumberedThrough     
\ECRepeatTheorems

\EquationsNumberedThrough    

\begin{document}


\RUNAUTHOR{Cory-Wright and Pauphilet}

\RUNTITLE{Sparse PCA With Multiple Components}

\TITLE{Sparse PCA With Multiple Components}

\ARTICLEAUTHORS{%
\AUTHOR{Ryan Cory-Wright}
\AFF{Department of Analytics, Marketing and Operations, Imperial Business School, London, UK\\
ORCID: \href{https://orcid.org/0000-0002-4485-0619}{$0000$-$0002$-$4485$-$0619$}\\ \EMAIL{r.cory-wright@imperial.ac.uk}} 
\AUTHOR{Jean Pauphilet}
\AFF{Management Science and Operations, London Business School, London, UK \\ORCID: \href{https://orcid.org/0000-0001-6352-0984}{$0000$-$0001$-$6352$-$0984$}\\\EMAIL{jpauphilet@london.edu}}
} 

\ABSTRACT{%
Sparse principal component analysis is a {\color{black}fundamental} technique for obtaining {\color{black}interpretable} combinations of features, or principal components (PCs), that explain the variance of high-dimensional datasets. This involves solving a sparsity- and orthogonality-constrained convex maximization problem, which is extremely computationally challenging. 
Most existing work addresses sparse PCA via methods---such as iteratively computing one sparse PC and deflating the covariance matrix---that do not guarantee the orthogonality, let alone the optimality, of the resulting solution when we seek multiple mutually orthogonal PCs. We challenge this status {\color{black}quo} by 
investigating three competitive approaches {\color{black}that} each {\color{black}produce} valid relaxations and high-quality solutions. Together,
these relaxations and their associated heuristic{\color{black}s} {\color{black}yield} solutions with an average bound gap on the order of 3\% 
for real-world datasets with {\color{black}hundreds or thousands} of features and $r \in \{2,3\}$ components. 
The first approach reformulates orthogonality conditions as rank constraints and thereby derives tight semidefinite relaxations, strengthened via additional second-order cone inequalities. 
{The second approach uses Lagrangian decompositions to relax the orthogonality constraints via a penalty term in the objective, thus solving the problem as a sequence of $r=1$ sparse PC problems. }
The third approach is based on a new combinatorial upper bound on the variance explained for a given support {\color{black}pattern, obtained by solving} a mixed-integer linear optimization problem.
Numerically, our algorithms match (and sometimes surpass) the best-performing methods in terms of {\color{black}the} fraction of variance explained and systematically return PCs that are sparse and orthogonal. In contrast, we find that existing methods, such as deflation, return solutions that violate the orthogonality constraints, even when the data {\color{black}are} generated from sparse orthogonal PCs. Altogether, our approaches solve sparse PCA problems with multiple components to certifiable near{\color{black}-}optimality in a practically tractable fashion.
}%

\KEYWORDS{Sparse Principal Component Analysis; Semidefinite Optimization; Practical Tractability} 

\maketitle
\vspace{-10mm}

\section{Introduction}
Principal Component Analysis {\color{black}(}PCA{\color{black})}, {\color{black}first} proposed by \citet{pearson1901liii}, is one of the most popular techniques {\color{black}for reducing} the dimension of a dataset
\citep[see also][]{hotelling1933analysis, eckart1936approximation}. Given a normalized and centered data matrix $\bm{A} \in \mathbb{R}^{n \times p}$, and its sample covariance matrix ${\bm{\Sigma}:=\frac{1}{n-1} \bm{A}^\top\bm{A}}$, the top $r$ principal components ($r \ll p$) of $\bm{\Sigma}$ {\color{black}are obtained} by solving:
\begin{align}\label{prob:pca}
    \max_{\bm{U} \in \mathbb{R}^{p \times r}} \quad & \langle \bm{U}\bm{U}^\top, \bm{\Sigma}\rangle \ \text{s.t.} \ \bm{U}^\top \bm{U}=\mathbb{I}.
\end{align}
As described by \citet{hotelling1933analysis}, the principal components of $\bm{\Sigma}$ correspond to its $r$ leading eigenvectors and can efficiently be obtained via a greedy procedure: at each iteration, {\color{black}one computes the leading eigenvector of $\bm \Sigma$, $\bm{u}$, (for instance, by solving \eqref{prob:pca} with $r=1$) and then updates (or \emph{deflates}) $\bm{\Sigma}$ to eliminate the influence of $\bm{u}$}.
PCA is now a {\color{black}fundamental} unsupervised learning paradigm that is practically useful across a range of fields, including pattern recognition \citep{naikal2011informative}, sequence classification \citep{tan2014classification}, factor models in finance \citep{fan2016projected}, and variable selection in statistics \citep{wang2023variable}.

Despite efficient modern implementations 
\citep{udell2016generalized, tropp2017practical}, PCA suffers from two limitations. First, it generates components that are dense linear combination{\color{black}s} of the original features and hence uninterpretable \citep{rudin2022interpretable}. Second, it yields inconsistent estimates in high-dimensional settings where $p/n \rightarrow \alpha > 0$ 
\citep{johnstone2009consistency}. 
Accordingly, several authors \citep[such as][]{jolliffe2003modified, d2007direct} have proposed sparse PCA, namely{\color{black},} augmenting Problem \eqref{prob:pca} with a sparsity constraint. 
When $r=1$, sparse PCA can be formulated as the following optimization problem:
\begin{align}\label{prob:spcarank1}
    \max_{\bm{u} \in \mathbb{R}^p} \quad & \langle \bm{u}\bm{u}^\top, \bm{\Sigma}\rangle \ \text{s.t.}\ \Vert \bm{u}\Vert_2^2=1, \Vert \bm{u}\Vert_0 \leq k,
\end{align}
where $\Vert \bm{u}\Vert_0$ denotes the cardinality of $\bm{u}$ or the size of its support: $\Vert \bm{u}\Vert_0 = | \operatorname{supp}(\bm{u})| = |\{ j : \bm{u}_j \neq 0 \}|$. 
From a statistical recovery perspective, {\color{black}suppose that} the data {\color{black}are} generated according to a `true' covariance matrix of the form $\bm{\Sigma}^\star = \beta \bm{v} \bm{v}^\top + \bm{I}_p$, for some $\bm{v}$ with $\|\bm{v}\|_0 \leq k$ and some $\beta > 0$. {\color{black}Then,} \citet{amini2008high} {\color{black}have shown} that an exhaustive search algorithm can reliably identify the support of $\bm{v}$ provided that the number of samples {\color{black}satisfies $n \gtrsim k \log(p)$,} which constitutes a significant improvement over the traditional PCA formulation that requires $n \gtrsim p$. {\color{black}Moreover, they have also shown that} no method can succeed when $n \lesssim k \log(p)$, because of information-theoretic limitations. \citet{berthet2013optimal} analyze the support recovery ability of polynomial-time algorithms and show that no polynomial-time algorithm can succeed when $n \lesssim k^2$ {\blue(assuming planted clique cannot be solved in randomized polynomial time)}. Hence, there exists a regime, $k \log(p) \lesssim n \lesssim k^2$, where exhaustive search successfully detects the support of $\bm{v}$, while no polynomial-time algorithm can {\color{black}do so}.

This gap motivated the development of tailored discrete optimization algorithms to efficiently implement exhaustive search for sparse PCA with a single PC. Indeed, Problem \eqref{prob:spcarank1}
can be formulated as a mixed-integer semidefinite optimization problem and solved via global optimization techniques such as branch-and-bound \citep{berk2019certifiably}, or branch-and-cut \citep{bertsimas2020solving}. Confirming the statistical theory, certifiably optimal methods for \eqref{prob:spcarank1} are often significantly more accurate than polynomial-time methods. 

To our knowledge, there is no widely {\color{black}agreed-upon} formulation in the literature that {\color{black}generalizes} Problem \eqref{prob:spcarank1} to $r>1$, and there are no practically relevant algorithms with optimality guarantees that successfully address this extension. 
Indeed, most algorithmic work for sparse PCA with one PC cannot be readily generalized to the case with $r > 1$. This is because, in the multiple component case, the sparsity constraint in \eqref{prob:spcarank1} causes sparse principal components to no longer be eigenvectors of $\bm{\Sigma}$ and deflation methods to no longer lead to orthogonal, let alone optimal, PCs \citep{mackey2008deflation}.
In other words, the conventional PCA wisdom that multiple components can be computed one by one fails as soon as sparsity is required. This observation calls for new methods for sparse PCA with multiple PCs that {\color{black}simultaneously} optimize the components. 
In response, we study a generic optimization formulation that extends Problem \eqref{prob:spcarank1} to $r>1$, and investigate three competitive approaches {\color{black}that produce provably near-optimal solutions with certificates of optimality.}

\subsection{A Generic Formulation for Sparse PCA with Multiple PCs}\label{ssec:genericformulation}
Perhaps the most natural extension of sparse PCA to multiple principal components, and the one which we advocate in this paper, is to augment Problem \eqref{prob:pca} with a constraint on the number of non-zero entries in the matrix $\bm{U}$, $$\| \bm{U} \|_0=\vert \mathrm{supp}(\bm{U})\vert=\vert\{(i,j) \in [p]\times [r]: U_{i,j}\neq 0\}\vert \leq k.$$
This gives a formulation which enforces two desirable properties on the matrix $\bm{U}$: orthogonality ($\bm{U}^\top \bm{U}=\mathbb{I}$), as present in the prototypical formulation of PCA \citep[see, e.g.,][]{johnson1985matrix}; and sparsity, to address interpretability and accuracy concerns \citep[cf.][]{rudin2022interpretable}.

Formally, introducing a binary matrix $\bm{Z}$ to encode the support of $\bm{U}$, we consider the problem:
\begin{align}\label{prob:spcaorth_compact}
    \max_{\substack{\bm{Z} \in \{0, 1\}^{p \times r}}} \max_{\bm{U} \in \mathbb{R}^{p \times r}} \quad \langle \bm{U}\bm{U}^\top, \bm{\Sigma}\rangle \quad
    \text{s.t.} \quad & \bm{U}^\top \bm{U}=\mathbb{I}, 
    \\ &  U_{i,t}=0 \ \text{if} \ Z_{i,t}=0, \ \forall i \in [p], \ \forall t \in [r], \nonumber 
    \\ & \sum_{i \in [p],t \in [r]} Z_{i,t} \leq k \mbox{ or } \sum_{i \in [p]} Z_{i,t} \leq k_t,\ \forall t \in [r], \nonumber
\end{align}
where {\color{black}$[r]$ denotes the {\color{black}index set} $\{1, \ldots, r\}$}. We either impose a global sparsity constraint on $\bm{U}$, $\| \bm{U} \|_0 \leq k$, via $\sum_{i \in [p],t \in [r]} Z_{i,t} \leq k$ (with $ p r>k>r$, in this case, for the problem to be well-posed and non-trivial), or a per-component sparsity requirement, $\| \bm{U}_t \|_0 \leq k_t$, via the constraints $\sum_{i \in [p]} Z_{i,t} \leq k_t$ (with $0 < k_t < p$). We address both modeling options in this paper, although the second option is more common in the literature 
\citep[e.g.,][]{hein2010inverse,deshpande2014information,berk2019certifiably}.

Observe that optimal solutions to Problem \eqref{prob:spcaorth_compact} may not be eigenvectors of a submatrix of $\bm{\Sigma}$, especially when the support of the columns of $\bm{Z}$ is partly overlapping. This fundamental difference from sparse PCA with one PC makes Problem \eqref{prob:spcaorth_compact} {\color{black}substantially more difficult} to solve.

We remark that although Problem \eqref{prob:spcaorth_compact} is a very natural extension of Problem \eqref{prob:spcarank1}, we are not aware of any work{\color{black}s} that explicitly formulate sparse PCA with multiple components as an orthogonality-constrained problem with logical constraints. {\color{black}Nor are we aware of any works that study} \eqref{prob:spcaorth_compact} at this level of generality.
As we review in the next section, many works {\color{black}implicitly} study \eqref{prob:spcaorth_compact} by computing PCs that are simultaneously sparse and explain most of the {\color{black}dataset's variance} without defining them as formal solutions of a cardinality- and orthogonality-constrained optimization problem. In the discrete optimization literature, several authors study special cases of \eqref{prob:spcaorth_compact} with additional constraints on $\bm{Z}$ that ensure it admits a mixed-integer semidefinite reformulation (row sparsity and disjoint sparsity, see below). 
{\blue In practice, however, the support of each PC may neither be disjoint nor fully overlapping. For example, in stock return data, each PC can be interpreted as a long-short portfolio of an industry sector \citep{avellaneda2010statistical}{\color{black},} so we should not expect the supports of each PC to be perfectly disjoint nor fully overlapping.}
Closest to our formulation are PCA problems with orthogonality constraints and an $\ell_1$ penalty term to induce sparsity \citep{zou2006sparse,lu2012augmented,vu2013fantope,benidis2016orthogonal}.

Finally, from a generative model perspective, Problem \eqref{prob:spcaorth_compact} is consistent with a spiked covariance model \citep[see, e.g.,][]{amini2008high,d2008optimal}, where the true covariance matrix $\bm{\Sigma}^\star$ can be decomposed as the sum of a sparse and low-rank term plus some noise:
\begin{align}\label{prob:genmodel}
    \bm{\Sigma}^\star = \sum_{t \in [r]}\beta_t \bm{v}_t \bm{v}_t^{\top} + \bm{N},
\end{align}
where $\bm{v}_t$ are sparse vectors that may have non-overlapping, partially overlapping, or completely overlapping support, and $\bm{N}$ is a noise matrix. This generative model is referred to as the spiked Wishart model when $\bm{N} = \bm{I}_p$, 
and {\color{black}as} the spiked Wigner model when $\bm{N}$ is drawn from the Gaussian orthogonal ensemble \citep[see][for a general theory of both families of models]{ding2023subexponential}.

\subsection{Literature Review}\label{sec:litrev}
To identify the extent to which the state-of-the-art for sparse PCA {\color{black}can} be improved, we now review methods {\color{black}for} approximately solv{\color{black}ing} sparse PCA with multiple PCs. 

\paragraph{Methods for sparse PCA with $r=1$:}
Several polynomial-time algorithms have been proposed to obtain high-quality solutions to \eqref{prob:spcarank1}, including greedy heuristics \citep{d2008optimal}, $\ell_1$ relaxations \citep{zou2006sparse,d2007direct, dey2021using}, linear regression-based estimators \citep{bresler2018sparse,behdin2021sparse}, or thresholding techniques \citep{johnstone2009consistency,deshpande2014sparse}. Note that the covariance thresholding method of \citet{deshpande2014sparse} provably recovers the support in the spiked Wishart model whenever $n \gtrsim k^2$, which is the best achievable rate for polynomial-time methods \citep{berthet2013optimal}.
Over the past decade, various authors including \cite{gally2016computing,bertsimas2020solving,kim2021convexification, li2020exact} have shown that Problem \eqref{prob:spcarank1} can be recast as a mixed-integer semidefinite optimization (MISDO) problem and have derived both high-quality solutions and valid dual bounds using this discrete optimization lens. 

\paragraph{Deflation methods for sparse PCA:} 
It is well known that an optimal solution to Problem \eqref{prob:pca} can be obtained via a greedy deflation procedure. 
Consequently, \cite{mackey2008deflation} propose{\color{black}s} a sparse extension of this deflation scheme where, at each iteration, a sparse PC is computed (e.g., by solving \eqref{prob:spcarank1} or a relaxation) and $\bm{\Sigma}$ is updated by projecting out the eigenspace modeled by $\bm{u}$: $\bm{\Sigma}_{\text{new}}=(\mathbb{I}-\bm{u}\bm{u}^\top)\bm{\Sigma}(\mathbb{I}-\bm{u}\bm{u}^\top)$. Empirically, this method often performs reasonably well \citep[see][Section 5.3]{berk2019certifiably}, particularly when Problem \eqref{prob:spcarank1} is solved to global optimality; see also \citet[][]{hein2010inverse, buhler2014flexible} for a related deflation-based scheme. 
However, unlike in the traditional case, deflation need not return an optimal solution to \eqref{prob:spcaorth_compact}---see \citet{asteris2015sparse} for a simple four-dimensional example.
{\color{black}Moreover}, deflation need not even return a feasible solution, since the orthogonality constraint is not explicitly imposed by the method, and is therefore often violated in practice. Because they cannot guarantee the feasibility (i.e., orthogonality) of the returned PCs, we consider deflation-based methods as heuristic methods for our sparse PCA problem with multiple PCs.
Furthermore, because of their iterative nature, deflation-based procedures can easily control the sparsity of each component but usually struggle to enforce global sparsity.

\paragraph{Methods for generic sparse PCA:}
A second approach for solving Problem \eqref{prob:spcaorth_compact} is to apply a heuristic {\color{black}that} approximately optimizes all $r$ PCs simultaneously, rather than sequentially. Among others, \citet{zou2006sparse} propose an alternating maximization scheme for an $\ell_1$ relaxation of Problem \eqref{prob:spcaorth_compact}, \citet{journee2010generalized} propose an iterative conditional gradient method to identify a local optimum of Problem \eqref{prob:spcaorth_compact} without the orthogonality constraints, \cite{lu2012augmented} apply an augmented Lagrangian method which solves an $\ell_1$ relaxation of \eqref{prob:spcaorth_compact}, \cite{vu2013fantope} solve a semidefinite relaxation of Problem \eqref{prob:spcaorth_compact}'s $\ell_1$ relaxation, and \cite{benidis2016orthogonal} adopt a minorization-maximization approach which also solves Problem \eqref{prob:spcaorth_compact} approximately. Unfortunately, these approaches are often suboptimal {\color{black}even when} $r=1$, and {\color{black}provide} no indication {\color{black}of} the optimality gap. Indeed, none of these approaches explicitly control both the sparsity and orthogonality constraints, and therefore none of the methods reviewed here are guaranteed to return feasible solutions to \eqref{prob:spcaorth_compact}. {\color{black}N}ote that the covariance thresholding method of \citet{deshpande2014sparse} can be applied to the $r>1$ case and can return PCs that are asymptotically orthogonal. However, for a given dataset, {\color{black}covariance thresholding} cannot guarantee the orthogonality of the returned solution, as we observe empirically in Sections \ref{ssec:feasiblemethods} and \ref{ssec:auc}.

\paragraph{Row Sparsity:} Motivated by tractability concerns, 
another line of work studies a special case of \eqref{prob:spcaorth_compact}, namely
row-sparse principal component analysis or principal component analysis with global support. This formulation replaces the sparsity constraint on $\bm{U}$ with one requiring {\blue that the matrix} $\bm{U}$ has at most $k$ non-zero rows, as advocated by \cite{boutsidis2011sparse,probel2011large, vu2013minimax}. This rewrites sparse PCA as performing a top-$r$ SVD on {\color{black}the $k \times k$ principal submatrix} of $\bm{\Sigma}$, i.e.,
\begin{align}\label{prob:globalsupport}
    \max_{\bm{z} \in \{0, 1\}^p: \bm{e}^\top \bm{z} \leq k}\ \max_{\bm{U} \in \mathbb{R}^{p \times r}} \quad & \langle \bm{U}\bm{U}^\top, \bm{\Sigma}\rangle\ \text{s.t.} \ \bm{U}^\top \bm{U}=\mathbb{I}_{r \times r}, \ U_{i,t}=0 \ \text{if} \ z_{i}=0\ \forall i \in [p], \forall t \in [r].
\end{align}
Problem \eqref{prob:globalsupport} is a special case of Problem \eqref{prob:spcaorth_compact} where each PC has the same support, which corresponds to constraining $\bm{Z}$ such that {\color{black}$Z_{i,t_1}=Z_{i,t_2}, \ \forall t_1, t_2 \in [r], \forall i \in [p]$} in \eqref{prob:spcaorth_compact}, and dividing the ``$k$'' in \eqref{prob:spcaorth_compact} by $r$. This restriction is advantageous from a computational perspective, but disadvantageous from a statistical one. Indeed, because of the global support assumption, {\blue Problem \eqref{prob:globalsupport} can be reformulated as a mixed-integer semidefinite problem \citep[see][for derivations]{bertsimas2020solving, li2021beyond} and solved via relax-and-round and local search strategies \citep{li2021beyond}, mixed-integer linear approximations \citep{li2021beyond,bertsimassparse2022}, mixed-integer second-order-cone approximations \citep{dey2020solving}, or branch-and-cut \citep{li2021beyond}. Although derived for theoretical purposes and not implemented (let alone evaluated computationally),  \citet{del2022sparse} also describes a partitioning procedure to solve \eqref{prob:globalsupport}.}
However, {\color{black}from the perspective of the generative model in} \eqref{prob:genmodel}, it is equivalent to making the very strong assumption that all leading eigenvectors $\bm{v}_t$ have the same sparsity pattern. Indeed, if one PC $\bm{v}_t$ is sparser than $k$ ($\| \bm{v}_t \|_0 \leq k_t < k$) or two PCs {\color{black}have partially overlapping support} then Problem \eqref{prob:globalsupport} will nonetheless estimate that the support of each {\blue column $\bm{U}_t$} is of size $k$. Therefore, Problem \eqref{prob:globalsupport} is vulnerable to consistently making false discoveries in the identification of relevant features. 

\paragraph{Disjoint Sparsity:}
Another relevant special case of sparse PCA with multiple PCs is when 
the supports of all $r$ columns of $\bm{U}$ are assumed mutually disjoint, as originally proposed by \cite{asteris2015sparse}. This gives rise to the formulation:
\begin{align}\label{prob:spcadisjoint}
    \max_{\substack{\bm{Z} \in \{0, 1\}^{p \times r}: \\\langle \bm{E},\ \bm{Z}\rangle \leq k, \bm{Z}\bm{e} \leq \bm{e}}} \max_{\bm{U} \in \mathbb{R}^{p \times r}} \quad & \langle \bm{U}\bm{U}^\top, \bm{\Sigma}\rangle\ \text{s.t.} \ \bm{U}^\top \bm{U}=\mathbb{I}_{r \times r}, \ U_{i,t}=0 \ \text{if} \ Z_{i,t}=0\ \forall i \in [p], \, t \in [r],
\end{align}
where $\bm{E}$ denotes a matrix of all ones of the appropriate dimension.
Note that \eqref{prob:spcadisjoint} is a special case of \eqref{prob:spcaorth_compact} where we additionally require that $\sum_{t \in [r]}Z_{i,t} \leq 1 \ \forall i \in [p]$, i.e., that each feature can be included in at most one PC. Interestingly, this restriction allows \eqref{prob:spcadisjoint} to be recast as {\color{black}an} MISDO and solved as such \citep[cf.][]{bertsimas2020solving}, although we are not aware of any work that has exploited this MISDO reformulation. 
An obvious criticism of formulation \eqref{prob:spcadisjoint} is that, in practice, we may wish to include a feature in multiple PCs. Therefore, \eqref{prob:spcadisjoint} is best thought of as a special case of sparse PCA. 
Indeed, if the true generative model involves vectors $\bm{v}_t$ with partially overlapping supports, then \eqref{prob:spcadisjoint} {\color{black}cannot} recover them. Nonetheless, as we explore in our numerical experiments (Section \ref{sec:numres}), disjoint solutions often perform well when $k$ is very small relative to $p$. This is because, when $k \ll p$, 
there are often several disjoint submatrices {\color{black}that} are near-optimal in the rank-one case, and selecting the leading PCs from each of them is often a reasonable approach in practice.

\paragraph{Summary of Sparsity Constraints:}
We now summarize and contrast the different types of sparsity patterns considered throughout the literature and in this paper {\color{black}(see Table \ref{tab:sparsitytypes})}. 
While special cases of Problem \eqref{prob:spcaorth_compact} admit formulations as mixed-integer semidefinite optimization problems (MISDO), sparse PCA with multiple PCs is a mixed-integer low-rank optimization problem (MIRO) in general, as we establish in the next section. {\color{black}Moreover, even when $\bm{Z}$ is fixed, we show that estimating the coefficients of $\bm{U}$ in Problem \eqref{prob:spcaorth_compact} is NP-hard (Theorem \ref{thm:resolve.nphard}), while disjoint sparse PCA and row sparse PCA can both be solved in polynomial time once $\bm{Z}$ is fixed. Furthermore, for a fixed $\bm{Z}$, semidefinite relaxations of disjoint or row sparse PCA are tight, while they need not be for full sparse PCA (Section \ref{sec:nonredundancyrankconstraints}). This suggests that (full) sparse PCA is computationally harder than its special cases commonly studied in the literature (assuming $P \neq NP$). 
}

\begin{table}[h]
\centering\footnotesize
\caption{Types of sparsity for sparse PCA with multiple principal components. We formulate the constraints on $\bm{Z}$ assuming a global sparsity budget $k$. Similar constraints with a per-component budget $k_t$ could be derived.} 
\label{tab:sparsitytypes}
\begin{tabular}{@{}l l l l l@{}} \toprule
Name & Constraints on $\bm{Z} \blue \in \{0, 1\}^{p \times r}$ &  Formulation & {\color{black}Complexity when $\bm{Z}$ fixed} & Reference\\
    \midrule
Disjoint Sparsity & $Z_{i,t_1}+Z_{i,t_2}\leq 1 \ \forall i \in [p], \ \forall t_1\neq t_2 \in [r]$ & MISDO & {\color{black}P-time}& \cite{asteris2015sparse}\\
& $\sum_{i \in [p], t \in [r]} Z_{i,t}\leq k$ & & \\[1em]
Row Sparsity & $Z_{i,t_1}=Z_{i,t_2} \ \forall i \in [p], \ \forall t_1, t_2 \in [r], $ & MISDO & {\color{black}P-time} & \cite{boutsidis2011sparse}\\
& $\sum_{i \in [p], t \in [r]} Z_{i,t}\leq k$ & & \\[1em]
Sparsity & $\sum_{i \in [p], t \in [r]} Z_{i,t}\leq k$  & MIRO & {\color{black}NP-hard} & This paper\\
\bottomrule
\end{tabular}
\end{table}

\subsection{Contributions and Structure}
To our knowledge, no existing algorithm solves sparse PCA problems with multiple components and obtains certificates of optimality, except in the aforementioned special cases of row sparsity or disjoint support. Accordingly, we undertake a detailed study of Problem \eqref{prob:spcaorth_compact} in its full generality. 

{\color{black}This paper develops and evaluates} three alternative strategies for obtaining both upper and lower bounds {\color{black}on Problem \eqref{prob:spcaorth_compact}}, i.e., feasible solutions with optimality certificates. 
\begin{itemize}
    \item In Section \ref{sec:relaxations}, we propose a semidefinite relax-and-round scheme (Algorithm \ref{alg:greedymethod2}). To do so, we first derive a sparse and low-rank reformulation of Problem \eqref{prob:spcaorth_compact}. This reformulation differs from the single sparse PC case, where the rank constraints are redundant. Further, we propose a semidefinite relaxation strengthened by second-order cone valid inequalities. Finally, we propose an algorithm {\blue(Algorithm~\ref{alg:greedymethod2})} that rounds an optimal solution to this semidefinite relaxation into $r$ sparse PCs with disjoint support, hence giving orthogonal PCs by default. {\blue Even when the sparsity pattern of each PC is fixed, we show that computing the coefficients of each PC can be NP-hard (Theorem~\ref{thm:resolve.nphard}), justifying the need to restrict our attention to disjoint supports for Algorithm~\ref{alg:greedymethod2} to be tractable.}
    \item In Section \ref{sec:deflation}, we develop a Lagrangian alternating maximization method (Algorithm \ref{alg:alternatingmin}). First, we derive a 
    relaxation of Problem \eqref{prob:spcaorth_compact} by {\blue penalizing the orthogonality constraints} and demonstrate that this gives rise to a{\color{black}n} $r$-factor upper bound on the overall objective, which is computable by solving $r$ different single component sparse PCA problems to optimality. Second, we leverage our relaxation to design an alternating maximization heuristic, {\blue which provides asymptotically feasible solutions}. 
    \item In Section \ref{sec.gen.gershgorin}, we propose a combinatorial relax-and-round approach (Algorithm \ref{alg:disjoint.linalg}). We first generalize the Gershgorin Circle Theorem to multiple PCs with sparsity constraints to obtain a mixed-integer linear upper bound on Problem \eqref{prob:spcaorth_compact}'s objective value as a function of the support of the PCs $\bm{Z}$ only. We then generate feasible solutions from this relaxation by focusing on solutions with disjoint supports. 
\end{itemize}

Finally, in Section \ref{sec:numres}, we thoroughly evaluate the quality of these three approaches on a collection of UCI and synthetic datasets. In terms of upper bounds, we find that the semidefinite approach provides the strongest bounds (within $2$--$5\%$ of the optimal solutions for the small instances that can be solved with a commercial solver) and can be solved in minutes for $p \approx 100$s. However, for larger problem sizes, the Lagrangian relaxation provides the best tightness/tractability trade-off. In terms of generating feasible solutions, our three algorithms generate systematically feasible solutions, while existing approaches routinely violate the orthogonality requirements. In addition, they also improve in terms of the fraction of variance explained, and collectively generate solutions with bound gaps of $<3\%$ for $p$ up to 1000. 

Thus, the methods developed in this paper provide near-optimal solutions of the sparse PCA problem with multiple components in minutes for problems with $p \approx 1000$s, together with certificate{\color{black}s} of their near optimality, providing significant benefit compared with the state-of-the-art. 
For example, on the \verb|pitprops| dataset, with 6 PCs, we explain $81\%$ of the variance with an overall sparsity of $24${\color{black}, whereas} previous studies \citep{lu2012augmented} could explain less than $70\%$ of the variance with twice as many variables.

\subsection{Preliminaries and Notation}
We let nonbold face characters such as $u$ denote scalars, lowercase bold-faced characters such as $\bm{u}$ denote vectors, uppercase bold-faced characters such as $\bm{U}$ denote matrices, and calligraphic uppercase characters such as $\mathcal{Z}$ denote sets. If $\bm{U}$ is a matrix then $\bm{U}_t$ denotes the $t$th column vector of $\bm{U}$, and $U_{i,t}$ denotes the $(i,t)$th entry of $\bm{U}$. We let $[p]$ denote the {\color{black}index set} $\{1, ..., p\}$. We let $\mathbf{e}$ denote a vector of all $1$'s, $\bm{0}$ denote a vector of all $0$'s, {\color{black}$\bm{E}$ denote a matrix of all ones, and} $\mathbb{I}$ denote the identity matrix, with dimension implied by the context.

We also use an assortment of matrix operators. We let $\langle \cdot,\cdot \rangle$ denote the Euclidean inner product between two matrices, $\Vert \cdot \Vert_F$ denote the Frobenius norm, and $\mathcal{S}_+^p$ denote the $p \times p$ positive semidefinite cone; see \citet{johnson1985matrix} for a general theory of matrix operators. 

Further, we use some basic properties of orthogonal projection matrices. Let $\mathcal{Y}_n:=\{\bm{Y} \in \mathcal{S}^n: \bm{Y}^2 =\bm{Y}\}$ denote the set of $n \times n$ projection matrices and $\mathcal{Y}_n^k:=\{\bm{Y} \in \mathcal{Y}_n: \mathrm{tr}(\bm{Y}) \leq k\}$ denote projection matrices with rank at most $k$: note that $\mathrm{Rank}(\bm{Y})=\mathrm{tr}(\bm{Y})$ for any projection matrix $\bm{Y}$. Among others, the convex hulls of $\mathcal{Y}_n$ and $\mathcal{Y}^k_n$ are well-studied, as we now remind the reader:

\begin{lemma}{\citep[Theorem 3 of][]{overton1992sum}}\label{lemma:convhull}
$\mathrm{Conv}(\mathcal{Y}_n)=\{\bm{P}: 0 \preceq \bm{P} \preceq \mathbb{I}\}$ and $\mathrm{Conv}(\mathcal{Y}_n^k)=\{\bm{P}: 0 \preceq \bm{P} \preceq \mathbb{I}, \mathrm{tr}({\bm{P}}) \leq k\}$. Moreover, the extreme points of $\mathrm{Conv}(\mathcal{Y}_n)$ are $\mathcal{Y}_n$, and the extreme points of $\mathrm{Conv}(\mathcal{Y}_n^k)$ are $\mathcal{Y}_n^k$.
\end{lemma}

\section{Semidefinite Optimization Approach}
\label{sec:relaxations}
In this section, we reformulate Problem \eqref{prob:spcaorth_compact} as a mixed-integer low-rank problem, study its semidefinite and second-order cone relaxations, and propose valid inequalities for strengthening them. 

\subsection{An Extended Formulation With Binary and Low-Rank Variables}\label{ssec:reform_lowrank}
Our sparse PCA formulation \eqref{prob:spcaorth_compact} is a mixed-integer quadratic optimization problem 
that exhibits three primary sources of difficulty. First, as {\color{black}is} typical in PCA problems, Problem \eqref{prob:spcaorth_compact} maximizes a convex quadratic function in the decision variable $\bm{U}$. Second, there is a sparsity constraint. 
Finally, and more consequentially, there is an orthogonality constraint. To our knowledge, existing generic non-convex solvers such as \verb|Gurobi| cannot optimize over such orthogonality constraints at a scale of {\color{black}$p$ in the hundreds of} features.

To address the three aforementioned difficulties, we now derive an orthogonality-free reformulation in five steps. First, we introduce the matrix $\bm{Y} = \bm{U} \bm{U}^\top \preceq \mathbb{I}$, thus linearizing the non-convex objective. Second, we introduce rank-one matrices $\bm{Y}^t$ to model the outer products of each column of $\bm{U}$, $\bm{U}_t$, with itself, $\bm{U}_t\bm{U}_t^\top$. Third, we reassign the indicator variable $Z_{i,t}$ to model whether $Y_{i,j}^t$, rather than $U_{i,t}$, is non-zero. Fourth, by letting $\bm{Y}=\sum_{t=1}^r \bm{Y}^t$, we observe that we can omit the matrix $\bm{U}$ (and the constraints involving $\bm{U}$) without altering the set of feasible $\bm{Y}$'s.
Finally, we use the fact that $\bm{Y}_{i,j}^t$ is only supported on indices $i,t$ where $Z_{i,t}>0$ to strengthen the constraint $\bm{Y} \preceq \mathbb{I}$ to $\bm{Y} \preceq \mathrm{Diag}\left(\min\left(\bm{e}, \sum_t \bm{Z}_t\right)\right)${\color{black}, where the minimization operator acts componentwise}. We have:

\begin{theorem}\label{thm:reform}
Problem \eqref{prob:spcaorth_compact} with the global sparsity constraint $\langle \bm{E}, \bm{Z}\rangle \leq k$ attains the same optimal objective value as the problem:
\begin{align}\label{prob:spca_extended}
    \max_{\substack{\bm{Z} \in \{0, 1\}^{p \times r}:\\ \langle \bm{E}, \bm{Z}\rangle \leq k}}\max_{\bm{Y} \in \mathcal{S}^p,
    \bm{Y}^t \in \mathcal{S}^p_+} \  \langle \bm{Y}, \bm{\Sigma}\rangle \ \text{s.t.} \quad & \bm{Y} \preceq \mathrm{Diag}\left(\min\left(\bm{e}, \sum_t \bm{Z}_t\right)\right), \bm{Y}=\sum_{t=1}^r \bm{Y}^t,\\
    & \mathrm{tr}(\bm{Y}^t)=1, \ \forall t \in [r], \ Y_{i,j}^t=0 \ \text{if} \ Z_{i,t}=0, \ \forall t \in [r], i,j \in [p], \nonumber\\
    & \mathrm{Rank}(\bm{Y}^t)=1, \ \forall t \in [r].\nonumber
\end{align}
\end{theorem}
\begin{remark}
We do not explicitly require $\bm{Y} \succeq \bm{0}$, as $\bm{Y}$ is the sum of positive semidefinite matrices. 
\end{remark}

\begin{remark}\label{rem:nonredundancyrankconstraints}
Problem \eqref{prob:spca_extended} is not a convex mixed-integer semidefinite optimization problem because of the presence of rank constraints which are not mixed-integer convex representable as proven by \cite{lubin2022mixed}. Indeed, the rank constraints 
are {not} redundant and {cannot} be dropped from the formulation \eqref{prob:spca_extended} without altering its optimal value. We demonstrate this via a simple example in Section \ref{sec:nonredundancyrankconstraints}. This is notably different from the special cases of \eqref{prob:spca_extended} reviewed in the introduction, where $\bm{Z}$ is restricted to have fully overlapping or 
fully disjoint supports and the rank constraints are redundant. 
\end{remark}

The proof of Theorem \ref{thm:reform} requires an intermediate result (Proposition \ref{prop:orth.to.rank}). Proposition \ref{prop:orth.to.rank} shows that by imposing rank-one constraints on each $\bm{Y}^t$, the condition that the $\bm{Y}^t$'s are mutually orthogonal can be reformulated as a linear semidefinite constraint
(proof deferred to Section \ref{ec.orth.to.rank}).

\begin{proposition}\label{prop:orth.to.rank} Consider $r$ matrices, $\bm{Y}^t \in \mathcal{S}^p_+$, such that $\operatorname{tr}(\bm{Y}^t) = 1$ and $\operatorname{Rank}(\bm{Y}^t) = 1$. Then, ${\sum_{t \in [r]} \bm{Y}^t \preceq \mathbb{I}}$ if and only if $\langle \bm{Y}^t, \bm{Y}^{t'} \rangle = 0\ \forall t, t' \in [r]:\ t \neq t'$.
\end{proposition}

\proof{Proof of Theorem \ref{thm:reform}}
It suffices to show that for any feasible solution to \eqref{prob:spcaorth_compact}, we can construct a feasible solution to Problem \eqref{prob:spca_extended} with an equal or greater payoff, and vice versa.
\begin{itemize}
    \item Let $(\bm{U}, \bm{Z})$ be a solution to Problem \eqref{prob:spcaorth_compact}. Then, since $U_{i,t}$ can only be non-zero if $Z_{i,t}=1$ and $\bm{U}^\top \bm{U}\preceq \mathbb{I}$, it follows that $\bm{U}\bm{U}^\top \preceq \mathrm{Diag}\left(\min\left(\bm{e}, \sum_t \bm{Z}_t\right)\right)$. Therefore, $(\bm{Y}:=\bm{U}\bm{U}^\top, \bm{Y}^t:=\bm{U}_t \bm{U}_t^\top, \bm{Z})$ is a feasible solution to \eqref{prob:spca_extended} with an equal cost. 
    \item Let $(\bm{Y}, \bm{Y}^t, \bm{Z})$ denote a feasible solution to Problem \eqref{prob:spca_extended}. Then, since each $\bm{Y}^t$ is symmetric and rank-one, we can decompose $\bm{Y}^t$ as $\bm{Y}^t=\bm{U}_t \bm{U}_t^\top$ for a vector $\bm{U}_t$ such that $U_{i,t}=0$ if $Z_{i,t}=0$, and concatenate these vectors $\bm{U}_t$ into a matrix $\bm{U}$ such that $(\bm{U}, \bm{Z})$ has the same cost in \eqref{prob:spcaorth_compact} as $(\bm{Y}, \bm{Y}^t, \bm{Z})$ does in \eqref{prob:spca_extended}. Therefore, it remains to show that $\bm{U}^\top \bm{U}=\mathbb{I}$. To see this, {observe that $\bm{Y} \preceq \mathbb{I}$ implies $(\bm{U}_t^\top \bm{U}_{t'})^2 = \langle \bm{Y}^t, \bm{Y}^{t'}\rangle = 0$ if $t\neq t'$ by Proposition \ref{prop:orth.to.rank}.}
    \hfill \Halmos
    \end{itemize}
\endproof

Theorem \ref{thm:reform} provides a formulation that is less compact than \eqref{prob:spcaorth_compact} but contains rank constraints rather than orthogonality constraints. Therefore, it is amenable to exact approaches for addressing sparsity \citep{bertsimas2019unified} and rank \citep{bertsimas2020mixed} constraints.

\subsection{Semidefinite Relaxation With Global Sparsity}\label{ssec:valid.inequalities}
We leverage this mixed-integer low-rank reformulation of Problem \eqref{prob:spcaorth_compact} to derive a semidefinite relaxation. The valid inequalities we derive in this section rely {\color{black}on} the presence of a global sparsity constraint $\langle \bm{E}, \bm{Z} \rangle \leq k$. In particular, per-component sparsity constraints ($\sum_{i \in [p]} Z_{i,t} \leq k_t$) imply a global sparsity constraint $\langle \bm{E}, \bm{Z} \rangle \leq \sum_{t \in [r]}k_t$, so our results apply to both cases. However, tighter relaxations can be obtained when per-component sparsity is specified, as we investigate in the following section.

By relaxing the {\color{black}rank-one constraints on $\bm{Y}^t$ and the integrality constraints on $\bm{Z}$} in \eqref{prob:spca_extended}, we get a first semidefinite relaxation:
\begin{align}\label{prob:disjunctiverelax}
    \max_{\substack{\bm{Z} \in [0, 1]^{p \times r}:\\ \langle \bm{E}, \bm{Z}\rangle \leq k}} \: \max_{\substack{\bm{Y} \in \mathcal{S}^p, \bm{Y}^t \in \mathcal{S}^p_+, \\ \bm{w} \in [0, 1]^p}} \ & \langle \bm{Y}, \bm{\Sigma}\rangle &\\
    \mbox{s.t.} \quad & \bm{Y} \preceq \mathrm{Diag}(\bm{w}), \bm{Y}=\sum_{t=1}^{r} \bm{Y}^t, \mathrm{tr}(\bm{Y}^t) = 1, \bm{w} \leq \bm{Z}\bm{e} &  
    \nonumber\\
    & \vert Y_{i,j}^t\vert  \leq M_{i,j} Z_{i,t} & \forall i, j \in [p], t \in [{\color{black}r}],\nonumber
\end{align}
where the variable $w_i$ models $\min(1, \sum_{t=1}^r Z_{i,t})$ and the $M$-constants are set to $M_{i,i}=1$ and $M_{i,j}=1/2$ if $i \neq j$---this is an upper bound on $\vert Y_{i,j}^t\vert$ because $\bm{Y}^t$ was, before relaxing the rank constraint, a rank one matrix \citep[cf.][]{bertsimas2020solving}. 

We propose to strengthen the relaxation \eqref{prob:disjunctiverelax} via the following valid inequalities:
\begin{theorem}\label{thm:validinequalities} 
Consider a feasible solution to the mixed-integer low-rank problem \eqref{prob:spcaorth_compact}. 
Then, the following inequalities hold:
\begin{align*}
   &  \sum_{j=1}^p (Y^{t}_{i,j})^2\leq Y^t_{i,i}Z_{i,t} & \forall i \in [p], t \in [r],\\
    &  \sum_{j=1}^p {Y^2_{i,j}}\leq r Y_{i,i}w_i & \forall i \in [p],\\
    &   \left(\sum_{j=1}^p \vert Y_{i,j}\vert\right)^2 \leq k Y_{i,i}w_i & \forall i \in [p],\\
    &  \sum_{i \in [p]: i \neq j}Y_{i,j}^2 \leq (k-r+1) w_j ( w_j-Y_{j,j}), & \ \forall j \in [p]. 
\end{align*}
Note that these constraints are rotated second-order cone constraints \citep[see, e.g.,][]{alizadeh2003second}.
\end{theorem}

The proof of Theorem \ref{thm:validinequalities} is detailed in Section \ref{sec:ec.validineq}. We remark that, out of the four groups of valid inequalities, the first inequality has been previously stated in the case of sparse PCA with one component by 
\cite{bertsimas2020polyhedral, bertsimas2020solving, li2020exact}, but the three other groups of valid inequalities are, to our knowledge, new. 

Based on Theorem \ref{thm:validinequalities}, we {\color{black}obtain} the following semidefinite relaxation:
\begin{align}\label{prob:relax_ext_str}
    \max_{\substack{\bm{Z} \in [0, 1]^{p \times r}:\\ \langle \bm{E}, \bm{Z}\rangle \leq k}} \: \max_{\substack{\bm{Y} \in \mathcal{S}^p, \bm{Y}^t \in \mathcal{S}^p_+, \\ \bm{w} \in [0, 1]^p}} \ & \langle \bm{Y}, \bm{\Sigma}\rangle &\\
    \mbox{s.t.} \quad & \bm{Y} \preceq \mathrm{Diag}(\bm{w}), \bm{Y}=\sum_{t=1}^{r} \bm{Y}^t, \mathrm{tr}(\bm{Y}^t) = 1, \bm{w} \leq \bm{Z}\bm{e} &  
    \nonumber\\
    & \vert Y_{i,j}^t\vert  \leq M_{i,j} Z_{i,t} & \forall i, j \in [p], t \in [{\color{black}r}],\nonumber\\
   &  \sum_{j=1}^p {Y^{t}_{i,j}}^2\leq Y^t_{i,i}Z_{i,t} & \forall i \in [p], t \in [r],\nonumber\\
    &  \sum_{j=1}^p {Y^2_{i,j}}\leq r Y_{i,i}w_i & \forall i \in [p],\nonumber\\
    &   \left(\sum_{j=1}^p \vert Y_{i,j}\vert\right)^2 \leq k Y_{i,i}w_i & \forall i \in [p],\nonumber\\
    &  \sum_{i \in [p]: i \neq j}Y_{i,j}^2 \leq (k-r+1) w_j ( w_j-Y_{j,j}), & \ \forall j \in [p],\nonumber 
\end{align}

\subsection{Stronger Relaxation with a Per-Component Sparsity Budget}\label{sec:validineqk}
In this section, we further strengthen our relaxation when specifying a sparsity budget $k_t$ for each component $\bm{Y}^t$. 
This can be understood as allocating the total sparsity on $\bm{U}$ between the different columns of $\bm{U}$ via $\bm{Y}_t$, which models the outer product $\bm{U}_t \bm{U}_t^\top$.

Formally, we have the following result (proof deferred to Section \ref{ssec:append.proof_prop_validineqskt}):
\begin{proposition}\label{prop:validineqskt}
Suppose that $\sum_{i \in [p]}\bm{Z}_{i,t} \leq k_t$ in Problem \eqref{prob:spca_extended}. Then, the following inequalities hold:
\begin{align}\label{ell1_split}
\left(\sum_{j=1}^p \vert Y_{i,j}^t\vert\right)^2 \leq k_t Y^t_{i,i}Z_{i,t} \quad & \forall i \in [p], \forall t \in [r],\\
   \sum_{i \in [p]: i \neq j}(Y_{i,j}^t)^2 \leq (k_t-1) Z_{j,t}( Z_{j,t}-Y_{j,j}^t)  \quad & \forall j \in [p]. \label{eqref:soc2}
\end{align}
\end{proposition}

Interestingly, as we observe in Section \ref{ssec:bounds}, combining Problem \eqref{prob:relax_ext_str} with constraints \eqref{ell1_split}-\eqref{eqref:soc2} often yields much tighter upper bounds than \eqref{prob:relax_ext_str} alone, even if we take the worst-case upper bound over all feasible splits {\color{black}$\{k_t\}_{t \in [r]}$} which sum to $k$. 

Now that we {\color{black}have} introduced the sparsity of each PC, $k_t$, we can
further tighten our semidefinite relaxations by leveraging valid inequalities obtained in the case of a single PC. For example,
for each component $t \in [r]$, \cite{kim2021convexification} observe that the feasible set of all $k_t$-sparse components $\{ \bm{u} \in \mathbb{R}^p \: : \: \| \bm{u} \|_0 \leq k_t\}$ is permutation and sign invariant, i.e., for any feasible vector $\bm{u}$, any vector obtained by permuting or changing the sign of the coordinates of $\bm{u}$ is also feasible. Based on this observation, they propose a lifted formulation for sparse PCA with a single PC, which, to the best of our knowledge, 
leads to the strongest known relaxation for sparse PCA with $r=1$ which can be solved in polynomial time. 
For the sake of concision, we denote $(\bm{Y}^t, \bm{Z}_t) \in \mathcal{T}(k_t)$ the set of valid inequalities comprised in their ``T-relaxation'' (see Problem \eqref{prob:disjunctiverelax_permutationinvariant.full} in Section \ref{ssec:ec.sdp} for an explicit formulation) and  consider the following relaxation:
\begin{align}
    \max_{\substack{\bm{Z} \in [0, 1]^{p \times r}:\\ \langle \bm{E}, \bm{Z}\rangle \leq k, \\ \bm{w} \in [0, 1]^p}
    }\max_{{\bm{Y} \in \mathcal{S}^p_+, \bm{Y}^t \in \mathcal{S}^p_+}} \quad & \langle \bm{Y}, \bm{\Sigma}\rangle \label{prob:disjunctiverelax_permutationinvariant}\\
    \text{s.t.} \quad & \bm{Y} \preceq \mathrm{Diag}(\bm{w}),\ \bm{Y}=\sum_{t=1}^{r} \bm{Y}^t,\ \mathrm{tr}(\bm{Y}^t)=1,\ \bm{w} \leq \bm{Z}\bm{e}, \nonumber \\ 
    & \sum_{j=1}^p Y_{i,j}^2\leq r Y_{i,i}w_i & \forall i \in [p],\nonumber\\
   & \left(\sum_{j=1}^p \vert Y_{i,j}\vert\right)^2 \leq k Y_{i,i}w_i & \forall i \in [p],\nonumber\\
   &  \sum_{i \in [p]: i \neq j}Y_{i,j}^2 \leq (k-r+1) w_j ( w_j-Y_{j,j}) \ & \forall j \in [p],\nonumber\\
   & (\bm{Y}^t, {\bm{Z}_t}) \in \mathcal{T}(k_t) \ & \forall t \in [r].\nonumber
\end{align}

\begin{remark}
The relaxation of Problem \eqref{prob:disjunctiverelax_permutationinvariant} dominates that of Problem \eqref{prob:relax_ext_str}, even when \eqref{prob:relax_ext_str} is strengthened with the inequalities proposed in \eqref{ell1_split}-\eqref{eqref:soc2}.  
Indeed, \citet[Theorem 13]{kim2021convexification} can be extended to show that \eqref{ell1_split}-\eqref{eqref:soc2} are redundant in \eqref{prob:disjunctiverelax_permutationinvariant}. 
\end{remark}

\subsection{Feasible Solutions via Greedy Disjoint Rounding}\label{sec:rounding}
We develop a rounding mechanism that converts an optimal solution to our semidefinite convex relaxations into a high-quality feasible solution. Historically, a useful strategy for similar integer optimization problems has been to $(a)$ solve a convex relaxation in $(\bm{Z}, \bm{Y})$, $(b)$ greedily round $\bm{Z}^\star$, the solution to the relaxation, to obtain a feasible binary matrix $\hat{\bm{Z}}$ that is close to $\bm{Z}^\star$, and $(c)$ resolve for $\bm{U}$ under the constraints $U_{i,t}=0 \ \text{if} \ \hat{Z}_{i,t}=0$ \cite[cf.][]{bertsimas2020solving}. 

Unlike in the case with a single PC, observe that the ``resolve'' step $(c)$ is non-trivial. Namely,
solving for $\bm{U}$ for some arbitrary and fixed sparsity pattern $\hat{\bm{Z}}$, i.e., solving 
\begin{align}
    \max_{\bm{U}\in \mathbb{R}^{p \times r}}\quad & \langle \bm{U}\bm{U}^\top, \bm{\Sigma}\rangle \ \text{s.t.} \ \bm{U}^\top \bm{U}=\mathbb{I}, \ U_{i,t}=0 \ \text{if} \ \hat{Z}_{i,t}=0, \ \forall i \in [p], t \in [r],\label{prob:resolve}
\end{align}
cannot be done in closed form in general.
{\blue 
Actually, we theoretically characterize the worst-case complexity of Problem~\eqref{prob:resolve}. 
\begin{theorem} \label{thm:resolve.nphard}
Problem \eqref{prob:resolve} is NP-hard. 
\end{theorem}
The proof of Theorem~\ref{thm:resolve.nphard} (deferred to Section~\ref{sec:a.nphard}) relies on a reduction from a variant of the exact cover by 3-sets problem, which is known to be NP-complete \citep[][Chapter A3.1]{garey1979computers}. The hard instances generated by this reduction have $r \geq p/2$, $k_t = 4$, and neither fully overlapping nor fully disjoint supports. Indeed,}
when $\hat{\bm{Z}}$ corresponds to fully overlapping or fully disjoint supports, we can obtain the solution of \eqref{prob:resolve} {\blue in polynomial time} via an eigenvalue decomposition of the corresponding submatrix/submatrices.

Therefore, we propose in Algorithm \ref{alg:greedymethod2} a relax-round-and-resolve strategy where the rounding step $(b)$ generates a solution with fully disjoint supports, i.e., where $\sum_{t \in [r]}\hat{Z}_{i,t}\leq 1, \ \forall i \in [p]$. 

The rounding step searches for the binary disjoint support vector $\bm{Z}$ that is aligned with the solution of the relaxation as much as possible. 
Note that maximizing the alignment between the rounded and the relaxed solution, $\langle \bm{Z}, \bm{Z}^\star \rangle$, is equivalent to minimizing the distance between the two since $\| \bm{Z} - \bm{Z}^\star \|^2 = k - 2 \langle \bm{Z}, \bm{Z}^\star \rangle + \| \bm{Z}^\star \|^2$.

Since $\hat{\bm Z}$ encodes for disjoint supports, any matrix $\bm{U}$ with support $\hat{\bm{Z}}$ automatically satisfies the orthogonality constraint. We can obtain $\bm{U}$ solution of \eqref{prob:resolve} by solving for each PC independently. For each $t \in [r]$, we consider the submatrix of $\bm{\Sigma}$ over the indices $\{i : \hat{Z}_{i,t}=1\}$, extract its leading eigenvector via SVD, and pad it with zeros to construct $\bm{U}_t$.
Although the restriction to disjoint support is not without loss of optimality, we will observe numerically in Section \ref{sec:numres} that disjoint solutions are not {particularly} suboptimal for Problem \eqref{prob:spcaorth_compact} when $k$ and $r$ are small relative to $p$---an observation already made by \cite{asteris2015sparse}. 

\setcounter{algorithm}{0}
\begin{algorithm*}[h]
\caption{A disjoint greedy rounding method of the semidefinite relaxation}
\label{alg:greedymethod2}
\begin{algorithmic}\normalsize
\REQUIRE Covariance matrix $\bm{\Sigma}$, rank parameter $r$, sparsity parameter $k$ or $k_t$
\STATE Compute $UB$ the objective value of \eqref{prob:relax_ext_str} or \eqref{prob:disjunctiverelax_permutationinvariant}
\STATE Compute $\bm{Z}^\star$ solution of \eqref{prob:relax_ext_str} or \eqref{prob:disjunctiverelax_permutationinvariant}
with the constraint $\sum_{t=1}^r Z_{i,t}\leq 1, \ \forall i \in [p]$
\STATE Construct $\hat{\bm{Z}} \in \{0, 1\}^{p \times r}$ solution of 
\begin{align*}
    \max_{\bm{Z} \in \{0, 1\}^{p \times r}} \ \langle \bm{Z}, \bm{Z}^\star\rangle \  
    \text{s.t.}  & \ \sum_{i=1}^p Z_{i,t}\geq 1, \ \forall t \in [r], \\
    & \ \sum_{i \in [p], t \in [r]} Z_{i,t} \leq k \mbox{ or }  \sum_{i \in [p]} Z_{i,t} \leq k_t, \ \forall t \in [r], \\
    & \ \sum_{t=1}^r Z_{i,t}\leq 1, \ \forall i \in [p].
\end{align*}
\STATE Compute $\bm{U}$ solution of \eqref{prob:resolve} via SVD
\RETURN $UB$, $\hat{\bm{Z}}, \bm{U}$.
\end{algorithmic}
\end{algorithm*}

\section{Lagrangian Relaxation Approach}\label{sec:deflation}
In this section, we investigate an approach for obtaining upper and lower bounds based on the theory of {\blue penalty methods and Lagrangian relaxations} \citep[e.g.,][\blue Chapter 17]{fisher1981lagrangian, bertsekas1996constrained,nocedal2006numerical},
which argues that if a non-convex problem is decomposable as a sum of easier (but still non-convex) subproblems with a coupling constraint, a good strategy is often to penalize the coupling constraint in the objective and iteratively solve the non-convex subproblems with different multipliers on the coupling constraint.
For example, \citet{lu2012augmented} proposes an augmented Lagrangian method for solving an $\ell_1$ relaxation of the sparse PCA, where the sparsity constraints are replaced by an $\ell_1$ penalty in the objective.

Unlike the other two approaches we investigate, the Lagrangian only applies to the case where a separate sparsity budget $k_t$ is imposed on each PC in Problem \eqref{prob:spcaorth_compact}.

\subsection{Lagrangian Relaxation and Upper Bound}\label{ssec:deflation.ub}
We consider the case where a separate sparsity budget $k_t$ is imposed on each PC so that the orthogonality constraints are the only coupling constraints. Namely we consider
\begin{equation}\label{prob:spca_extended_lagrangean}
\begin{aligned}  
    \max_{\bm{U}_t \in \mathbb{R}^{p}, t \in [r]} \quad \sum_{t \in [r]} \langle  \bm{U}_t\bm{U}_t^\top, \bm{\Sigma}\rangle \quad 
    \text{s.t.} \quad & \langle \bm{U}_t, \bm{U}_{t'} \rangle = \delta_{t,t'},\ \forall t,t' \in [r], 
    \\ & \| \bm{U}_t \|_0 \leq k_t,\ \forall t \in [r].
\end{aligned}
\end{equation}
We introduce a penalty parameter 
$\lambda_{t,t'} > 0$ for each orthogonality constraint $\langle \bm{U}_t, \bm{U}_{t'} \rangle = 0$ for $t \neq t'$. Observing that 
$(\langle \bm{U}_t, \bm{U}_{t'} \rangle)^2 = \langle \bm{U}_t \bm{U}_t^\top, \bm{U}_{t'}\bm{U}_{t'}^\top \rangle$, we obtain the following {\blue relaxation:}
\begin{equation}\label{prob:spca_extended_lagrangean_2}
\begin{aligned}  
    \max_{\bm{U}_t \in \mathbb{R}^{p}, t \in [r]} \quad \sum_{t \in [r]} \langle  \bm{U}_t\bm{U}_t^\top, \bm{\Sigma}\rangle  - \sum_{t, t' \in [r]: t\neq t'} \lambda_{t,t'} \langle \bm{U}_t \bm{U}_t^\top, \bm{U}_{t'}\bm{U}_{t'}^\top \rangle \quad 
    \text{s.t.} \quad & \langle \bm{U}_t, \bm{U}_{t} \rangle = 1,\ \forall t \in [r], 
    \\ & \| \bm{U}_t \|_0 \leq k_t,\ \forall t \in [r].
\end{aligned}
\end{equation}
{\blue Note that the penalty in the objective of Problem \eqref{prob:spca_extended_lagrangean_2} can be interpreted as a quadratic penalty on the orthogonality constraints $\langle \bm{U}_t, \bm{U}_{t'} \rangle = 0, t \neq t'$ or, equivalently, as the Lagrangian of \eqref{prob:spca_extended_lagrangean} where the orthogonality constraints are reformulated as $\langle \bm{U}_t, \bm{U}_{t'} \rangle^2 = 0, t \neq t'$.}
For a fixed value of the 
parameters $\bm{\lambda}$, Problem \eqref{prob:spca_extended_lagrangean_2} provides an upper bound on the objective value of \eqref{prob:spca_extended_lagrangean} by the Lagrangian duality theorem. Optimizing for all PCs $\bm{U}_t$ simultaneously in Problem \eqref{prob:spca_extended_lagrangean_2} is challenging due to the non-convex objective. For a given index $t$, however, optimizing for $\bm{U}_t$ (with all other $\bm{U}_{t'}$, $t' \neq t$, and $\bm{\lambda}$ fixed) is equivalent to finding the leading $k_t$-sparse PC of the matrix $\bm{\Sigma}-\sum_{t' \neq t} \lambda_{t,t'} \bm{U}_{t'}\bm{U}_{t'}^\top$, for which efficient algorithms have been developed recently (see Section \ref{sec:litrev})\endnote{Most sparse PCA algorithms for $r=1$ require the input covariance matrix to be semidefinite. Accordingly, since $\mathrm{tr}(\bm{U}_t \bm{U}_t^{\top})=1$, we can add a constant term $\lambda_{\text{offset}} \mathrm{tr}(\bm{U}_t \bm{U}_t^{\top})$ to the objective without impacting the optimal solution. {\color{black}Thus}, we can pick $\lambda_{\text{offset}} > 0$ to be sufficiently large that the matrix $\bm{\Sigma}-\sum_{t' \neq t} \lambda_{t,t'} \bm{U}_{t'}\bm{U}_{t'}^\top + \lambda_{\text{offset}} \mathbb{I}$ and use these algorithms off-the-shelf.}. 
We use this recent technology to derive an upper bound on \eqref{prob:spca_extended_lagrangean}. Namely, for $\lambda = 0$, we solve \eqref{prob:spca_extended_lagrangean_2} by solving $r$ single-PC sparse PCA problems, or equivalently
\begin{align}  
    \max_{\bm{U} \in \mathbb{R}^{p\times r}} \quad \sum_{t \in [r]}\bm{U}_t^\top \bm{\Sigma} \bm{U}_t \quad 
    \text{s.t.} \quad \| \bm{U}_t \|_2 = 1,\ \| \bm{U}_t \|_0 \leq k_t, \forall t \in [r]\label{eqn:lagrangianbound}
\end{align}
using a certifiably optimal method \citep[][in our implementation]{berk2019certifiably}. Note that $\bm{U}_t$ denotes the $t$th column of $\bm{U}$ and the overall objective value provides an upper bound on \eqref{prob:spca_extended_lagrangean}. Actually, \eqref{eqn:lagrangianbound} provides an $r$-factor approximation of \eqref{prob:spca_extended_lagrangean}, as formally stated below.
\begin{proposition}\label{prop:lagrangean.factor} The objective value of \eqref{eqn:lagrangianbound} is at most $r$ times the objective value of the original sparse PCA problem \eqref{prob:spca_extended_lagrangean}.
\end{proposition}
\proof{Proof of Proposition \ref{prop:lagrangean.factor}} We already {\color{black}k}now that \eqref{eqn:lagrangianbound} provides an upper bound on \eqref{prob:spca_extended_lagrangean}, $\eqref{prob:spca_extended_lagrangean} \leq \eqref{eqn:lagrangianbound}$ in short. Furthermore, \eqref{eqn:lagrangianbound} is bounded above by
\begin{align*}  
\quad\quad\quad\quad\quad    
 r \: \times \: \max_{t \in [r]} \left\lbrace \max_{\bm{u} \in \mathbb{R}^{p}} \quad \bm{u}^\top \bm{\Sigma} \bm{u} \quad 
    \text{s.t.} \quad \| \bm{u} \|_2 = 1,\ \| \bm{u} \|_0 \leq k_t \right\rbrace \quad \leq r \: \times \: \eqref{prob:spca_extended_lagrangean}. 
\quad\quad\quad\quad\quad \Halmos
\end{align*}
\endproof

\subsection{A Heuristic Based on Iterative Deflation}\label{ssec:deflation.heuristic}
When solving \eqref{prob:spca_extended_lagrangean_2} for a fixed $\bm{\lambda}$, the returned PCs may not be feasible for the original Problem \eqref{prob:spcaorth_compact} because they can violate the orthogonality constraints. A strategy to obtain a feasible solution is thus to solve a sequence of problems of the form \eqref{prob:spca_extended_lagrangean_2} with an increasing penalty parameter value, as we propose in Algorithm \ref{alg:alternatingmin}.

At each iteration, Algorithm \ref{alg:alternatingmin} optimizes for each $\bm{U}_t$ in \eqref{prob:spca_extended_lagrangean_2} sequentially (instead of simultaneously).
{\blue Consequently, Algorithm \ref{alg:alternatingmin} is not guaranteed to converge to a near-optimal solution. 
However, when the updates of each column $\bm{U}_t$ are $\epsilon$-optimal \citep[e.g., when using the methods of][]{berk2019certifiably,bertsimas2020solving}, we show in Proposition \ref{prop.ec.convsubsequence} that the limiting solutions are feasible. 
For tractability considerations, we use the truncated power method of \citet{yuan2013truncated} (instead of exact methods) for which we cannot guarantee $\epsilon$-optimality. However, we can still show convergence to some stationary solutions (see Proposition~\ref{prop:penalty.guarantee}). Unfortunately, however, we cannot rule out convergence to infeasible solutions, which is a known limitation of the quadratic penalty method \citep[see][Theorem 17.2 and discussion therein]{nocedal2006numerical}. In practice, we guarantee $\eta$-feasibility of the returned solution by recording all $\eta$-feasible solutions generated and returning the best one.}


\setcounter{algorithm}{1}
\begin{algorithm*}[h]
\caption{Lagrangian Alternating Maximization for Problem \eqref{prob:spcaorth_compact}}
\label{alg:alternatingmin}
\begin{algorithmic}\normalsize
\REQUIRE Matrix $\bm{\Sigma}$, rank parameter $r$, sparsity parameters $k_1, \ldots, k_r$, number of iterations $L$, feasibility tolerance $\eta$
\STATE Compute $UB$ the optimal value of \eqref{eqn:lagrangianbound}
\STATE Initialize $\bm{U}_t^{\blue (0)} \leftarrow \bm{0}$
\FOR{$\ell = 1,\dots,L$}
\FOR{$t=1,\dots,r$}
\STATE Compute $\bm{U}_t^{\blue (\ell)}$ an approximate solution of 
\begin{align*} 
    \max_{\bm{u} \in \mathbb{R}^{p}} \quad \left\langle  \bm{u}\bm{u}^\top, \bm{\Sigma} - \sum_{t' \in [r]: t' \neq t} \lambda_{t,t'} \bm{U}_{t'}^{\blue (\ell-1)}\bm{U}_{t'}^{\blue (\ell-1) \top} \right\rangle \ 
    \text{s.t.} \ \| \bm{u} \|_2 = 1,\ \| \bm{u} \|_0 \leq k_t.
\end{align*}
\ENDFOR
\IF{$\sum_{t,t'} | \langle \bm{U}_t^{\blue (\ell)}, \bm{U}_{t'}^{\blue (\ell)} \rangle - \delta_{t,t'} | \leq \eta$}
\STATE Consider $\left\{ \bm{U}_t^{\blue (\ell)}, t \in [r] \right\}$ feasible and record its objective value
\ENDIF
\STATE Update $\bm{\lambda}$ value.
\ENDFOR
\RETURN $UB$ and the best feasible solution $\{ \bm{U}_t, t \in [r] \}$ found
\end{algorithmic}
\end{algorithm*}

{\blue In our block computations (column-wise) of $\bm{U}^{(\ell)}$, we use the previous values of $\bm{U}^{(\ell-1)}_{t'}$, even for the columns $t' < t$, which have already been updated (this is the so-called block-Jacobi update rule). Instead, for memory efficiency, one could use $\bm{U}^{(\ell)}_{t'}$ for $t' < t$ and $\bm{U}^{(\ell-1)}_{t'}$ for $t' > t$ (known as Gauss-Seidel updates; see Remark~\ref{rem:jacobi}).}
In our implementation, we consider a feasibility tolerance (e.g., $10^{-4}$) and return the best solution among those satisfying the orthogonality constraints up to this tolerance. Note that our definition of constraint violation, $\sum_{t,t'} | \langle \bm{U}_t, \bm{U}_{t'} \rangle - \delta_{t,t'} |$, includes the constraints $\| \bm{U}_t \|_2 = 1$, so that the initialization (all-zero solution) is correctly considered as infeasible. 
{\color{black}For the termination criterion, we set a limit on the total number of iterations.} We can also terminate the algorithm when it stalls, i.e., when the difference in objective value between two consecutive feasible solutions fall{\color{black}s} below a certain threshold $\epsilon$, or impose a time limit. 
Regarding the update of the penalty parameter $\bm{\lambda}$, we increase it progressively as is standard in the Lagrangian relaxation literature \citep[][Chapter 17]{fisher1981lagrangian,bertsekas1996constrained,nocedal2006numerical}. 
Typically, one can take increments proportional to the gradient of the Lagrangian with respect to $\bm{\lambda}$, $\langle \bm{U}_t, \bm{U}_{t'} \rangle^2$, or as a scaling factor between the two concurrent objectives, $\sum_{t \in [r]} \langle \bm{U}_t \bm{U}_t^\top, \bm{\Sigma}\rangle \left/ \sum_{t, t' \in [r]: t\neq t'}\langle \bm{U}_t, \bm{U}_{t'} \rangle^2 \right.$. In our implementation, we use the former increment (i.e., gradient-based) during the first iterations, and the latter (which is typically larger) at later stages, the intuition being that we slowly increase $\bm{\lambda}$ in the beginning to enable exploration and switch to being more aggressive at later iterations to ensure that we are generating feasible solutions. Precisely, we use parameters of the form $\lambda_{t,t'} = w_t \Tilde{\lambda}$, where $w_t$ is equal to the objective value obtained by the PC $\bm{U}_t$ found at the first iteration ($\ell = 1$), and update $\Tilde{\lambda}$ according to the following rule
\begin{align*}
    \mbox{initialization: } & \quad  \Tilde{\lambda} = 0, \\
    \mbox{iteration } \ell = 1, \dots, \lceil 0.15 L \rceil - 1 & \quad \Tilde{\lambda} \leftarrow \Tilde{\lambda} + \alpha \times \sum_{t \neq t'} \langle \bm{U}_t, \bm{U}_{t'} \rangle^2, \\
    \mbox{iteration } \ell = \lceil 0.15 L \rceil, \dots, L & \quad \Tilde{\lambda} \leftarrow \Tilde{\lambda} + \alpha \times \sum_{t \in [r]} \langle \bm{U}_t \bm{U}_t^\top, \bm{\Sigma}\rangle \left/ \sum_{t, t' \in [r]: t\neq t'}\langle \bm{U}_t, \bm{U}_{t'} \rangle^2 \right. ,
\end{align*}
with $\alpha$ a fixed stepsize parameter.

\section{Combinatorial Approach}
\label{sec.gen.gershgorin}
The Gershgorin circle theorem \citep[cf.][Chapter 6]{johnson1985matrix} bounds the largest eigenvalue of a matrix $\bm{\Sigma} \in \mathcal{S}^p_+$ via the combinatorial function: 
\begin{align} \label{eqn:gershgorin}
    \lambda_{\max}(\bm{\Sigma}) \leq \max_{i \in [p]} \sum_{j \in [p]} \vert \Sigma_{i,j}\vert.
\end{align}
For sparse PCA with a single PC, several authors \citep{berk2019certifiably,bertsimas2020solving} have leveraged this result to derive an upper bound on \eqref{prob:spcarank1} that depends on the support of the sparse PC. 
Motivated by their observations, we now bound the objective value of \eqref{prob:spcaorth_compact} as a function of the support matrix $\bm{Z}$, and leverage this bound to generate upper bounds and feasible solution{\color{black}s} to Problem \eqref{prob:spcaorth_compact}.

\subsection{Generalizing Gershgorin Circle Theorem}\label{ssec:gengershgorin}
One naive upper bound is to apply the {\color{black}Gershgorin} circle theorem to each PC separately and bound the objective of \eqref{prob:spcaorth_compact} via the sum of the largest eigenvalues of each sub-matrix of $\bm{\Sigma}$ induced by a column of $\bm{Z}$. However, this approach is too conservative for non-disjoint PCs, since it does not take into account any information about {\color{black}the} overlap between the support of each PC. Indeed, with fully overlapping support, this approach bounds the sum of the $r$ largest eigenvalues of the relevant submatrix by $r$ times {\color{black}its} largest eigenvalue. Instead, the valid inequalities we derive in this section rely on a new and non-trivial bound on the variance collectively explained by $r$ orthogonal PCs. Incidentally, our result leads to the following bound o{\color{black}n} the sum of the $r$ largest eigenvalues of a semidefinite matrix $\bm{\Sigma}$:
\begin{align} \label{eqn:bound.topr}
    \sum_{t \in [r]} \lambda_t (\bm{\Sigma}) \leq \max_{\bm{\mu} \in \{0, 1\}^p: \bm{e}^\top \bm{\mu} 
    \leq r}\sum_{i,j \in [p]} \mu_i \vert \Sigma_{i,j}\vert,
\end{align}
which strictly generalizes \eqref{eqn:gershgorin} and could be of independent interest.\endnote{As suggested by an anonymous reviewer, \eqref{eqn:bound.topr} could be derived from linear algebra principles directly, instead of as a corollary of Theorem \ref{thm.MILPbound_multiPCs}. Denoting $r_i := \sum_{j \in [p]} |\Sigma_{i,j}|$, the matrix $\operatorname{Diag}(\bm{r}) - \bm{\Sigma}$ is diagonally dominant, hence {\color{black}positive semidefinite}. Therefore, we have $\operatorname{Diag}(\bm{r}) \succeq \bm{\Sigma}$. 
In particular, the ordered vector of eigenvalues of $\operatorname{Diag}(\bm{r})$ majorizes that of $\bm{\Sigma}$ \citep[section L.1 of][]{marshall1979inequalities}, which leads precisely to \eqref{eqn:bound.topr}. 
}

Formally, we derive the following generalization of the circle theorem bound, which holds with multiple PCs and a fixed but arbitrary support pattern $\bm{Z} \in \{0, 1\}^{p \times r}$ {satisfying} $ \sum_{i \in [p]}Z_{i,t} \geq 1 \ \forall t \in [r]$. Subsequently, we derive a mixed-integer linear representation of this bound:
\begin{theorem}\label{thm.MILPbound_multiPCs}
For any feasible support pattern $\bm{Z} \in \{0,1\}^{p \times r}: \sum_{i \in [p]}Z_{i,t} \geq 1 \ \forall t \in [r]$, an upper bound on the objective value attained by any matrix $\bm{U}$ such that $U_{i,t}=0 \ \text{if} \ Z_{i,t}=0 \ \forall i \in [p], \ \forall t \in [r]$ in Problem \eqref{prob:spcaorth_compact} is given by:
\begin{align}\label{prob.linearalgebraicbound}
    \max_{\substack{\bm{\mu} \in \{0, 1\}^{p\times r}: \\
    \sum_{i \in [p]}\mu_{i,t}= 1, \forall t \in [r], \\
    \sum_{t \in [r]}\mu_{i,t}\leq 1 \forall i \in [p]}} \quad & 
    \sum_{i,j \in [p]}\sum_{t \in [r]} \mu_{i,t} Z_{i,t} Z_{j,t} \vert \Sigma_{i,j}\vert.
\end{align}
\end{theorem}

\begin{remark}
    Taking $\bm{Z}=\bm{E}$ in Theorem \ref{thm.MILPbound_multiPCs} yields \eqref{eqn:bound.topr}.
If $r=1$, this bound is equivalent to the circle theorem, and if $r=p$ it is equivalent to the (known) fact that $\sum_{t \in [p]}\lambda_t(\bm{\Sigma})=\mathrm{tr}(\bm{\Sigma}) \leq \sum_{i,j \in [p]}\vert \Sigma_{i,j}\vert$. More generally, it shows that the eigenvalues of a positive semidefinite matrix are majorized by the absolute column sums \citep[see][for a general theory of majorization]{marshall1979inequalities}.
\end{remark}
\begin{remark}
If $\bm{\Sigma}$ is a diagonally dominant matrix, i.e., $\Sigma_{i,i}\geq \sum_{j \neq i} \vert \Sigma_{i,j}\vert \ \forall i$, then 
{\blue for any $(i,t)$ with $Z_{i,t}=1$, $\sum_{j \in [p]} Z_{j,t} \vert \Sigma_{i,j}\vert = \Sigma_{i,i}   + \sum_{j \neq i}Z_{j,t} \vert \Sigma_{i,j}\vert \leq 2 \Sigma_{i,i}$. Hence, \eqref{prob.linearalgebraicbound} is lower than the maximum of $2 \times \sum_{t} \Sigma_{i_t,i_t}$ over distinct indices $i_t$ such that $Z_{i_t,t} =1$. Observe that $\sum_{t} \Sigma_{i_t,i_t} = \langle \bm{U}\bm{U}^\top, \bm{\Sigma}\rangle$ for the matrix $\bm{U}$ whose columns are the canonical basis vectors $\bm{e}_{i_t}$. 
Consequently,} \eqref{prob.linearalgebraicbound} provides a $2$-factor approximation {\blue of the optimal} objective value \blue in this case.     
\end{remark}
\begin{remark}
    \blue For a fixed $\bm Z$, Problem \eqref{prob.linearalgebraicbound}  is a matching problem, where each PC $t \in [r]$ has to be matched to exactly one coordinate $i \in [p]$, and each coordinate can be matched to at most one PC. The Hopcroft-Karp algorithm 
    \citep{hopcroft1973n} solves it in $O(rp\sqrt{r+p})$ time.
\end{remark}
When the supports of each PC are disjoint, \eqref{prob.linearalgebraicbound}'s upper bound is equivalent to applying the circle theorem to each PC separately. Alternatively, with fully overlapping support, it reduces to bounding the variance explained by the largest $r$ column sums of the submatrix selected by $\bm{Z}$. In the case with partially overlapping support, it systematically interpolates between these bounds.

\proof{Proof of Theorem \ref{thm.MILPbound_multiPCs}}
Fix $\bm{Z} \in \{0, 1\}^{p \times r}$ in Problem \eqref{prob:spcaorth_compact}. Then, it follows directly from Theorem \ref{thm:reform} that an upper bound on the objective value attained by any orthogonal matrix $\bm{U}$ such that $U_{i,t}=0 \ \text{if} \ Z_{i,t}=0$ is given by the following maximization problem:
\begin{align*}
    \max_{\bm{Y}^t \in \mathcal{S}^p_+ 
     \ \forall t \in [r]
    } \ \sum_{t 
    \in [r]} \langle \bm{Y}^t, \bm{\Sigma} \rangle \
    \text{s.t.} \ \mathrm{tr}(\bm{Y}^t)=1 \ \forall t \in [r], \sum_{t \in [r]} \bm{Y}^t \preceq \mathbb{I}, Y_{i,j}^t=0 \ \text{if} \ Z_{i,t}=0 \ \forall i \in [p], t \in [r].
\end{align*}
To obtain a non-trivial mixed-integer linear representable upper bound as a function of the support pattern $\bm{Z}$, we now relax this problem. First, we observe that $Y_{i,j}^t=0$ if $Z_{i,t}=0$ and therefore we can replace $\bm{\Sigma}$ in the objective with $\mathrm{Diag}({\bm{Z}_t})\bm{\Sigma}\mathrm{Diag}({\bm{Z}_t})$ without loss of generality, where $\mathrm{Diag}({\bm{Z}_t})$ is a diagonal matrix with on-diagonal entries specified by the $t$th column of $\bm{Z}$. Further relaxing the problem by omitting the logical constraints then gives the following semidefinite upper bound:
\begin{align*}
    \max_{\bm{Y}^t \in \mathcal{S}^p_+ 
     \ \forall t \in [r]
    } \ \sum_{t 
    \in [r]} \langle \bm{Y}^t, \mathrm{Diag}({\bm{Z}_t})\bm{\Sigma}\mathrm{Diag}({\bm{Z}_t}) \rangle \
    \text{s.t.} \ \mathrm{tr}(\bm{Y}^t)=1 \ \forall t \in [r], \sum_{t \in [r]} \bm{Y}^t \preceq \mathbb{I}.
\end{align*}

Moreover, strong duality holds between this problem and its dual problem, namely:
\begin{align*}
    \min_{\bm{U} \in \mathcal{S}^p_+, \bm{s} \in \mathbb{R}^r} \quad
    \mathrm{tr}(\bm{U})+\bm{e}^\top \bm{s}\quad \text{s.t.} \quad \bm{U}+s_t \mathbb{I}\succeq \mathrm{Diag}({\bm{Z}_t})\bm{\Sigma}\mathrm{Diag}({\bm{Z}_t}) \ \forall t \in [r].
\end{align*}

To obtain a linear, rather than semidefinite, upper bound from this problem, we restrict $\bm{U}$ to be a diagonal matrix and $\bm{U}+s_t \mathbb{I}- \mathrm{Diag}({\bm{Z}_t})\bm{\Sigma}\mathrm{Diag}({\bm{Z}_t})$ to be contained within the cone of diagonally dominant matrices, which is an inner approximation of the positive semidefinite cone \citep[see also][for detailed studies of this inner approximation]{barker1975cones, ahmadi2017optimization}. This gives the following upper bound:
\begin{align*}
    \min_{\bm{u} \in \mathbb{R}^p_+, s} \quad & 
    \bm{e}^\top \bm{u}+\bm{e}^\top \bm{s}\ \text{s.t.} \ u_i+s_t \geq \sum_{j \in [p]}Z_{i,t} Z_{j,t} \vert \Sigma_{i,j}\vert \ \forall i \in [p], \forall t \in [r].
\end{align*}

Finally, we invoke strong duality and use the fact that some (binary) extreme point in the dual problem must be dual-optimal, to verify that the above problem attains the same value as:
\begin{align*}
    \max_{\substack{\bm{\mu} \in \{0, 1\}^{p\times r}: \\
    \sum_{i \in [p]}\mu_{i,t}= 1 \ \forall t \in [r], \\ 
    \sum_{t \in [r]}\mu_{i,t}\leq 1 \ \forall i \in [p]}} \quad & 
    \sum_{i,j \in [p]}\sum_{t \in [r]} \mu_{i,t} Z_{i,t} Z_{j,t} \vert \Sigma_{i,j}\vert. \quad \quad \Halmos
\end{align*}
\endproof

Observe that the proof of Theorem \ref{thm.MILPbound_multiPCs} involves invoking strong duality and taking a finitely generated inner approximation of the positive semidefinite cone. This is quite different {\color{black}from} existing proofs of the Gershgorin circle theorem, which usually leverage properties of eigenvectors and therefore cannot easily be generalized.\endnote{Actually, one could derive an alternative proof of Theorem \ref{thm.MILPbound_multiPCs} that relies on linear algebra principles, in much the same way as the alternative proof of \eqref{prob.linearalgebraicbound}. We detail this alternative proof technique in Section \ref{sec:ec.combinatorial.linalgproof}.
}
Thus, our proof technique could also be useful in other contexts, e.g., in sparse canonical correlation analysis \citep{witten2009penalized}.

\subsection{A Mixed-Integer Linear Relaxation}\label{ssec.gershgorin.ub}
Since Theorem \ref{thm.MILPbound_multiPCs} provides an upper bound on the objective value (i.e., the fraction of variance explained) for a given support $\bm{Z}$, we can obtain an upper bound on \eqref{prob:spca_extended} by optimizing \eqref{prob.linearalgebraicbound} over all possible supports, $\bm{Z} \{ 0,1 \}^{p \times r}$: 
\begin{align} \label{eqn:relaxation.linalg}
    \max_{\bm{Z} \{ 0,1 \}^{p \times r}} \: \eqref{prob.linearalgebraicbound} \quad \mbox{s.t.} \quad \sum_{i,t} Z_{i,t} \leq k \mbox{ or } \sum_{i} Z_{i,t} \leq k_t, \forall t \in [r],
\end{align}
depending on the nature of the sparsity budget (global or per-component). Let us observe that Problem \eqref{eqn:relaxation.linalg} can be reformulated as a mixed-integer linear optimization problem by introducing auxiliary variables $\rho_{i,t} \ \forall i \in [p], t \in [r]$ to model the column sum $\sum_{j \in [p]}Z_{j,t} \vert \Sigma_{i,j}\vert$ if $\mu_{i,t}=1$ and equal $0$ if $\mu_{i,t}=0$ in \eqref{prob.linearalgebraicbound}. This allows us to represent Theorem \ref{thm.MILPbound_multiPCs}'s upper bound via the system:
\begin{equation}\label{prob.gen.gershgorin.formulation}
\begin{aligned}
    \quad & \theta=\sum_{i, t} \rho_{i,t}, & \\
    & \rho_{i,t} = \begin{cases} \sum_{j \in [p]} Z_{j,t}\vert \Sigma_{i,j}\vert \ & \text{ if} \ \mu_{i,t}=1 \\ 
    0 \ & \text{ if} \ \mu_{i,t}=0 \end{cases} & \forall i \in [p], t \in [r], \\ 
    & \sum_{i \in [p]} \mu_{i,t} = 1 \ & \forall t \in [r], \\
    & \sum_{t \in [r]} \mu_{i,t} \leq 1 \ & \forall i \in [p],\\
    & \mu_{i,t}\leq Z_{i,t} \ & \forall i \in [p], t \in [r], \\
    & \mu_{i,t} \in \{ 0,1 \} & \forall i \in [p], t \in [r],\\[1em] 
\end{aligned}
\end{equation}
in the variables $(\theta, \bm{\mu},\bm{\rho})$.
We omit the term $Z_{i,t}$ from the partial sum $\sum_{j \in [p]} Z_{j,t}\vert \Sigma_{i,j}\vert$ by imposing the constraint $\mu_{i,t} \leq Z_{i,t}$ because if $Z_{i,t}=0$ then the $(i,t)$th column sum is zero and can be omitted from the bound without loss of generality. 

\subsection{Feasible Solution with Disjoint Supports}\label{ssec.gershgorin.rounding}
Unfortunately, as explained in Section \ref{sec:rounding}, the solution to the mixed-integer linear relaxation \eqref{eqn:relaxation.linalg}, $\bm{Z}^\star \in \{0,1\}^{p \times r}$, cannot be used to generate a feasible solution $\bm{U}$ directly, because estimating the optimal $\bm{U}$ for a given $\bm{Z}$ in non-trivial. To alleviate this issue, as in Algorithm \ref{alg:greedymethod2}, we restrict our attention to feasible solutions with disjoint supports. 

In Algorithm \ref{alg:disjoint.linalg}, we first solve Problem \eqref{eqn:relaxation.linalg}, with the additional constraint that $\bm{Z}$ encodes for disjoint supports, and estimate the optimal PCs with sparsity pattern $\hat{\bm{Z}}$. This gives both a feasible solution $\bm{U}$ and a high-quality upper bound, after solving two MIOs and an eigenproblem.

\setcounter{algorithm}{2}
\begin{algorithm*}[h!]
\caption{A disjoint-support solution from the combinatorial upper bound  rounding \eqref{prob.linearalgebraicbound}}
\label{alg:disjoint.linalg}
\begin{algorithmic}\normalsize
\REQUIRE Covariance matrix $\bm{\Sigma}$, rank parameter $r$, sparsity parameter $k$
\STATE Compute $UB$ the optimal value of \eqref{eqn:relaxation.linalg}
\STATE Construct $\hat{\bm{Z}} \in \{0, 1\}^{p \times r}$ solution of 
\begin{align*}
    \max_{\bm{Z} \in \{0, 1\}^{p \times r}} \ \theta \  \text{s.t.}  & \ \sum_{i,t} Z_{i,t} \leq k \mbox{ or } \sum_{i} Z_{i,t} \leq k_t, \forall t \in [r], \\
    & \ \sum_{t=1}^r Z_{i,t}\leq 1, \ \forall i \in [p], \\ 
    & \theta \mbox{ satisfying } \eqref{prob.gen.gershgorin.formulation}.
\end{align*}
\STATE Compute $\bm{U}$ solution of \eqref{prob:resolve} via SVD
\RETURN $UB$, $\hat{\bm{Z}}, \bm{U}$.
\end{algorithmic}
\end{algorithm*}

\section{Numerical Results}\label{sec:numres}
In this section, we evaluate the algorithmic strategies derived in the previous three sections {\blue for solving Problem \eqref{prob:spcaorth_compact}}, implemented in \verb|Julia| $1.9$ using \verb|JuMP.jl| $1.12.0$, \verb|Gurobi| version $10.0.0$ to solve all non-convex quadratically constrained problems, and \verb|Mosek| $10.1.11$ to solve all conic relaxations. For the sake of conciseness, we defer full details of our experimental setup to Section \ref{sec:ec.supp_expr_setup}. Moreover, for the purpose of averaging results across datasets with different $p$'s, we report the proportion of variance explained whenever we report an objective value. For a correlation matrix, this corresponds to dividing by $p$, the number of features. We make our code available on \verb|GitHub| at \href{https://github.com/ryancorywright/MultipleComponentsSoftware}{github.com/ryancorywright/MultipleComponentsSoftware}.

\paragraph{Description of Data Sources:} We perform experiments on eleven datasets from the UCI database in Sections \ref{ssec:bounds}-\ref{ssec:feasiblemethods} and \ref{ssec:symm}, and experiments on synthetic data in Section \ref{ssec:auc}. Of the eleven datasets (described in detail in Section \ref{sec:ec.datadescr}), six datasets are overdetermined (meaning $n >p$), while five datasets are underdetermined (meaning $p>n$). 
Because some of our algorithms {\blue can only} return PC{\color{black}s} with disjoint supports, we sometimes report results over instances where $\sum_t k_t > p$ and where $\sum_t k_t \leq p$ separately. Given the values of $(k_t,r)$ we consider, datasets with $p \geq 60$ are systematically in the second category.

\subsection{Performance of Upper Bounds}\label{ssec:bounds}
In this section, we compare {the tightness and scalability of} the upper bounds {obtained from our three approaches. First, {\color{black}we consider} two semidefinite relaxations} 
(namely, \eqref{prob:relax_ext_str} with \eqref{ell1_split}-\eqref{eqref:soc2}, hereafter ``Extended-Ineq'', and \eqref{prob:disjunctiverelax_permutationinvariant}
, hereafter ``Perm-Ineq''). {This allows us to determine the merits of specifying $k_t$ versus only specifying $k$, in terms of whether it leads to tighter relaxations.} {\color{black}Second}, the Lagrangian upper bound obtained from solving  \eqref{eqn:lagrangianbound} with the exact method of \cite{berk2019certifiably}, hereafter ``Lagrangian''. Finally, the bound attained by maximizing our combinatorial bound \eqref{prob.linearalgebraicbound} in Theorem \ref{thm.MILPbound_multiPCs}, hereafter ``Combinatorial''. We also considered running \verb|Gurobi|'s non-convex branch-and-bound solver directly on this dataset (see Section \ref{sec:ec.branchandbound}). 
However, experiments on the same dataset show that its upper bound does not scale as well as any of our relaxations when the total sparsity of the PCs is $12$ or more (see Tables \ref{tab:comparison_gurobi_pitprops1}-\ref{tab:comparison_gurobi_pitprops2}). 
Thus, we do not report \verb|Gurobi|'s performance in terms of generating upper bounds in the main paper.

\paragraph{Benchmarking on Pitprops Data} We first compare the bounds generated by each method---in terms of proportion of correlation explained---on the \verb|pitprops| dataset ({$n=180$,} $p=13$) as we vary $r \in \{2,3\}$ and $k$, in Table \ref{tab:comparison_bounds_new}. We consider both imposing an overall sparsity budget alone (denoted by ``$k, -$'') and imposing a separate budget for each PC, denoted by ``$k, (k_1, \ldots, k_t)$''. We denote the performance of methods that require {\color{black}$\{k_t\}_{t \in [r]}$} by a dash when only $k$ is provided, to indicate they are not applicable. 

\begin{table}
\centering\footnotesize
\begin{tabular}{@{}l l r r r r r r r r r r r r r r @{}} \toprule
Rank ($r$) &  Sparsity ($k$, $k_t$) & \multicolumn{2}{c@{\hspace{0mm}}}{Extended-Ineq} & \multicolumn{2}{c@{\hspace{0mm}}}{Perm-Ineq} & \multicolumn{2}{c@{\hspace{0mm}}}{Lagrangian} & \multicolumn{3}{c@{\hspace{0mm}}}{Combinatorial}\\
\cmidrule(l){3-4} \cmidrule(l){5-6} \cmidrule(l){7-8} \cmidrule(l){9-11}   &  & UB & T(s) & UB & T(s) & UB & T(s) & UB & Nodes & T(s)\\\midrule
2 &	4, - &	0.297 &	20.28 &	- &	- &	- &	- & 0.301 & 248 & 0.048 \\
2 &	4, (1, 3) &	\textbf{0.267} &	20.58 &	\textbf{0.267} &	1.97 &	\textbf{0.267} & 0.22 & 0.277 &  1 & 0.007\\
2 &	4, (2, 2) &	\textbf{0.295} &	0.37 &	\textbf{0.295} &	0.47 &	0.301 & 0.03  & 0.301 & 1 & 0.007\\\midrule
2 &	6, &	\textbf{0.384} &	0.78 &	- &	- &	- &	- & 0.396 &  195 & 0.046\\
2 &	6, (1, 5) &	\textbf{0.339} &	0.44 &	\textbf{0.339} &	0.25 &	\textbf{0.339} & 0.02  & 0.360 &  1 & 0.006 \\
2 &	6, (2, 4) &	\textbf{0.371} &	0.47 &	\textbf{0.371} &	0.57 &	0.376 & 0.02  & 0.394 &  1 & 0.005\\
2 &	6, (3, 3) &	0.361 &	0.42 &	\textbf{0.360} &	0.45 &	0.381 & 0.03  & 0.396 & 1 & 0.009 \\\midrule
2 &	8, - &	\textbf{0.451} &	0.75 &	- &	- &	- &	- & 0.482 & 182 & 0.046\\
2 &	8, (1, 7) &	\textbf{0.384} &	0.46 &	\textbf{0.384} &	0.40 &	\textbf{0.384} & 0.02 & 0.415 &  1 & 0.006\\
2 &	8, (2, 6) &	\textbf{0.435} &	0.43 &	\textbf{0.435} &	0.48 &	0.440 & 0.03   & 0.465 &  1 & 0.005\\
2 &	8, (3, 5) &	0.420 &	0.52 &	\textbf{0.418} &	0.65 &	0.452 & 0.02  & 0.478  & 1 & 0.005\\
2 &	8, (4, 4) &	0.412 &	0.48 &	\textbf{0.408} &	0.54 &	0.452 & 0.02  & 0.482 & 1 & 0.007\\\midrule
2 &	10, - &	\textbf{0.490} &	0.72 &	- &	- &	- &	- & 0.559 & 213 & 0.045\\
2 &	10, (1, 9) & \textbf{0.395} &	0.43 &	\textbf{0.395} &	0.3 &	\textbf{0.395} & 0.02  & 0.465 & 1 & 0.006\\
2 &	10, (2, 8) &	\textbf{0.457} &	0.53 &	\textbf{0.457} &	0.5 &	0.463 & 0.02  & 0.516 & 1 & 0.007\\
2 &	10, (3, 7) &	0.461 &	0.41 &	\textbf{0.459} &	0.5 &	0.498 & 0.03  & 0.537 & 1 & 0.005\\
2 &	10, (4, 6) &	0.458 &	0.44 &	\textbf{0.455} &	0.6 & 0.516 & 0.02  & 0.553 & 1 & 0.005\\
2 &	10, (5, 5) &	0.453 &	0.55 &	\textbf{0.449} &	0.45 &	0.524 & 0.02  & 0.559 & 1 & 0.008\\\midrule
3 &	6, - &	\textbf{0.443} &	25.76 &	- &	- &	- & - & 0.445 & 1872 & 0.147\\
3 &	6, (1, 1, 4) &	\textbf{0.380} &	42.92 &	\textbf{0.380} &	5.02 &	\textbf{0.380} & 0.05  & 0.398 & 1 & 0.008\\
3 &	6, (1, 2, 3) &	\textbf{0.412} &	0.94 &	\textbf{0.412} &	1.54 &	{0.418} & 0.04  & 0.427 & 1 & 0.007\\
3 &	6, (2, 2, 2) &	\textbf{0.435} &	0.68 &	\textbf{0.435} &	2.65 &	0.451 & 0.03  & 0.445 & 1 & 0.009\\\midrule
3 &	9, - &	\textbf{0.570} &	2.10 &	- &	- &	- &	- & 0.588 & 1759 & 0.118\\
3 &	9, (1, 1, 7) &	\textbf{0.461} &	0.80 &	\textbf{0.461} &	1.21 &	\textbf{0.461} & 0.03  & 0.492 & 1 & 0.008\\
3 &	9, (1, 2, 6) &	\textbf{0.512} &	0.84 &	\textbf{0.512} &	0.59 &	0.517 & 0.03  & 0.542 & 1 & 0.007\\
3 &	9, (1, 3, 5) &	0.497 &	0.71 &	\textbf{0.495} &	0.80 &	0.529 & 0.04  & 0.555 & 1 & 0.008\\
3 &	9, (1, 4, 4) &	0.489 &	0.65 &	\textbf{0.485} &	0.54 &	0.529 & 0.03  & 0.559 & 1 & 0.010\\
3 &	9, (2, 2, 5) &	\textbf{0.539} &	0.78 &	\textbf{0.539} &	0.73 &	0.563 & 0.03  & 0.578 & 1 & 0.010\\
3 &	9, (2, 3, 4) &	0.532 &	0.77 &	\textbf{0.531} &	0.98 &	0.567 & 0.04  & 0.586 & 1 & 0.007\\
3 &	9, (3, 3, 3) &	0.520 &	0.88 &	\textbf{0.512} &	0.58 &	0.571 & 0.04 & 0.588 & 1 & 0.012\\
\bottomrule
\end{tabular}
\caption{Performance of upper bounds on the pitprops dataset ($p=13$), as we vary the overall sparsity ($k$), the number of PCs ($r$), and the allocation of a sparsity budget to the different PCs. We use a 7,200-second time limit. We denote the best-performing solution (least upper bound) in bold. Note that all results are normalized by dividing by the trace of $\bm{\Sigma}$, i.e., $p$, the number of features, to report results in terms of the proportion of variance explained.} 
\label{tab:comparison_bounds_new}
\end{table}%

On the instances presented in Table \ref{tab:comparison_bounds_new}, we observe that the semidefinite relaxations often terminate in less than a second and ``Perm-Ineq'' is uniformly the strongest relaxation. In particular, when individual sparsity budgets $k_t$ are given, Perm-Ineq provides uniformly and sometimes significantly tighter bounds than Extended-Ineq. However, the Lagrangian bound is often only weaker at the third decimal point and takes one order of magnitude less time to compute, {\color{black}suggesting it may scale better}. Finally, the combinatorial upper bound is worse than the bounds from all other methods considered for each instance.

We remind the reader that Perm-Ineq {and the Lagrangian bound} cannot compute a bound if we {do not} specify {component-specific 
sparsity budgets $k_t$}. Nonetheless, {when only an overall sparsity $k$ is imposed,} Extended-Ineq {and the combinatorial bound} are much weaker than the worst-case bound over all possible allocations {\color{black}$\{k_t\}_{t \in [r]}$ with $\sum_{t \in [r]}k_t=k$}. 
So, time permitting, we recommend computing the {\color{black}upper} bound by solving the relaxations for all possible allocations of $k_t$ and taking the worst-case bound. For instance, for $r=2, k=10$ with Perm-Ineq, this approach would give an upper bound of 0.459, which is $6.3\%$ better than Extended-Ineq's bound with $k$ alone.
Accordingly, in the rest of the paper, we only consider instances of Problem \eqref{prob:spcaorth_compact} where we know both $k$ and $k_t$ and consider Perm-Ineq and its second-order cone relaxations, but not Extended-Ineq.

\paragraph{Benchmarking on Larger-Scale Datasets} We now investigate the scalability of the SDP relaxation Perm-Ineq (``PSD''), its second-order cone relaxation as described in Section \ref{ssec:ec.soc} (``SOC'' as in Equation \eqref{prob:disjunctiverelax_permutationinvariant_soc}), and the aforementioned combinatorial and Lagrangian bounds on larger UCI datasets in Table \ref{tab:comparison_bounds_ucinew}. Note that for the micromass dataset, we only include sets of second-order cone constraints with fewer than $O(p^2)$ members in the SOC relaxation, to avoid excessively memory-intensive problems.

\begin{table}[h!]
\centering\footnotesize
\begin{tabular}{@{}l r r l r r r r r r r r r r r r@{}} \toprule
Dataset & Dim. ($p$) & Rank ($r$) &  Sparsity ($k$, $k_t$) & \multicolumn{2}{c@{\hspace{0mm}}}{PSD} & \multicolumn{2}{c@{\hspace{0mm}}}{SOC} & \multicolumn{2}{c@{\hspace{0mm}}}{Lagrangian} & \multicolumn{2}{c@{\hspace{0mm}}}{Combinatorial} \\
\cmidrule(l){5-6} \cmidrule(l){7-8} \cmidrule(l){9-10} \cmidrule(l){11-12}  &  & & & UB & T(s) & UB & T(s) & UB & T(s) & UB & T(s) \\\midrule
Pitprops &	13 &	2 &	10, (5, 5) &	0.449 &	17.81 & 	0.524 & 	0.22 & 0.524 & 1.44 & 0.559 & 2.12\\ 
 &	 &	2 &	20, (10, 10) &	0.507 &	0.59 &	0.672 & 	0.23 & 0.642 & 0.02 & 0.803 & 0.01\\
 &	 &	3 &	15, (5, 5, 5) &	0.616 &	0.87 &	0.761 & 	0.42  & 0.786 & 0.05 & 0.827 & 0.01\\
 &	 &	3 &	30, (10, 10, 10) &	0.652 &	0.72 &	1.007 & 	0.46  & 0.963 & 0.04 & 1.198 & 0.01\\ \midrule
Wine &	13 &	2 &	10, (5, 5) &	0.458 &	0.57 &	0.529 & 	0.25  & 0.529 & 0.03 & 0.579 & 0.01 \\ 
 &	 &	2 &	20, (10, 10) &	0.554 &	0.55 & 	0.722 & 	0.22  & 0.707 & 0.03 &  0.876 & 0.01 \\ 
 &	 &	3 &	15, (5, 5, 5) &	0.632 &	0.79 & 	0.762 & 	0.52  & 0.794 & 0.05 & 0.853 & 0.01\\ 
 &	 &	3 &	30, (10, 10, 10) &	0.665 &	0.74 &	1.083 & 	0.57  & 1.060 & 0.03 & 1.296 & 0.01\\ \midrule
Ionosphere &	34 &	2 &	10, (5, 5) &	0.209 &	8.22 & 	0.221 & 	1.45  & 0.221 & 0.10 &  0.228 & 0.03\\ 
 &	 &	2 &	20, (10, 10) &	0.305 &	8.35 & 	0.363 & 	1.53  & 0.361 & 0.04 & 0.401 & 0.02\\ 
 &	 &	2 &	40, (20, 20) &	0.378 &	9.84 & 	0.504 & 	1.85  & 0.500 & 0.05 &  0.618 & 0.02\\ 
 &	 &	3 &	15, (5, 5, 5) &	0.297 &	12.74 &	0.331 & 	4.39  & 0.331 & 0.05 & 0.340 & 0.04\\ 
 &	 &	3 &	30, (10, 10, 10) &	0.411 &	11.93 & 	0.545 & 	4.49  & 0.542 & 0.07 & 0.597 & 0.04\\ 
 &	 &	3 &	60, (20, 20, 20) &	0.464 &	13.38 & 	0.757 & 	5.95  & 0.749 & 0.08 & 0.920 & 0.03\\ \midrule
Communities &	101 &	2 &	10, (5,5) & 0.095 &	373.2 & 	0.096 & 	29.02  & 0.096 & 0.21 & 0.097 & 1.33\\ 
 &    & 2 & 20, (10, 10) &	0.169 &	307.7 & 	0.175 & 	44.83  & 0.175 & 0.29 & 0.180 & 0.89\\ 
 &	 &	2 &	40, (20, 20) &	0.263 &	373.9 & 	0.286 & 	40.13  & 0.284 & 0.52 & 0.320 & 0.72\\ 
 &	 &	3 &	15, (5, 5, 5) &	0.141 &	441.6 & 	0.144 & 	80.30  & 0.144 & 0.26 & 0.146 & 0.25\\ 
 &	 &	3 &	30, (10, 10, 10) &	0.245 &	564.2 & 	0.262 & 	65.25  & 0.262 & 0.27 & 0.270 & 0.24\\ 
 &	 &	3 &	60, (20, 20, 20) &	0.378 &	447.3 & 	0.429 & 	54.71  & 0.425 & 0.30 & 0.476 & 0.26\\\midrule
Arrhythmia &	274 &	2 &	10, (5, 5) &	-&	- & 	0.031 & 	54.71  & 0.031 & 0.60 & 0.032 & 10.83\\ 
 &	 &	2 &	20, (10, 10) &	- &	- & 	0.055 & 	456.0  & 0.055 & 0.31 & 0.060 & 7.69\\ 
 &	 &	2 &	40, (20, 20) &	- &	- & 	0.086 & 	445.8  & 0.084 & 0.39 & 0.105 & 5.49\\ 
 &	 &	3 &	15, (5, 5, 5) &	- &	- & 	0.047 & 	671.8  & 0.046 & 0.43 & 0.049 & 30.80\\ 
 &	 &	3 &	30, (10, 10, 10) &	- &	- & 	0.083 & 	812.2  & 0.083 & 0.57 & 0.089 & 11.32\\ 
 &	 &	3 &	60, (20, 20, 20) &	- &	- & 	0.129 & 	803.7  & 0.126 & 0.64 & 0.155 & 8.54\\ \midrule
 Micromass &	1300 &	2 &	10, (5, 5) &	- &	- &	0.008 & 	1089  & 0.008 & 2.98 & 0.008 & 854.6\\ 
 &	 &	2 &	20, (10, 10) &	- &	- & 	0.015 & 	13620 & 0.014 & 4.86 & 0.014 & 178.3\\ 
 &	 &	2 &	40, (20, 20) &	- &	- &		0.027 & 	9213  & 0.023 & 4.03 & 0.025 & 154.5\\ 
 &	 &	3 &	15, (5, 5, 5) &	- &	- &	0.012 & 	7953  & 0.011 & 5.40 & 0.012 & $>7200$\\ 
 &	 &	3 &	30, (10, 10, 10) &	- &	- &		0.023 & 	19640  & 0.021 & 3.51 & 0.021 & 1092\\ 
 &	 &	3 &	60, (20, 20, 20) &	- &	- &		0.043 & 	18630  & 0.034 & 4.39 & 0.038 & 649.0\\ 
\bottomrule
\end{tabular}
\caption{Performance of conic and combinatorial bounds across UCI datasets. All bounds are normalized by dividing by $p=\mathrm{tr}(\bm{\Sigma})$, the number of features, to report in terms of proportion of correlation explained. The notation ''-'' denotes that an instance could not be solved using the provided memory budget, namely $32$ GB for instances where $p \leq 101$, $100$ GB for instances where $p \in [102, 250]$, $370$ GB for instances where $p>250$.} 
\label{tab:comparison_bounds_ucinew}
\end{table}%

We observe that the PSD relaxation uniformly dominates all other upper bounds and can be solved within a few minutes for $p \approx 100$ but quickly requires a prohibitive amount of memory and time at higher dimensions. The SOC relaxation scales up to $p \approx 1000$ in hours but provides a weaker upper bound, within $1$--$20\%$ of PSD. On the other hand, the Lagrangian relaxation scales to $p \approx 1000$ in seconds, and outperforms the SOC relaxation's bound on every instance where $\sum_t k_t \leq p$. 
Finally, the combinatorial bound is uniformly worse than the Lagrangian bound in terms of both runtime and the quality of the bound. 
Thus, in practice, we recommend using the SDP bound if it can be computed (say $p\leq 300$), and the Lagrangian bound otherwise.

\subsection{Performance of Feasible Methods}\label{ssec:feasiblemethods}
In this section, we numerically evaluate the quality of our three algorithms in terms of their ability to recover approximately orthogonal and high-quality principal components on real-world datasets. We first validate that our methods are capable of recovering feasible and near-optimal solutions to small-scale sparse PCA problems with multiple PCs. We then compare our algorithms with five state-of-the-art techniques on eleven UCI datasets. 
Since Algorithm \ref{alg:alternatingmin} and most of the benchmarked algorithms from the literature require that the sparsity of each PC is specified separately, we consider this formulation in this section (and fix $k_t = k/r$ for concision).

\paragraph{Benchmarking on Pitprops Data:}
We first investigate the performance of Algorithms \ref{alg:greedymethod2}--\ref{alg:disjoint.linalg} on the \verb|pitprops| dataset because 
\cite{lu2012augmented} used this dataset to extensively benchmark several algorithms for sparse PCA with $r=6$ PCs. All six methods they benchmarked required an overall sparsity of around $45$--$60$ to explain less than $70\%$ of the {\color{black}correlation}. The best performing method could explain 69.55\% of the {\color{black}correlation} with an overall sparsity of 46 \citep[Table 11]{lu2012augmented}. They concluded that ``there do not exist six highly sparse, nearly orthogonal and uncorrelated PCs while explaining most of variance'' {\color{black}(note that the pitprops dataset is normalized, and hence variance should be read as correlation here)}.

As reported in Table \ref{tab:comparison_feasible_pitprops}, with $r=6$ PCs and $k_t = 2$, hence an overall sparsity of $12$, solutions returned by Algorithms \ref{alg:alternatingmin}--\ref{alg:disjoint.linalg} explain $74\%$--$75\%$ of the {\color{black}correlation}, and our SDP upper bound from Algorithm \ref{alg:greedymethod2} proves that it is not possible to explain more than $75\%$ of the {\color{black}correlation} in this case. Moreover, Algorithm \ref{alg:alternatingmin} even provides a solution that explains $81\%$ of the {\color{black}correlation} with an overall sparsity of $6 \times 4=24$. {\color{black}In other words, the solutions that were previously believed to be unattainable at this sparsity level are achievable with our methods, and can be computed within seconds on this instance.}

\paragraph{Benchmarking on Larger-Scale Datasets} We now investigate the performance of Algorithms \ref{alg:greedymethod2}--\ref{alg:disjoint.linalg} on eleven UCI datasets summarized in Table \ref{tab:dataset_summary}, whose dimension{\color{black}s} range from $p=13$ (\verb|pitprops|) to $p=1300$ (\verb|micromass|). We compare with four state-of-the-art methods. Namely, 
\begin{itemize}
    \item The branch-and-bound method of \cite{berk2019certifiably} for optimally computing one sparse PC, combined with the deflation scheme of \cite{mackey2008deflation} to obtain multiple PCs, implemented in \verb|Julia| and made available at \verb|github.com/lauren897/Optimal-SPCA|. According to \citet[]{berk2019certifiably}, this method outperformed four others across three UCI datasets ($r=3, k=5$).
    \item The deflation method of \cite{hein2010inverse}, using the custom deflation method developed in \cite{buhler2014flexible}, implemented in \verb|Matlab| and made publicly available at \verb|github.com/tbuehler/sparsePCA|, using default parameters. This approach was found by \citet[Table 9]{berk2019certifiably} to be second-best of the methods in their comparison.
    \item The \verb|Lasso|-inspired method of \cite{zou2006sparse}, using the \verb|spca| function in the \verb|elasticnet| package version $1.3$, using default parameters. This approach is perhaps the most commonly used one in practice, since it is distributed via the ubiquitous \verb|elasticnet| package.
    \item The covariance thresholding method of \cite{deshpande2014sparse}, building upon the works of \cite{krauthgamer2015semidefinite}, which relies on applying a soft-thresholding operation first on the entries of the covariance matrix and then on its $r$ leading eigenvectors. We implemented this method natively in \verb|Julia| and release it as part of our codebase.
\end{itemize}
Finally, we also compare with Gurobi's non-convex branch-and-bound solver as a benchmark, as described in Section \ref{sec:ec.branchandbound}. We do not report upper bounds for this benchmark, because they are uniformly above $1$, as shown for the pitprops dataset in Tables \ref{tab:comparison_gurobi_pitprops1}--\ref{tab:comparison_gurobi_pitprops2}, and sometimes no upper bound can be computed within the time limit, giving an average upper bound of $+\infty$.

Algorithm \ref{alg:greedymethod2} involves solving a convex relaxation from Section \ref{sec:relaxations}. Based on the scalability results presented in the previous section, we use different convex relaxations depending on the dimensionality of the problem. Namely, for Algorithm \ref{alg:greedymethod2}, we use the full semidefinite relaxation \eqref{prob:disjunctiverelax_permutationinvariant} if $p \leq 101$, and the SOC relaxation \eqref{prob:disjunctiverelax_permutationinvariant_soc} otherwise. 

We report summary results in Tables \ref{tab:ucilargersummarized}-\ref{tab:ucilargersummarized2} (see also Tables \ref{tab:comparison_feasiblemethods_uci}--\ref{tab:comparison_feasiblemethods_uci3_part4} in Section \ref{ssec:largeruci} for instance-wise results). In particular, we report the average objective value attained by each method, the average upper bound attained by each method (where applicable), and the average gap between the feasible solution and the best upper bound found by \textit{any} of our three methods. Our results are broken down by whether $\sum_t k_t\leq p$ or $\sum_t k_t>p$, since Algorithms \ref{alg:greedymethod2} and \ref{alg:disjoint.linalg} are incapable of returning PCs that have a collective sparsity larger than $p$.  Finally, when the returned solution violates the orthogonality condition, its objective value is not necessarily a valid bound on the objective of \eqref{prob:spcaorth_compact} and the average reported gap is thus an optimistic estimate.

\begin{table}[h!]
\centering\footnotesize
\begin{tabular}{@{}l r r r r r@{}} \toprule
Method & Obj. & UB & Rel. gap* ($\%$) & Viol. & T(s)  \\\midrule
Algorithm \ref{alg:greedymethod2} & 0.135 & 0.168 & $17.99\%$ & 0 & 2220\\
Algorithm \ref{alg:alternatingmin} & 0.144 & 0.177 & $8.82\%$ & 0 & 70.66\\
Algorithm \ref{alg:disjoint.linalg} & 0.151 & 0.188 & $3.11\%$& 0 & 126.9\\
Branch-and-bound & 0.137 & - & $14.73\%$ & 0 & 6945\\
\cite{berk2019certifiably} & 0.152 & - & $2.63\%$ & 0.012 & 29.81\\
\cite{deshpande2014sparse} & 0.144 & - & $7.73\%$ & 0.063 & 14.09\\
\cite{hein2010inverse} & 0.124 & - & $25.75\%$ & 0.016 & 0.23 \\
\cite{zou2006sparse} & 0.033 & - & $81.68\%$ & 1.340 & 5.66\\
\bottomrule
\end{tabular}
\caption{Average performance of methods across the 11 UCI datasets described in Table \ref{tab:dataset_summary} with $k \in \{5, 10, 20\}, r \in \{2, 3\}: k \leq p$, conditioning on instances where $kr \leq p$ (in particular, $kr$ ranges from 10 to 60 in our experiments). {\blue*}Optimality gaps computed using the best UB from Algorithms \ref{alg:greedymethod2}--\ref{alg:disjoint.linalg}. Methods that violate orthogonality can achieve a negative gap.}

\label{tab:ucilargersummarized}
\end{table}

\begin{table}[h!]
\centering\footnotesize
\begin{tabular}{@{}l r r r r r@{}} \toprule
Method & Obj. & UB & Rel. gap* ($\%$) & Viol. & T(s)  \\\midrule
Algorithm \ref{alg:greedymethod2} & 0.422 & 0.532 & $21.49\%$ & 0 & 7.84 \\
Algorithm \ref{alg:alternatingmin} & 0.519 & 0.733 & $2.82\%$ & 0 & 6.71\\
Algorithm \ref{alg:disjoint.linalg} & 0.477 & 0.878 & $10.43\%$ & 0 & 0.76\\
Branch-and-bound & 0.485 & - &  $10.35\%$ & 0 & $>7200$ \\
\cite{berk2019certifiably} & 0.519 & - & $2.78\%$ & 0.146 & 0.59\\
\cite{deshpande2014sparse} & 0.510 & - & $4.41\%$ & 0.279 & 0.01\\
\cite{hein2010inverse} & 0.507 & - & $5.26\%$ & 0.066 & 0.032\\
\cite{zou2006sparse} & 0.146 & - & $73.13\%$ & 2.168 & 0.76\\
\bottomrule
\end{tabular}
\caption{Average performance of methods across the 11 UCI datasets described in Table \ref{tab:dataset_summary} with $k \in \{5, 10, 20\}, r \in \{2, 3\}: k \leq p$, conditioning on instances where $kr > p$ (in particular, this criterion excludes datasets with $p \geq 60$). {\blue *}Optimality gaps computed using the best UB from Algorithms \ref{alg:greedymethod2}--\ref{alg:disjoint.linalg}. Methods that violate orthogonality can achieve a negative gap.}
\label{tab:ucilargersummarized2}
\end{table}

First, we observe that our three methods and branch-and-bound are the only ones to return PCs that are systematically orthogonal (with an average orthogonality violation $< 10^{-3}$). Among them, we observe that Algorithm~\ref{alg:alternatingmin} performs the best on datasets where ${\blue \sum_{t \in [r]} k_t = } kr > p$ and Algorithm~\ref{alg:disjoint.linalg} performs best on datasets where $kr \leq p$, with average relative optimality gaps of $2.78\%$ and $3.11\%$ on these respective instances, and average runtimes in the tens or hundreds of seconds for both methods. On the other hand, Algorithm \ref{alg:greedymethod2} uniformly provides the best upper bounds of all methods, although it is dominated by Algorithms~\ref{alg:alternatingmin}--\ref{alg:disjoint.linalg} in terms of solution quality and requires two orders of magnitude more runtime to compute. 

Of the remaining methods, the method of \cite{berk2019certifiably} performs the next best, in terms of explaining a large proportion of the correlation ($0.137$ when $kr \leq p$ and $0.485$ when $kr > p$) while not violating feasibility too significantly ($0.012$ when $kr \leq p$ and $0.146$ when $kr > p$).  However, its higher objective value than Algorithm \ref{alg:disjoint.linalg} on instances where $kr \leq p$ is driven by instances where it achieves an objective value higher than Algorithm \ref{alg:greedymethod2}'s upper bound by violating the orthogonality constraint. {Runtimes of Algorithm{\color{black}s}~\ref{alg:alternatingmin}-\ref{alg:disjoint.linalg} and of the method of \cite{berk2019certifiably} are also of similar order of magnitude, although the method of \cite{berk2019certifiably} is faster on average.}
Finally, the methods of \cite{hein2010inverse, deshpande2014sparse,zou2006sparse} repeatedly violate the orthogonality constraint and explain less correlation than Algorithm~\ref{alg:alternatingmin} on average. 

Finally, we remark that no one method performs best on every instance. Algorithm \ref{alg:alternatingmin} appears to perform the best overall on instances where $kr >p$ and Algorithm \ref{alg:disjoint.linalg} on instances where $kr \leq p$. However, some other methods \citep[e.g.][]{berk2019certifiably} perform best on some instances. These results suggest that both $k$ and the amount of overlap between the optimal PCs impact the performance of each method, and motivate a comparison on synthetic data, where we control the ground truth, in the next section.
 
\subsection{Statistical Recovery on Synthetic Data}\label{ssec:auc}
To evaluate the support recovery ability of each method, we now 
compare the performance of Algorithms \ref{alg:greedymethod2}, \ref{alg:alternatingmin}, {and \ref{alg:disjoint.linalg}} against the same four methods from the literature on synthetic data. 
We use a spiked Wishart model with multiple spikes, as in \cite{deshpande2014information, ding2023subexponential}, to generate data. 
Namely, we consider an underlying covariance matrix of the form $\bm{\Sigma} = \bm{I}_p + \beta \bm{x}_1 \bm{x}_1^\top + \beta \bm{x}_2\bm{x}_2^\top$, where $\beta=1$ is the signal-to-noise ratio. The vectors $\bm{x}_1, \bm{x}_2 \in \{-1,0,1\}^p$ are random $k_{\text{true}}$-sparse orthogonal vectors. We also control the proportion of overlap between the supports of $\bm{x}_1$ and $\bm{x}_2$, $q \in [0,1]$ ($q=0$ corresponds to disjoint support while $q=1$ corresponds to row sparsity). In our experiments, we take $p=50$, $k_{true}=20$ and vary $q \in \{ 0.1, 0.5, 0.9\}$.
Finally, we sample $n$ observations from a multivariate centered normal distribution with covariance matrix $\bm{\Sigma}$ and construct the empirical covariance matrix $\hat{\bm{\Sigma}}$. We investigate the performance of different methods as $n$ increases (so that the empirical covariance matrix converges to the underlying truth, $\bm{\Sigma}$).
We impose a limit of $100$ iterations for Algorithm \ref{alg:alternatingmin}.

\begin{figure}
\centering
\begin{subfigure}[t]{\linewidth}
    \centering
    \includegraphics[width=0.45\textwidth]{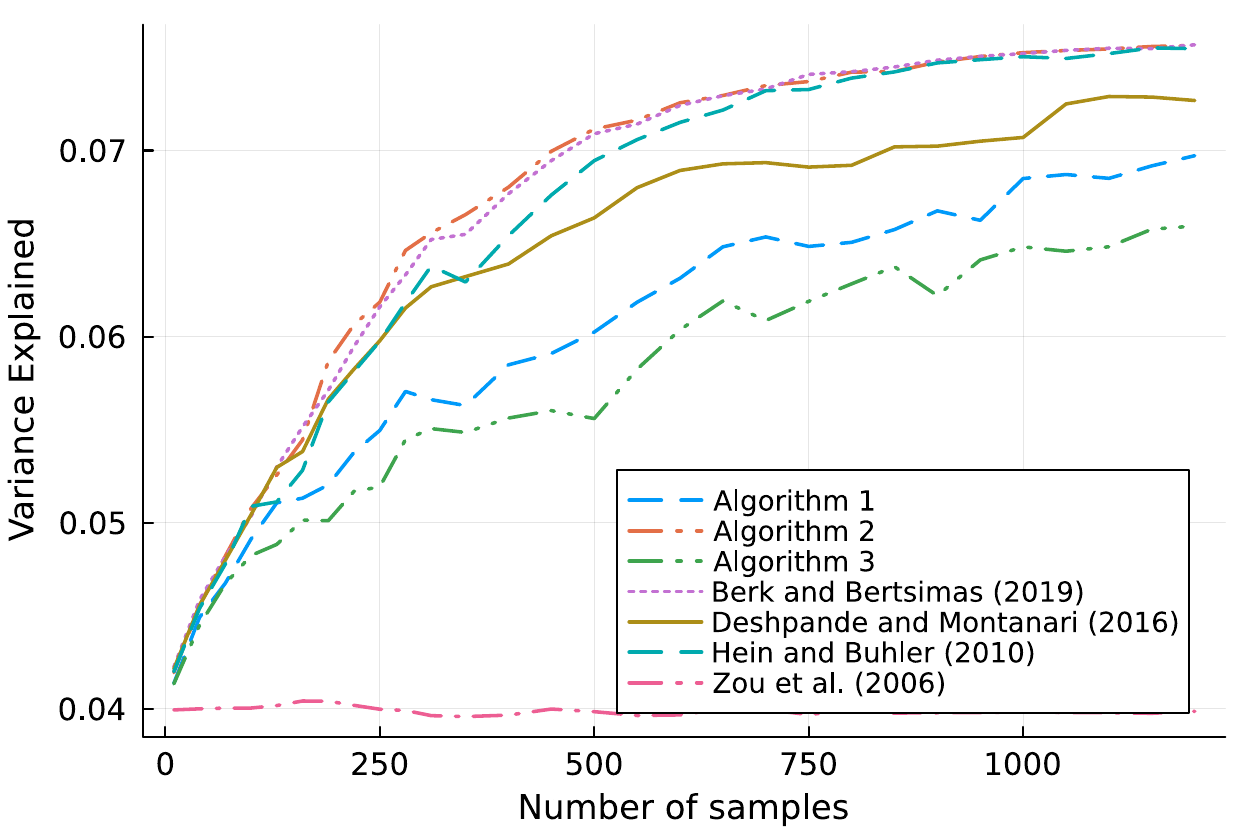} %
    \includegraphics[width=0.45\textwidth]{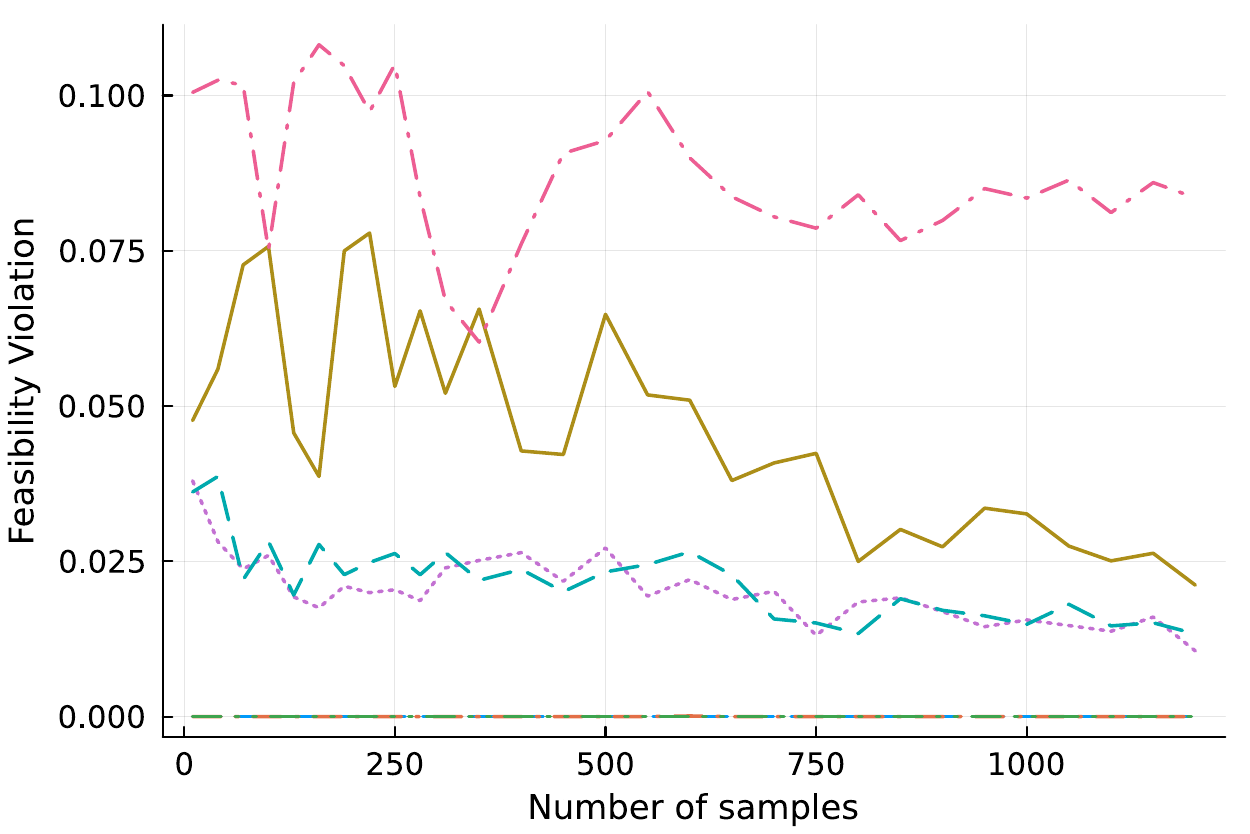}
    \subcaption{Support overlap $q=0.1$}
\end{subfigure}
\begin{subfigure}[t]{\linewidth}
    \centering
    \includegraphics[width=0.45\textwidth]{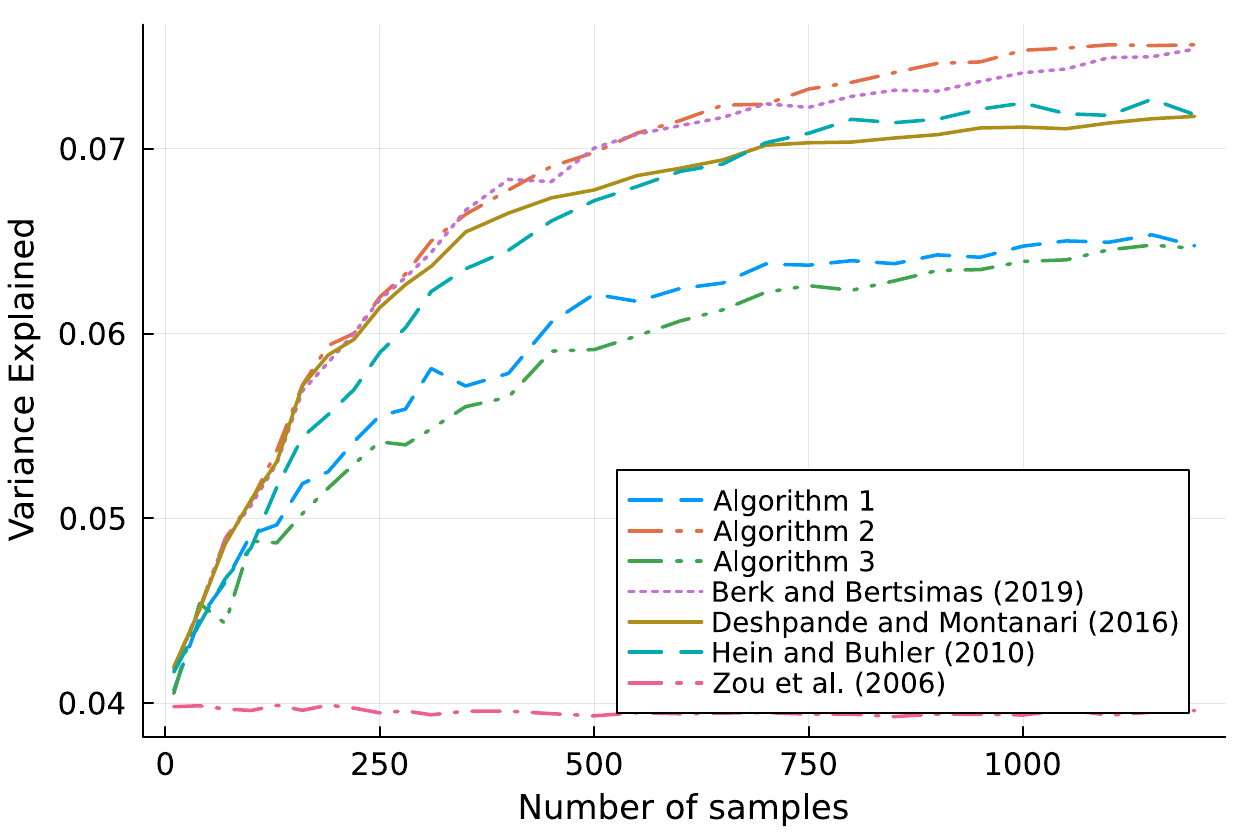} %
    \includegraphics[width=0.45\textwidth]{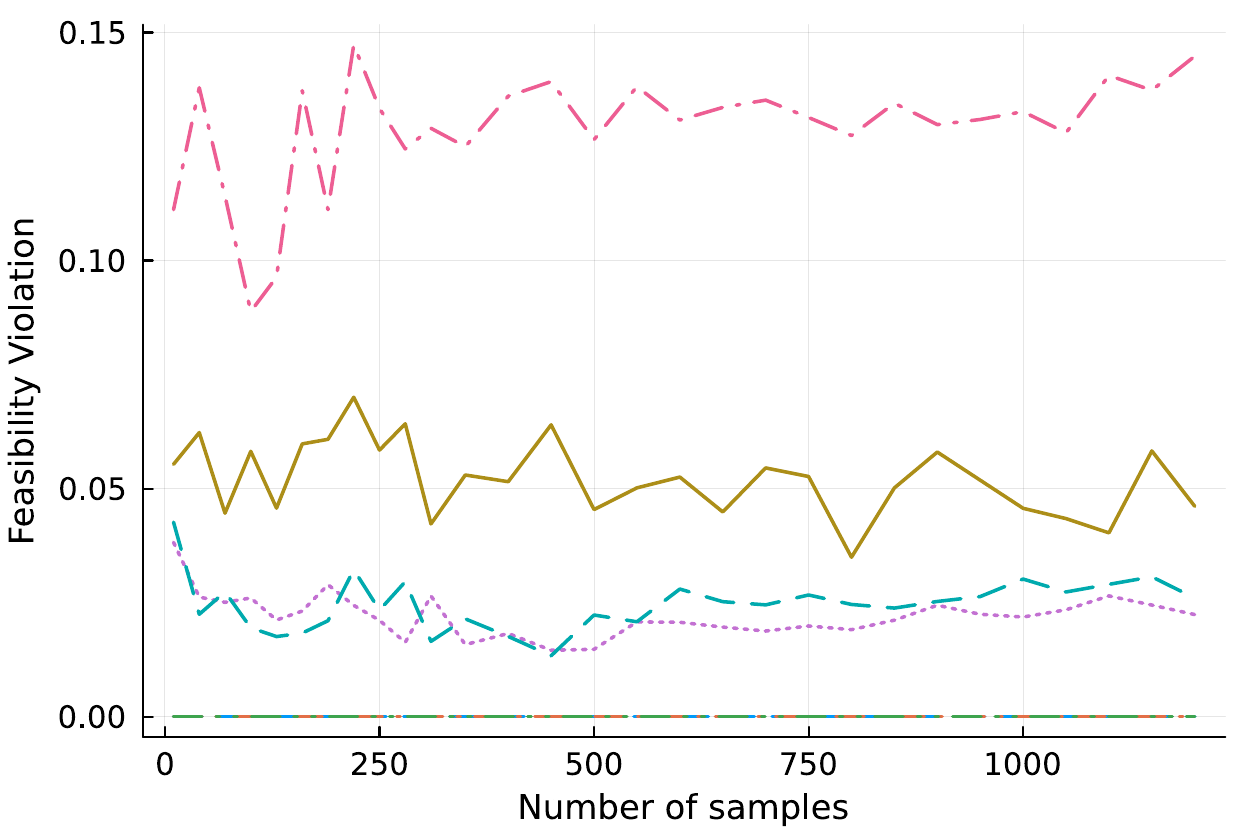}
    \subcaption{Support overlap $q=0.5$}
\end{subfigure}
\begin{subfigure}[t]{\linewidth}
    \centering
    \includegraphics[width=0.45\textwidth]{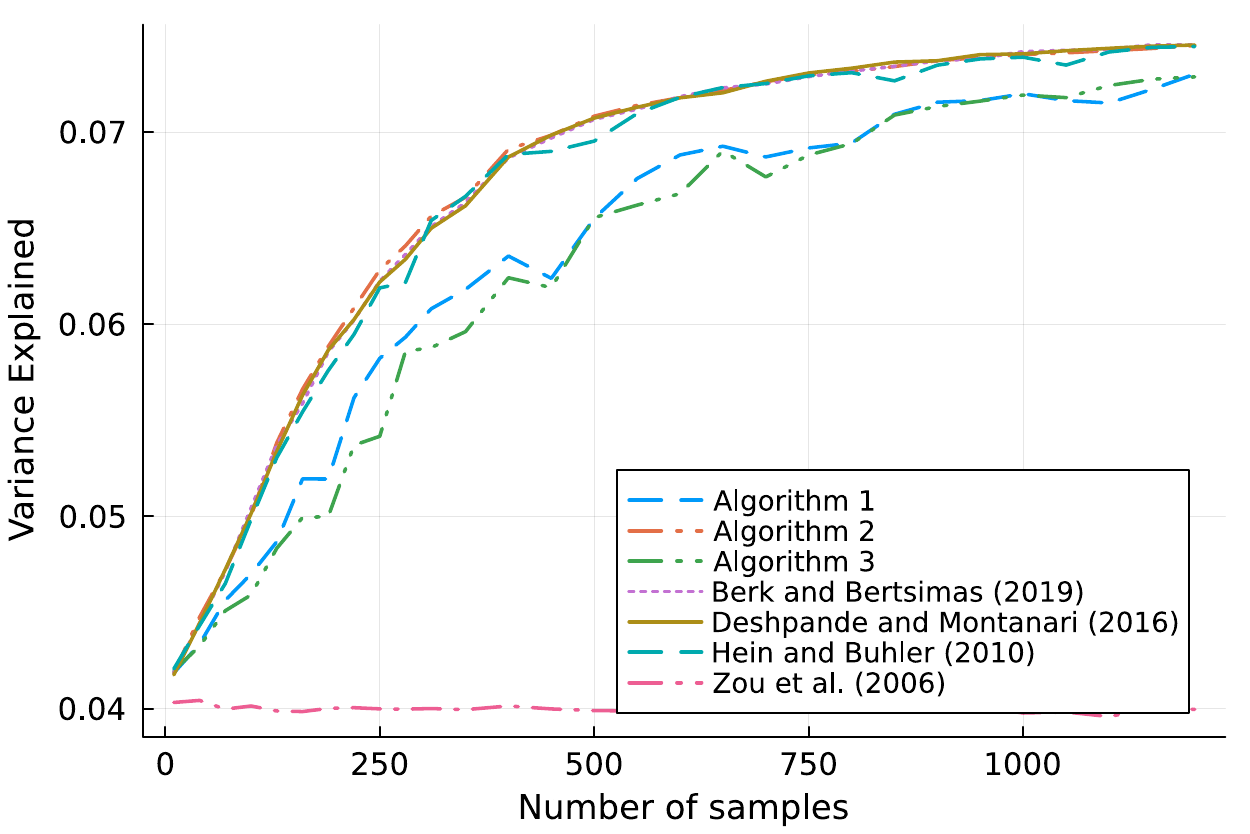} %
    \includegraphics[width=0.45\textwidth]{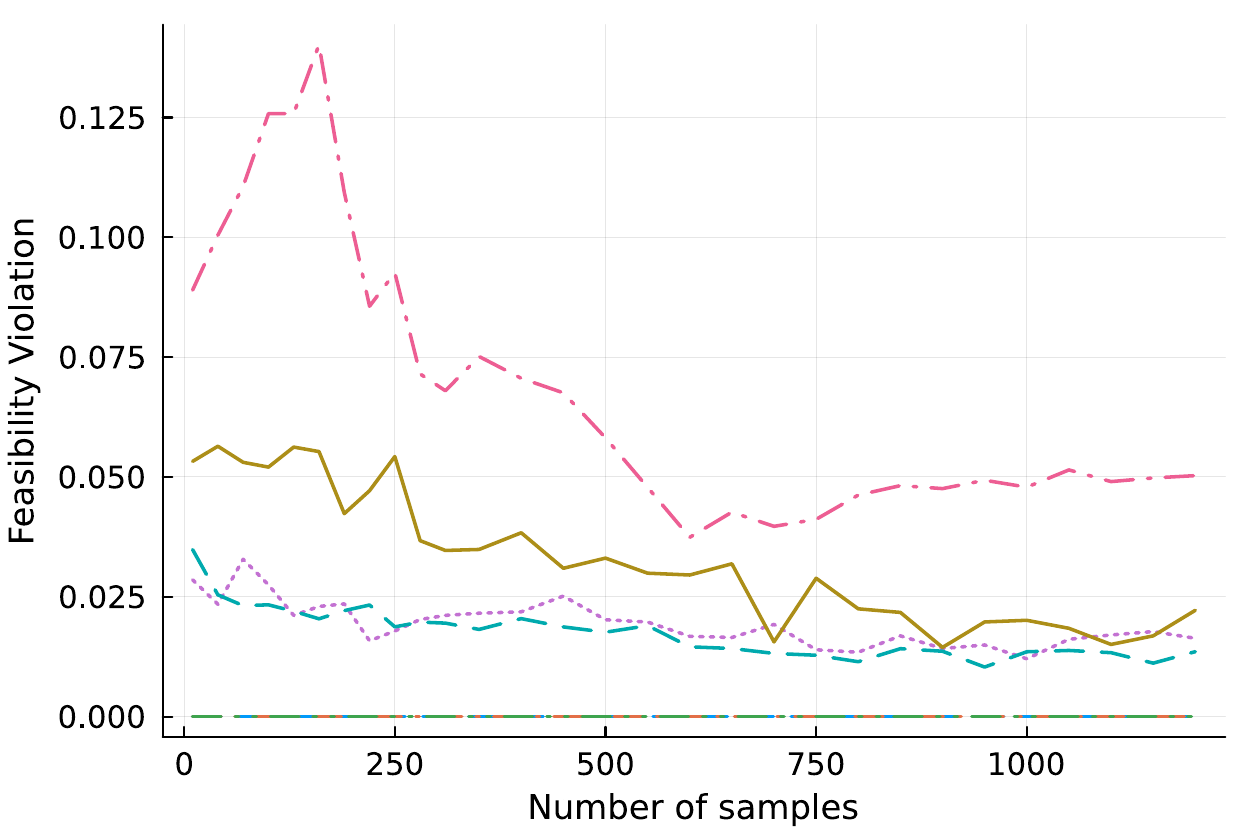}
    \subcaption{Support overlap $q=0.9$}
\end{subfigure}
\caption{Variance explained (left panel) and feasibility violation (right panel) on synthetic instances of sparse PCA with two 20-sparse PCs with partially overlapping support. Results are averaged over 20 replications.}
   \label{fig:syndata.pred}
\end{figure}

For each algorithm, we compute the fraction of variance explained and the feasibility violation (i.e., the inner product between the two PCs computed). We average these performance metrics over 20 random instances and report them in Figure \ref{fig:syndata.pred}. The method of \citet{zou2006sparse} is clearly dominated by all other methods since it explains a significantly lower fraction of the variance, while returning the least orthogonal vectors. 
In terms of objective value (left panel), we observe that Algorithm \ref{alg:alternatingmin}, \citet{berk2019certifiably}, \citet{hein2010inverse} perform almost identically, followed closely by covariance thresholding. They all explain a larger fraction of the variance than Algorithms \ref{alg:greedymethod2} {and \ref{alg:disjoint.linalg}}.
However, we observe on the right panel that Algorithms \ref{alg:greedymethod2}, \ref{alg:alternatingmin}, {and \ref{alg:disjoint.linalg}} are the only methods to return orthogonal PCs, across all values of $n$. In addition, the gap between the methods (and especially the gap between the four best performing methods) seems to shrink as $q$ increases, i.e., when the overlap between the support increases. 

\begin{figure}
\centering
\begin{subfigure}[t]{\linewidth}
    \centering
    \includegraphics[width=.45\textwidth]{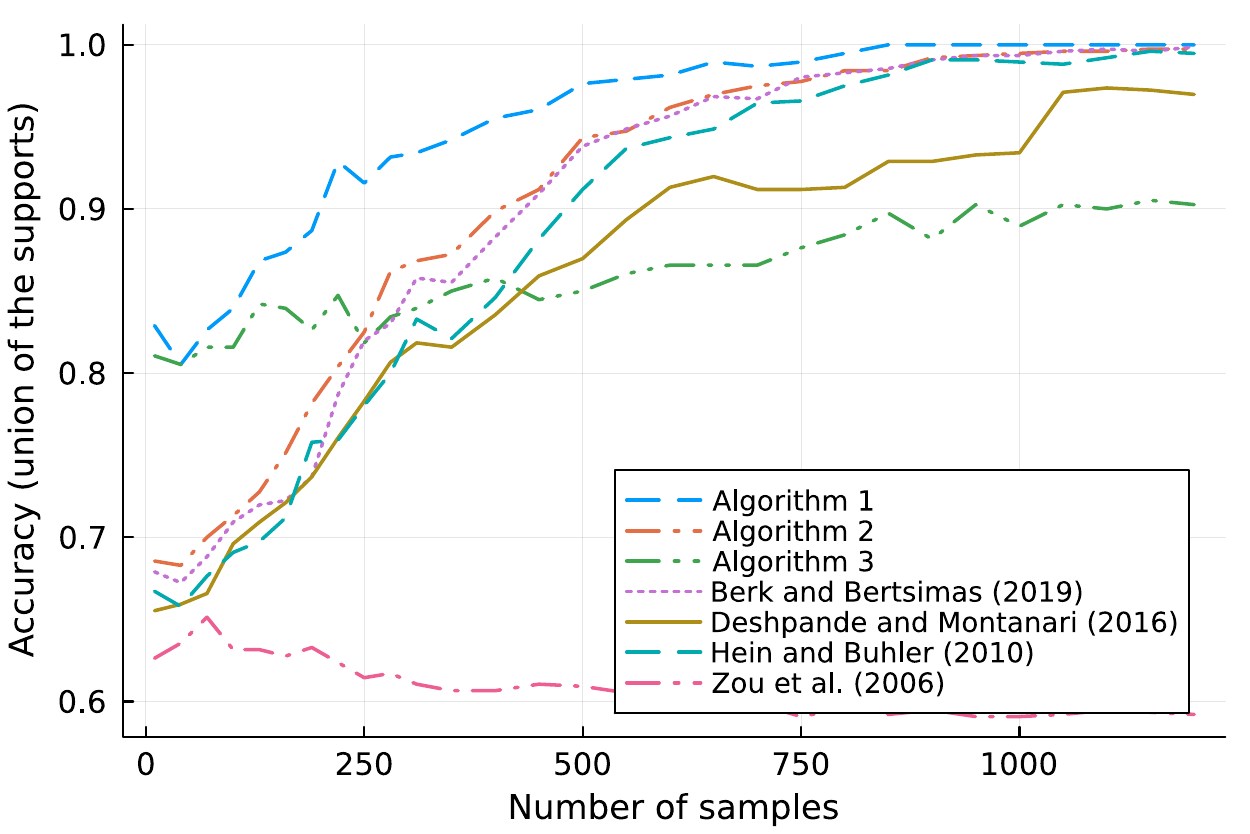} %
    \includegraphics[width=.45\textwidth]{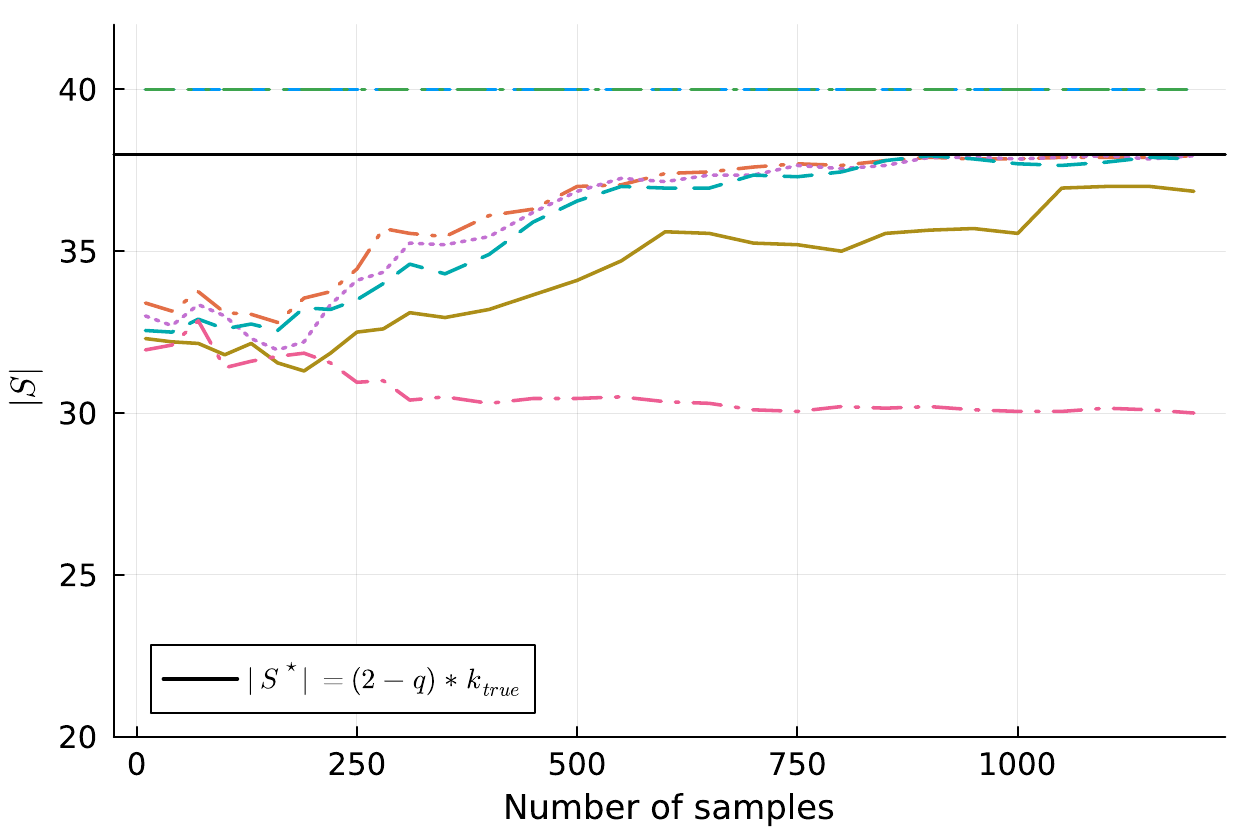}
    \subcaption{Support overlap $q=0.1$}
\end{subfigure} 
\begin{subfigure}[t]{\linewidth}
\centering
    \includegraphics[width=.45\textwidth]{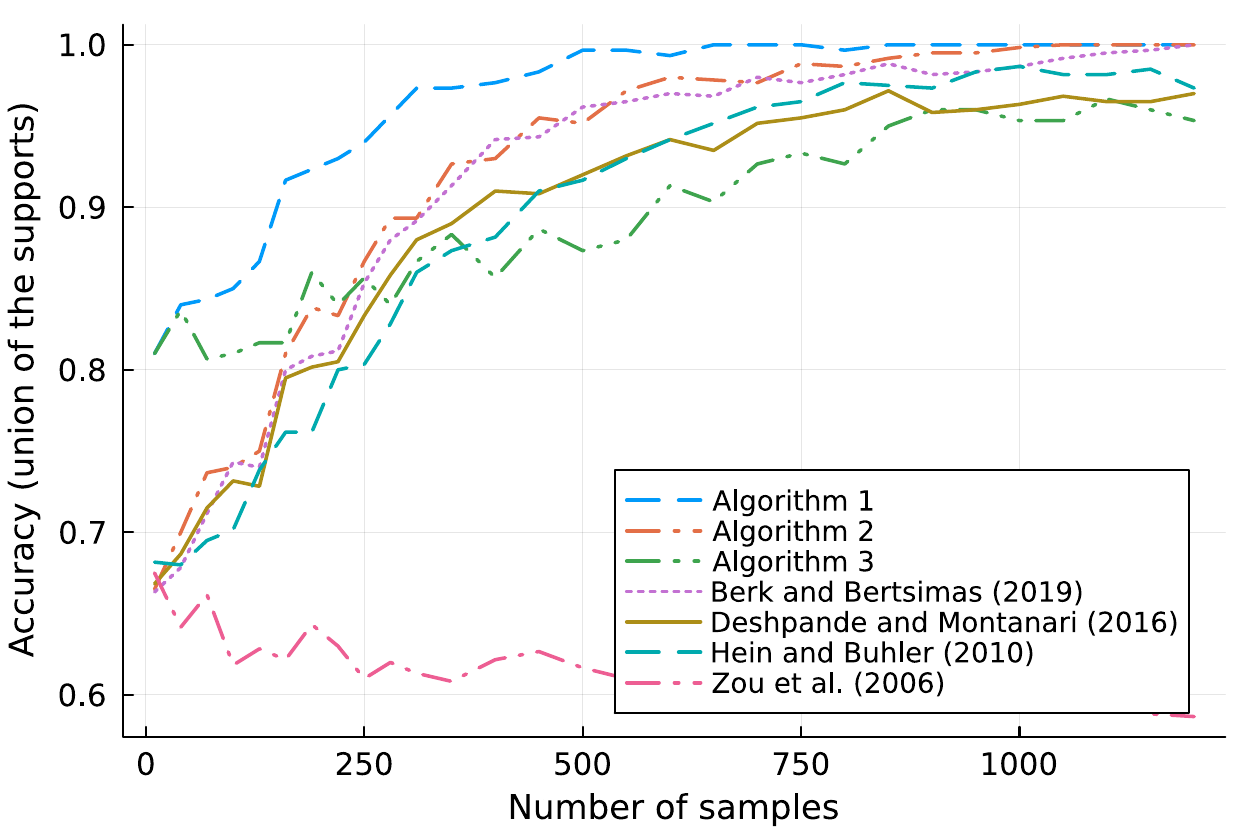} %
    \includegraphics[width=.45\textwidth]{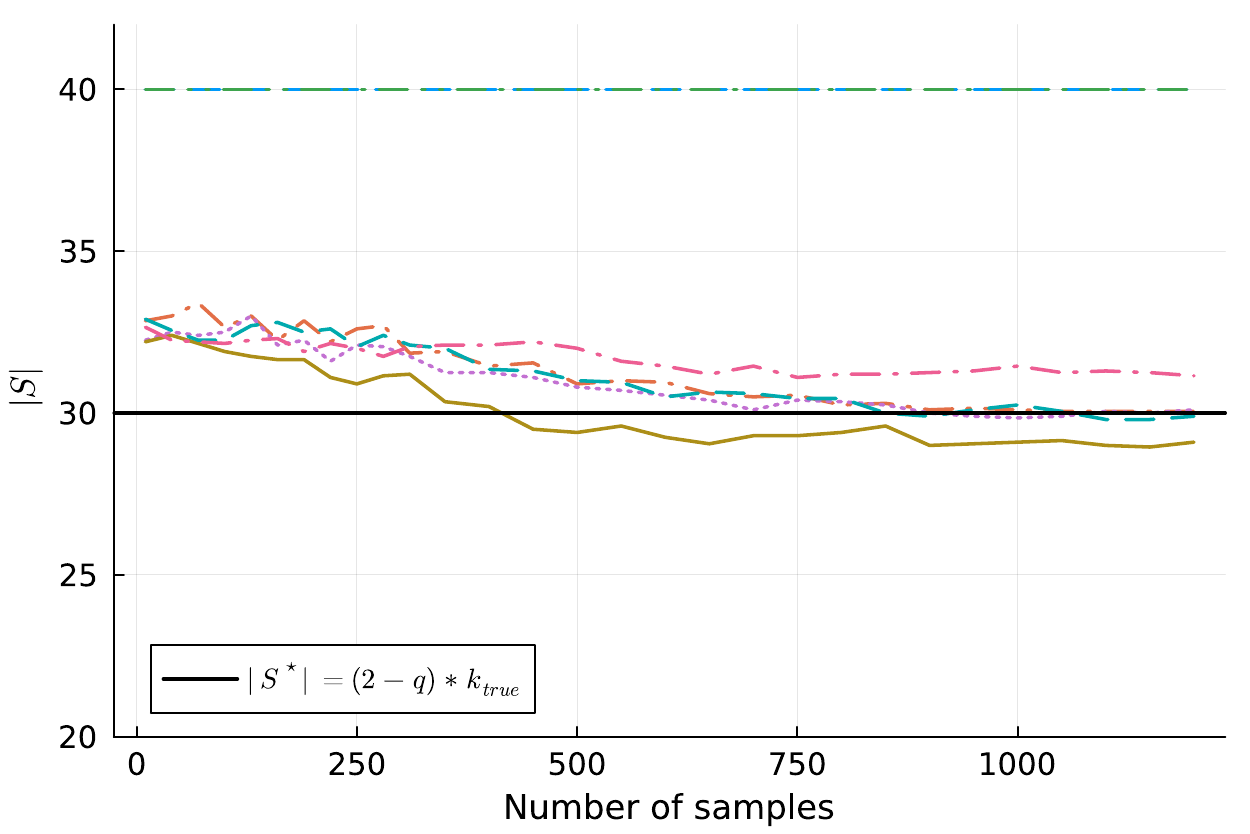} %
    \subcaption{Support overlap $q=0.5$}
\end{subfigure} 
\begin{subfigure}[t]{\linewidth}
    \centering
    \includegraphics[width=.45\textwidth]{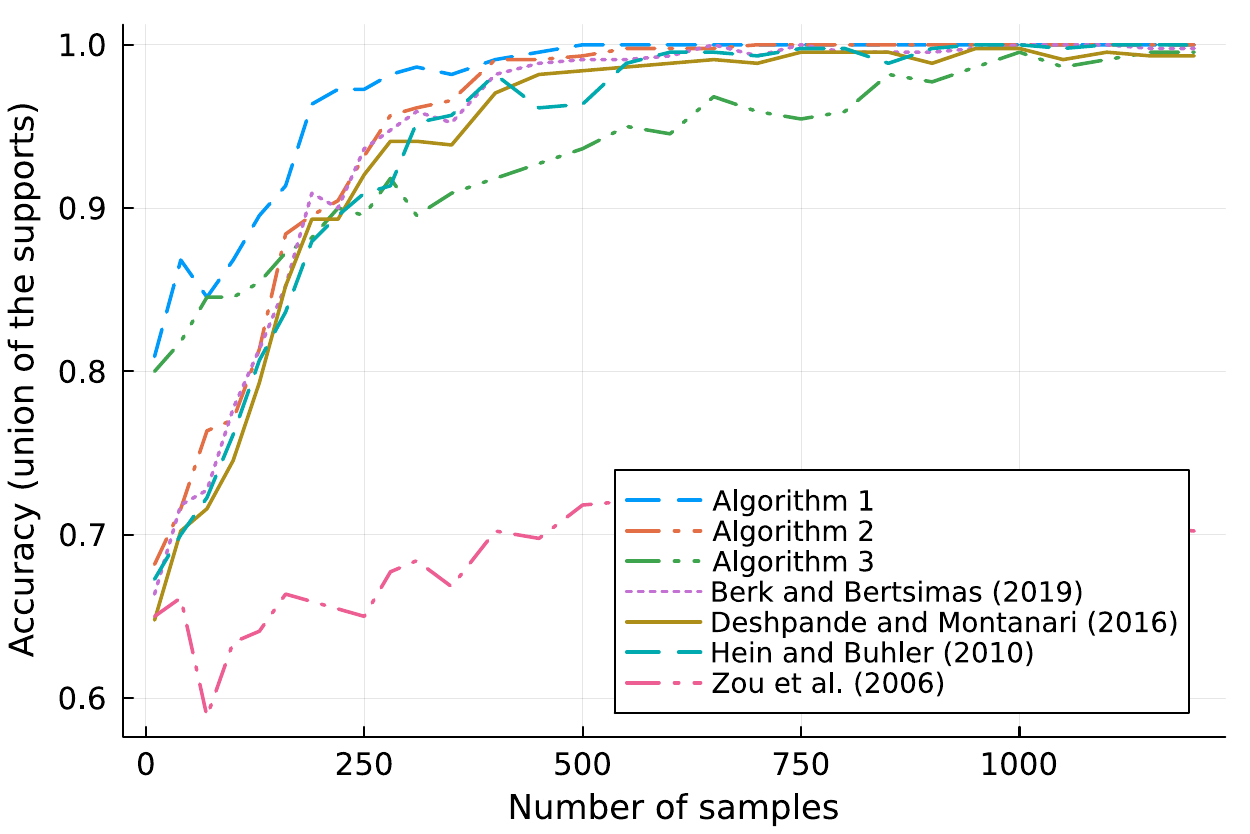}
    \includegraphics[width=.45\textwidth]{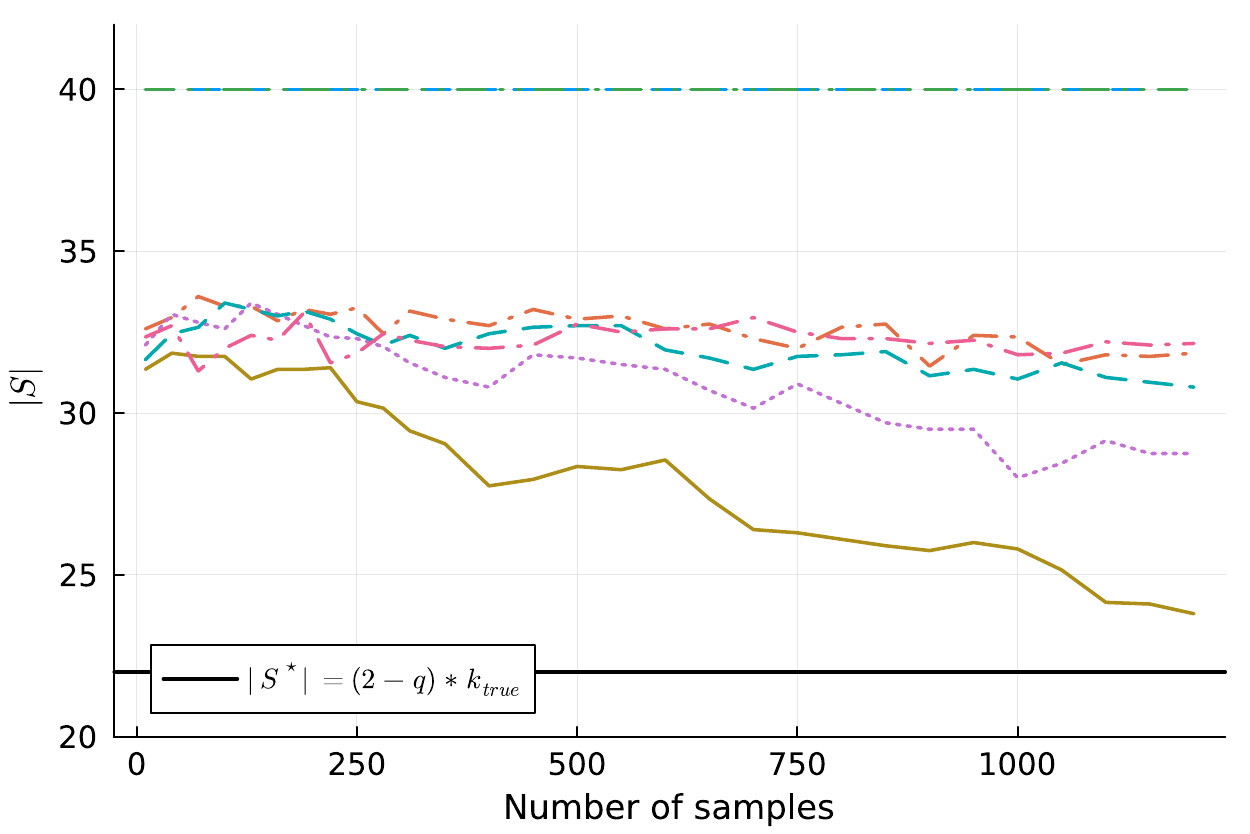}
    \subcaption{Support overlap $q=0.9$}
\end{subfigure}
\caption{Accuracy (left panel) and joint support size (right panel) for the recovery of $\operatorname{supp}(\bm{x}_1) \cup \operatorname{supp}(\bm{x}_2)$, on synthetic instances of sparse PCA with two 20-sparse PCs with partially overlapping support. Results are averaged over 20 replications.}
   \label{fig:syndata.unionrecovery}
\end{figure}

Since the data {\color{black}are} synthetically generated, we can also evaluate the ability to recover the true support. For two $k_{true}$-sparse candidate PCs $\bm{u}_1$ and $\bm{u}_2$, we measure how well $\operatorname{supp}(\bm{u}_1) \cup  \operatorname{supp}(\bm{u}_2)$ recovers the support of $\operatorname{supp}(\bm{x}_1) \cup \operatorname{supp}(\bm{x}_2)$ in terms of accuracy and false detection rate:
\begin{align*}
    A := \: \dfrac{\left| S \cap S^\star \right| }{|S^\star|}, 
\end{align*}
with $S = \operatorname{supp}(\bm{u}_1) \cup \operatorname{supp}(\bm{u}_2)$ and $S^\star = \operatorname{supp}(\bm{x}_1) \cup \operatorname{supp}(\bm{x}_2)$. 
This definition of support recovery corresponds to the one used in statistical studies for sparse PCA with multiple PCs \citep[e.g.,][]{deshpande2014sparse}.
Figure \ref{fig:syndata.unionrecovery} (left panel) reports the value of $A$ for the different algorithms 
as $n$ increases. 
Since we do not explicitly control for the overlap between the returned PCs (except for Algorithms \ref{alg:greedymethod2} {and \ref{alg:disjoint.linalg}}), methods might differ in the size of the support they return, $|S|$. Hence, to allow for a fair comparison, we also report $|S|$ in the right panel of Figure \ref{fig:syndata.unionrecovery}.
Regarding $A$, we observe that Algorithm \ref{alg:greedymethod2} detects a noticeably higher fraction of the true features than other methods, which is not surprising given the fact that it returns disjoint supports. For the remaining methods, their relative performance is aligned with their performance in terms of fraction of variance explained (Figure \ref{fig:syndata.pred}, left panel). 

\newpage \FloatBarrier
\subsection{Specifying the Sparsity Pattern: The Benefits of Asymmetry}\label{ssec:symm}
While Problem \eqref{prob:spcaorth_compact} only requires a bound on the total sparsity, thus allowing flexibility on how this budget is allocated across PCs, 
the worst-case semidefinite upper bound over all sparsity patterns $(k_1,\ldots, k_r): \sum_{t \in [r]}k_t=k$ is often significantly tighter than  the semidefinite relaxation  of \eqref{prob:spcaorth_compact}  with a sparsity budget of $k$ alone, as demonstrated in Section \ref{ssec:bounds}. Moreover, Algorithm \ref{alg:alternatingmin}, which as demonstrated in Sections \ref{ssec:feasiblemethods}--\ref{ssec:auc} is currently the best performing method for obtaining feasible solutions to Problem \eqref{prob:spcaorth_compact}, requires that $(k_1, \ldots, k_t)$ are individually specified. Collectively, these observations suggest that it may be necessary to enumerate all allocations of the sparsity budget $k$, which could be expensive. In our experiments, as is often done in practice, we restricted our search to symmetric allocations. In this section, we revisit the symmetry assumption, investigate when it is justified, and study the relative benefits of asymmetric sparsity budget allocations in terms of obtaining equally sparse sets of PCs that explain more variance.

\begin{figure}[h]
    \centering
    \begin{subfigure}[t]{.45\linewidth}
    \centering
            \includegraphics[width=\textwidth]{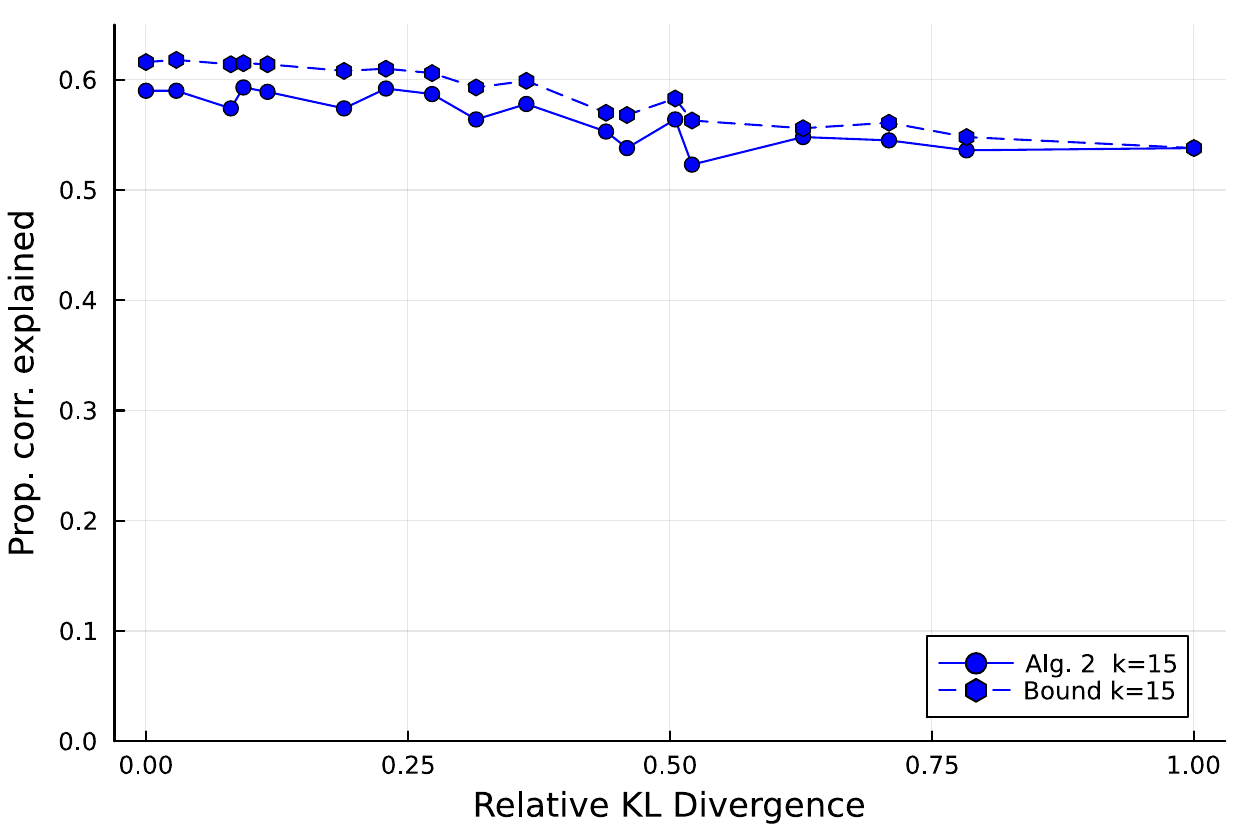}
            \caption{Pitprops data ($k=15$, $r=3$)}
    \end{subfigure}
        \begin{subfigure}[t]{.45\linewidth}
            \includegraphics[width=\textwidth]{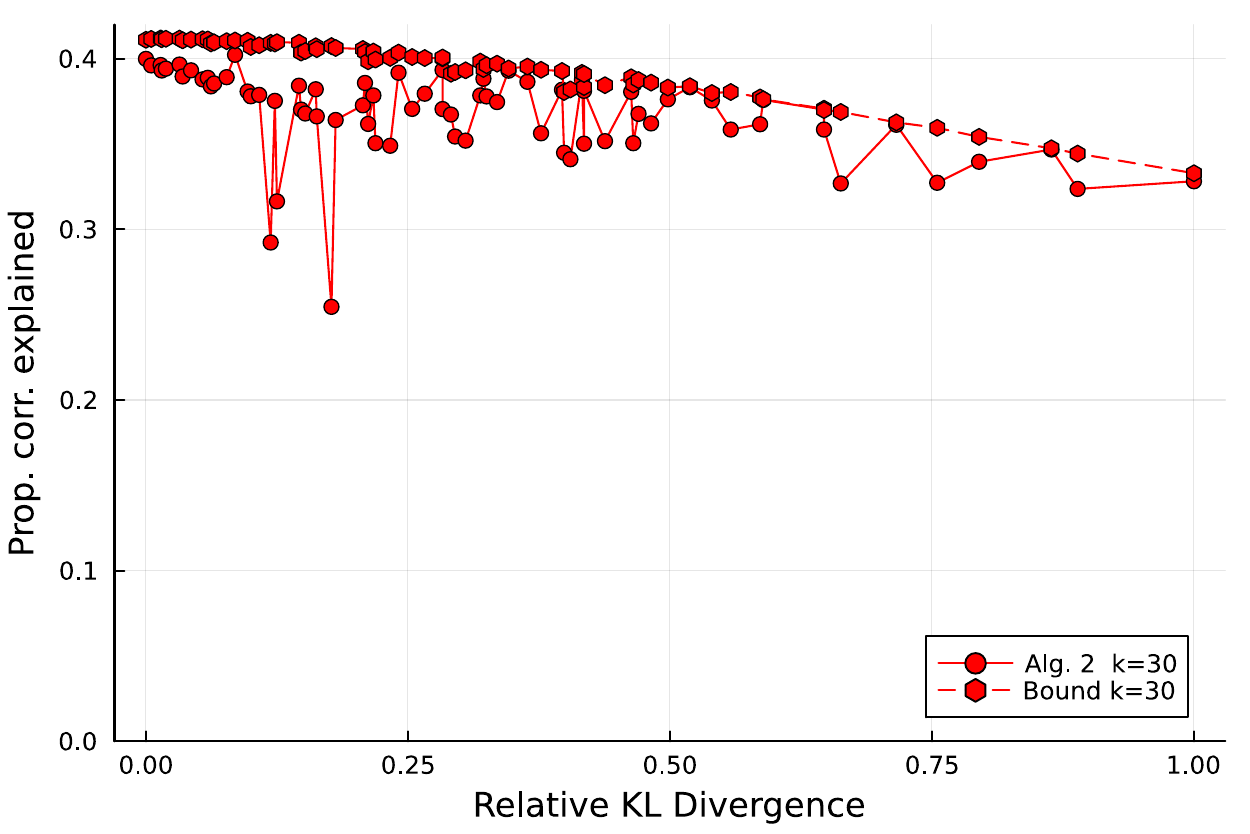}
            \caption{Ionosphere data ($k=30$, $r=3$)}
    \end{subfigure}
   \caption{Asymmetry of sparsity budget allocation (the higher the relative KL divergence, the further away $(k_1,k_2,k_3)$ is from a symmetric allocation) vs. proportion of correlation explained on the pitprops (left panel, $k=15$, $r=3$) and ionosphere (right panel, $k=30$, $r=3$) datasets. For the proportion of correlation explained, we report both an upper bound (obtained from solving our semidefinite relaxation) and a lower bound (obtained from the solution of Algorithm 2). Note that the normalizing constant for the KL divergence is different for each dataset, as the value of $k$ is different. 
   }
   \label{fig:sparsitysymmetry}
\end{figure}

\begin{figure}[h]
    \centering
    \begin{subfigure}[t]{.45\linewidth}
    \centering
            \includegraphics[scale=0.35]{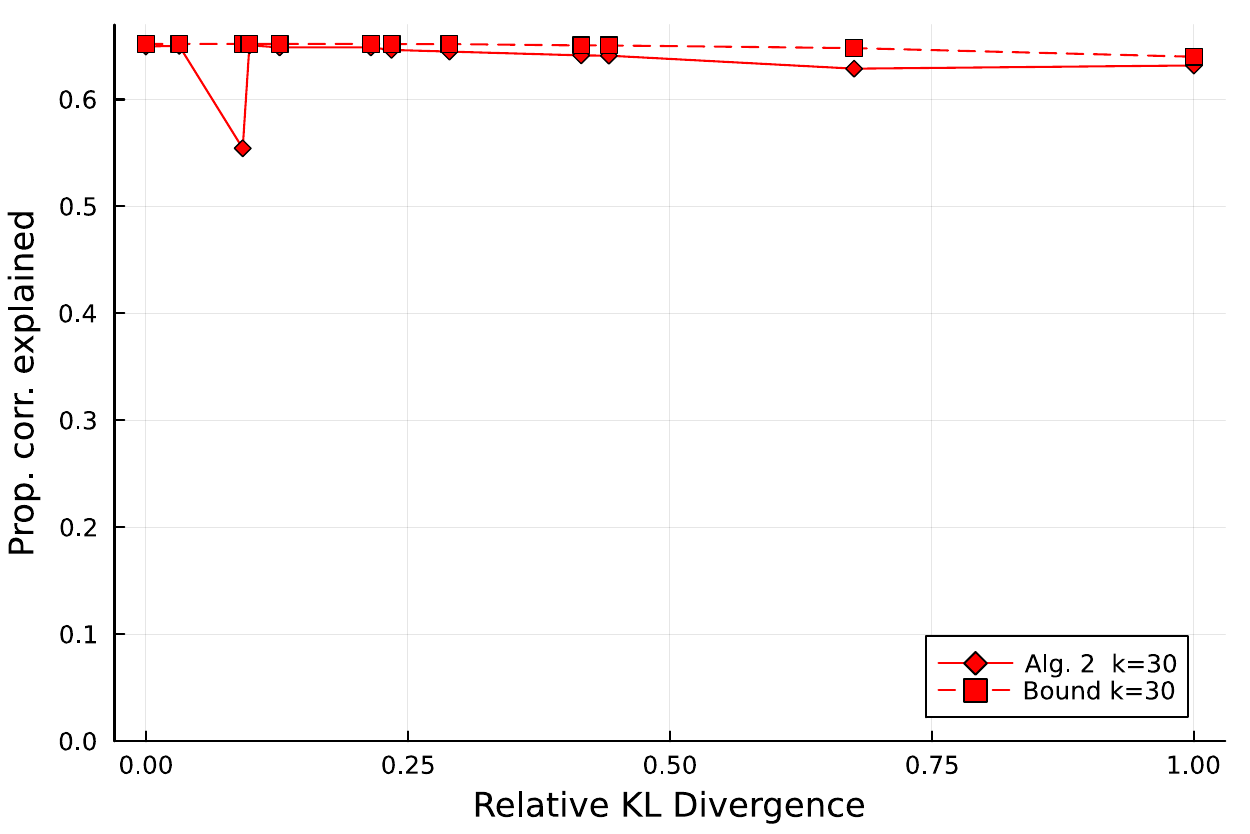}
    \end{subfigure}
        \begin{subfigure}[t]{.45\linewidth}
            \includegraphics[scale=0.35]{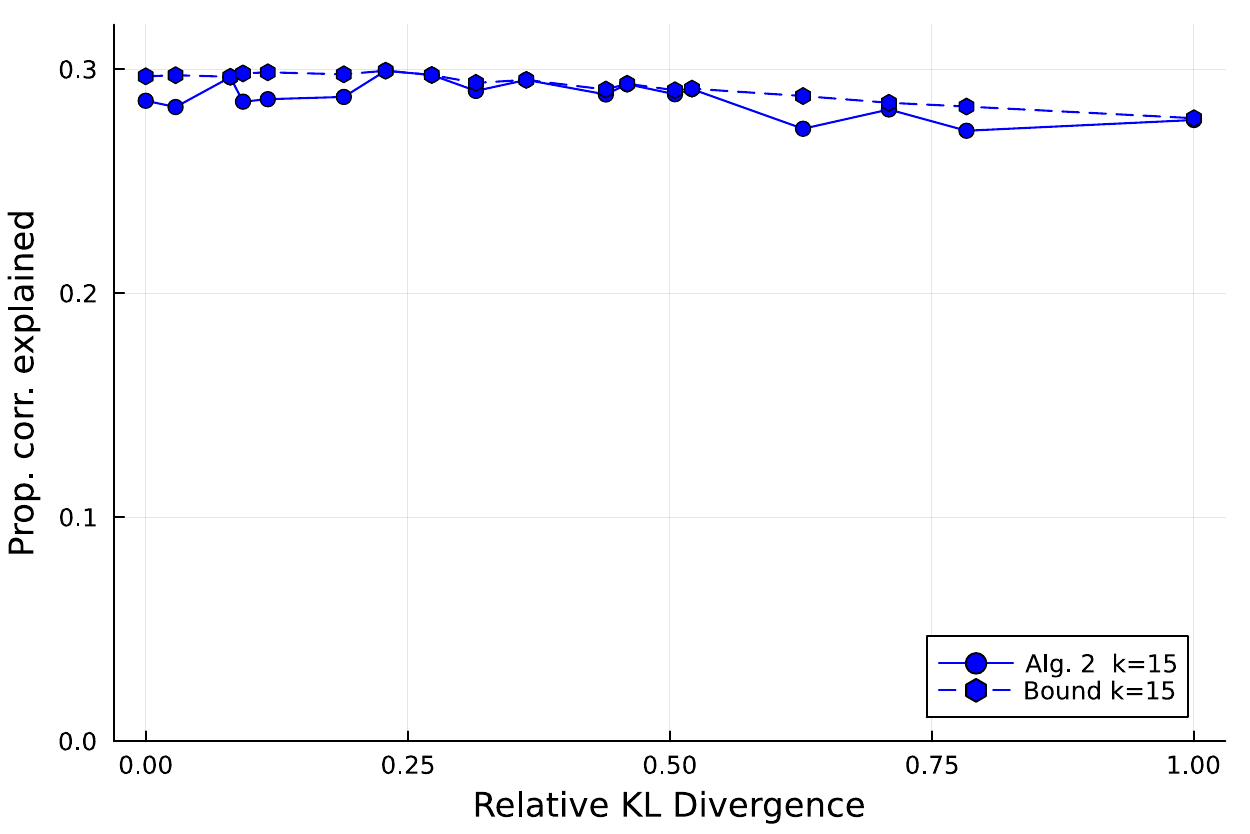}
    \end{subfigure}\\
        \begin{subfigure}[t]{.45\linewidth}
            \includegraphics[scale=0.35]{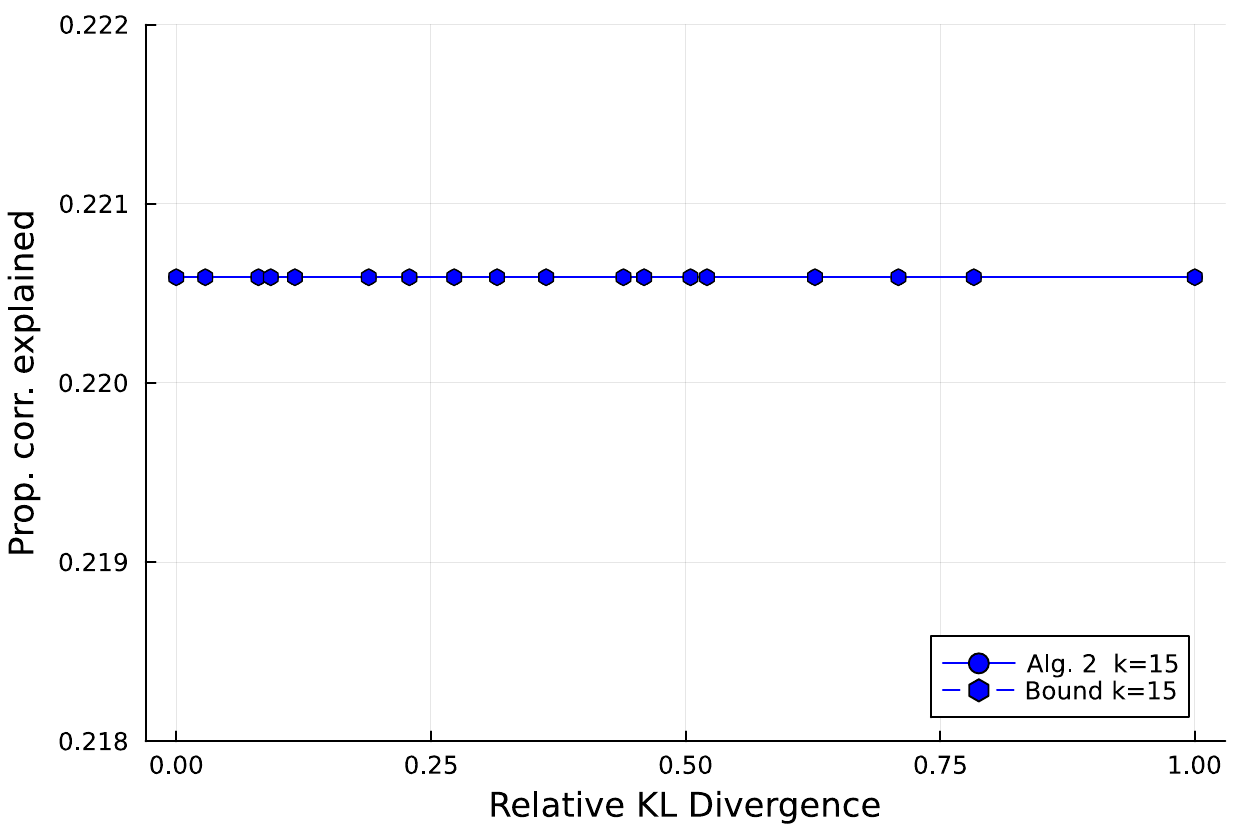}
    \end{subfigure}
    \begin{subfigure}[t]{.45\linewidth}
            \includegraphics[scale=0.35]{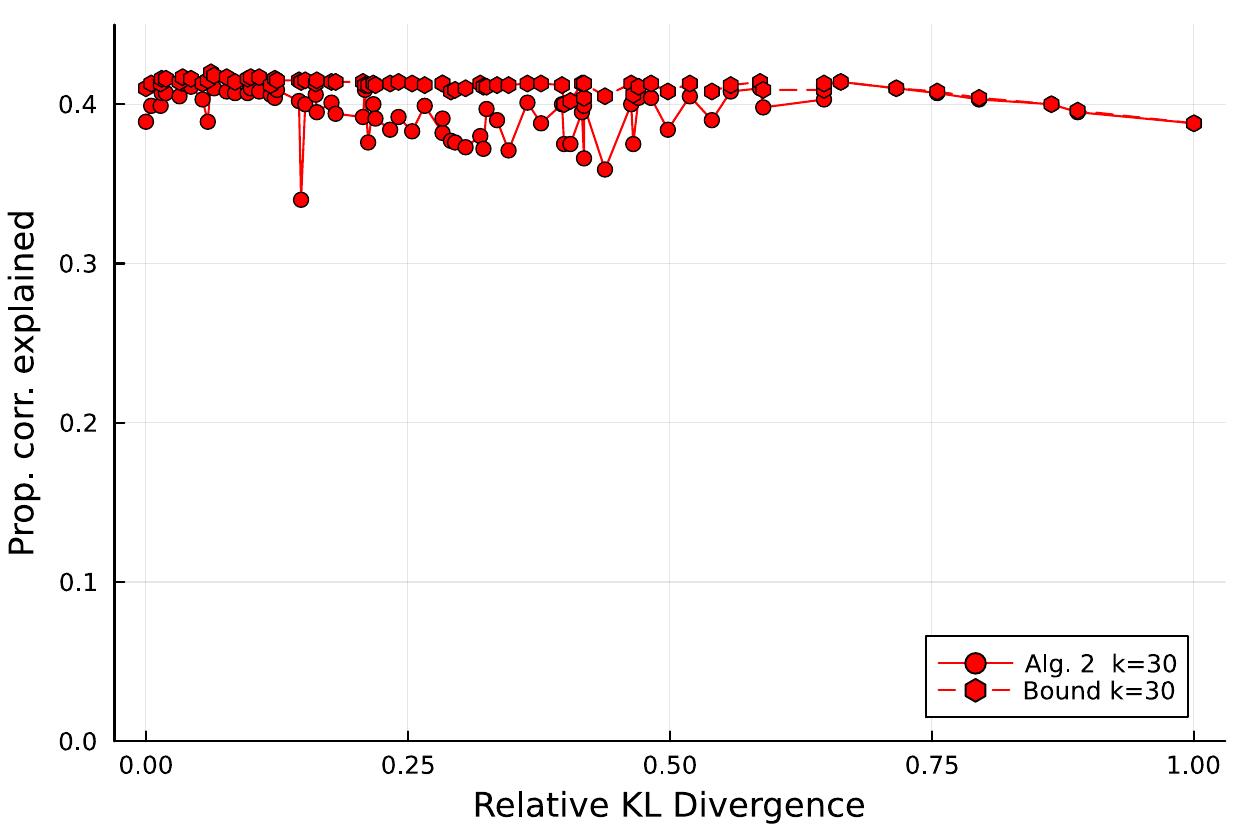}
    \end{subfigure}\\
    \begin{subfigure}[t]{.45\linewidth}
            \includegraphics[scale=0.35]{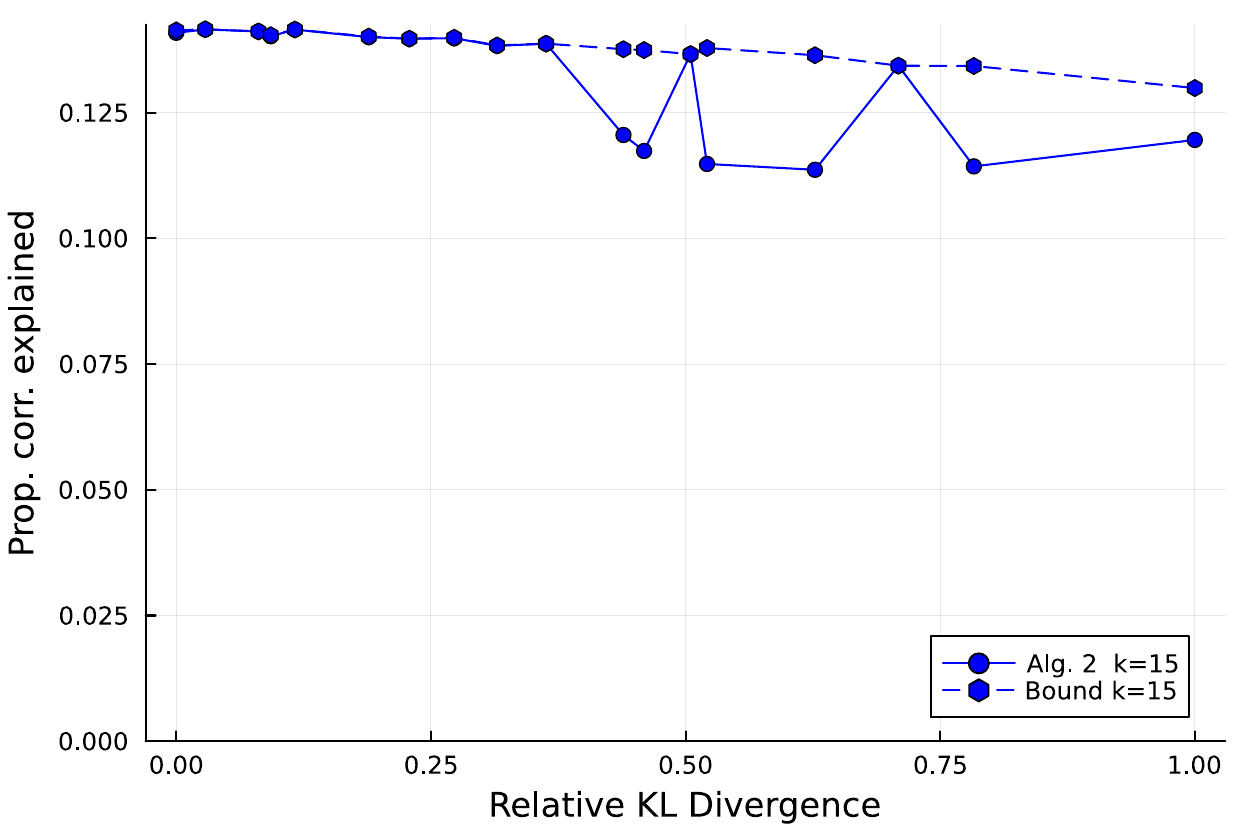}
    \end{subfigure}
        \begin{subfigure}[t]{.45\linewidth}
            \includegraphics[scale=0.35]{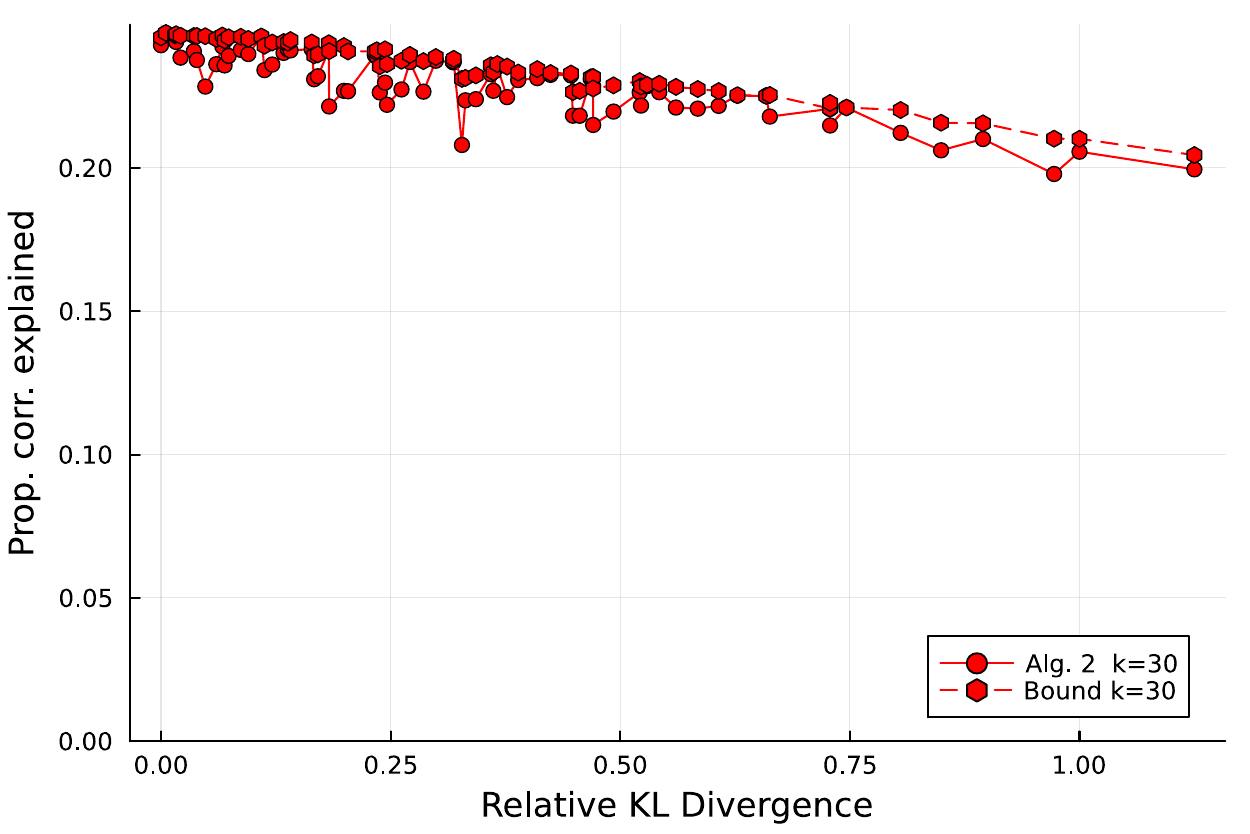}
    \end{subfigure}\\
   \caption{Symmetry of sparsity budget allocation vs. proportion of correlation in the dataset explained for pitprops $k=30$ (top left), ionosphere $k=15$ (top right), geographical $k=15$ (middle left), geographical $k=30$ (middle right), communities $k=15$ (bottom left), and communities $k=30$ (bottom right). Note that we normalize the KL divergence for $k=15$ and $k=30$ separately.} 
   \label{fig:sparsitysymmetry2}
\end{figure}

We consider the \verb|pitprops|, \verb|ionosphere|, \verb|geographical|, and \verb|communities| UCI datasets with a fixed number of PCs $r=3$ and a given overall sparsity budget $k \in \{15, 30\}$. Accordingly, in Figures \ref{fig:sparsitysymmetry} and \ref{fig:sparsitysymmetry2} of Section \ref{ssec:symmetryplots_appended}, we depict the relationship between the proportion of correlation explained in the data for each possible allocation of the sparsity budget $(k_1, k_2, k_3): k_1+k_2+k_3=k, p \geq k_1 \geq k_2 \geq k_3 \geq 1$, as computed by Algorithm \ref{alg:alternatingmin} with a limit of $200$ iterations and the same setup as in Section \ref{ssec:feasiblemethods} (and the corresponding upper bound computed by Algorithm \ref{alg:greedymethod2}), against the relative asymmetry in the sparsity budget, as measured by $$\frac{\mathrm{KL}((k_1, k_2, k_3) || (k/3, k/3, k/3))}{\max_{p \geq k_1 \geq k_2 \geq k_3: k_1+k_2+k_3=k}\mathrm{KL}((k_1, k_2, k_3) || (k/3, k/3, k/3))},$$ where $KL(p || q):=\sum_i p_i \log(p_i/q_i)$ denotes the KL divergence.

We observe a general trend that more symmetric sparsity budget allocations tend to explain more of the correlation in the data (both in terms of actual correlation explained by a solution from Algorithm \ref{alg:alternatingmin} and in terms of the upper bound). 
This suggests that, when time is a concern, requiring that all PCs are equally sparse is a reasonable approach. 

Table \ref{tab:asymmetrysummary} compares the quality of Algorithm \ref{alg:alternatingmin}'s solution (a) when all PCs have the sparsity budget of $k/r$ and (b) the maximum possible correlation explained over all feasible allocations of the sparsity budget $k$ (computed by enumerating all possible allocations of the sparsity budget-$18$ such allocations for $k=15, p \geq k$ and $74$ allocations for $k=30, p \geq 30$), together with the upper bound on the proportion of correlation explained obtained in each case. We observe that in several instances a perfectly symmetric allocation of the sparsity budget yields the highest quality solution, and in all instances a perfectly symmetric allocation is within $7\%$ in the worst-case (and within $1.75\%$ {\color{black}on} average) of the best solution. In the enumerated case, we also compute the optimality gap between the worst-case upper bound over all sparsity budget allocations, and the best solution found, and observe that on average it is less than $1\%$ over the instances considered. Note that this is a different gap to the one reported in Section \ref{ssec:feasiblemethods}, where the upper bound is computed after assuming that all PCs are equally sparse. 

\begin{table}[h]
\centering\footnotesize
\begin{tabular}{@{}l r r r r r r r r r r r r@{}} \toprule
Dataset & $p$ & $r$ & $k$ & \multicolumn{3}{c}{Symmetric} & \multicolumn{5}{c}{Enumerated} & Improvement ($\%$)\\
\cmidrule(l){5-7} \cmidrule(l){8-12} &   &                   & $k$  & UB & Obj. & Viol. & UB & $k_t$ & Obj. & Rel. gap $(\%)$ & Viol. &\\\midrule
Pitprops & 13 & 3 &	15 &	0.616 &	0.590 &	0.000 & 	0.618 & (6, 6, 3) &	0.593 &	$4.06\%$ & 0.000 &	$0.47\%$\\
 &	& & 	30 &	0.652 &	0.650 &	0.000	 & 0.652 &	(10, 10, 10) &	0.650	&$0.19\%$ & 0.000 &	$0\%$\\
 Ionosphere & 34 & 3 & 15 & 0.297	& 0.286 & 0.000 & 0.299 & (7, 6, 2) & 0.299 & $0\%$ & 0.001 & $4.48\%$\\
  &  &  & 30 & 0.411	& 0.400 & 0.000 & 0.412 & (15, 8, 7) & 0.402 & $2.34\%$ & 0.000 & $0.60\%$\\
 Geographical & 68 & 3 & 15 & 0.221	& 0.221 & 0.000 &  0.221 & (5, 5, 5) & 0.221 & $0\%$ & 0.000 & $0.00\%$\\
  &  &  & 30 & 0.410	& 0.389 & 0.000 & 0.420 & (12, 12, 6) & 0.415 & $1.07\%$ & 0.000 & $6.29\%$\\
 Communities & 101 & 3 & 15 & 0.141 & 0.141 & 0.000 & 0.142 & (6, 5, 4) & 0.142 & $0.02\%$ & 0.000 & $0.48\%$\\
 & & & 30 & 0.246 & 0.243 & 0.000 & 0.247 & (11, 10, 9) & 0.247 & $0.02\%$ & 0.000 & $1.71\%$\\
\bottomrule
\end{tabular}
\caption{Comparison of symmetric and enumerated solution computed by Algorithm \ref{alg:alternatingmin} for a given sparsity budget $k_{\text{total}}$ and given number of PCs $(r=3)$. We report the largest upper bound over all sparsity budget allocations as our enumerated upper bound, and report the relative optimality gap between the best solution and the worst-case bound. On average, considering asymmetric sparsity budget allocations improves the proportion of correlation explained by $1.65\%$.}
\label{tab:asymmetrysummary}
\end{table}

In summary, more symmetric allocations of the sparsity budget tend to perform better on average. Therefore, for a given sparsity budget $k$, a reasonable strategy could be to ($a$) run Algorithm \ref{alg:alternatingmin} to compute a perfectly symmetric allocation, ($b$) compute an upper bound using \eqref{prob:disjunctiverelax_permutationinvariant} across all possible allocations, and ($c$) run Algorithm \ref{alg:alternatingmin} only on the asymmetric allocations for which the upper bound from ($b$) allows for a significant potential improvement upon the symmetric solution.


\subsection{Summary and Guidelines From Numerical Experiments}

In summary, our main findings from our numerical experiments are as follows:
\begin{itemize}
\item As reflected in Section \ref{ssec:bounds},
the semidefinite relaxation \eqref{prob:disjunctiverelax_permutationinvariant} provides the tightest upper bound and can currently be computed within minutes for $p \leq 300$. At larger scales, the Lagrangian bound we develop in Algorithm \ref{alg:alternatingmin} is competitive when $kr \ll p$ and can typically be computed in seconds even when {\color{black}$p$ is in the thousands}.

\item In terms of solution quality, Algorithms \ref{alg:alternatingmin}--\ref{alg:disjoint.linalg} substantially outperform existing methods, obtaining higher-quality solutions that satisfy the orthogonality constraints (Section \ref{ssec:feasiblemethods}). Moreover, both algorithms require $100$ seconds or less to run, on average. 
Therefore, they should be considered a{\color{black}s} viable and accurate approach{\color{black}es} for sparse PCA problems with multiple PCs.

\item As demonstrated in Section \ref{ssec:feasiblemethods}, Algorithm \ref{alg:disjoint.linalg} gives certifiably near-optimal sparse PCs when $kr \leq p$ ($3.11\%$ bound gap across UCI instances on average), while Algorithm \ref{alg:alternatingmin} gives certifiably near-optimal sparse PCs when $kr > p$ ($2.82\%$ bound gap across UCI instances on average). Thus, using the techniques in the paper, it is possible to compute sparse PCs that are near-optimal in a practical amount of time, even when $p$ {\color{black}is in the thousands}. 

\item However, we should note that these optimality gap{\color{black}s} are computed using the tightest (i.e., typically Algorithm \ref{alg:greedymethod2}'s) upper bound. Using its own upper bound instead, Algorithm \ref{alg:disjoint.linalg} (resp. Algorithm \ref{alg:alternatingmin}) returns solutions with an average bound gap of 11.32\% (resp. 28.21\%) on datasets where $kr \leq p$ (resp. $kr > p$). So, while Algorithms \ref{alg:alternatingmin}--\ref{alg:disjoint.linalg} find these solutions in minutes when $p=1000$s, it currently takes hours to compute the conic bound and certify their near optimality. This observation emphasizes the benefit of combining these different approaches together, and motivates future research to improve the scalability of the conic bounds.

\item If practitioners have an overall sparsity budget but are agnostic about the sparsity of each column, a reasonable strategy is to require that all columns are equally sparse (i.e., set $k_t=k/r$), as shown in Section \ref{ssec:symm}. {\color{black}In our experiments, considering asymmetrically sparse sets of PCs increases the amount of {\color{black}correlation} explained by around $2\%$ on average, at the cost of an order-of-magnitude increase in total runtime.}
\end{itemize}

\section{Conclusion}
In this paper, we studied the problem of selecting a set of mutually orthogonal sparse principal components and proposed three algorithms which, collectively, for the first time, allow this problem to be solved to certifiable near-optimality with $100$s or $1000$s of features in minutes or hours. We propose a semidefinite relaxation (Section \ref{sec:relaxations}) which generates high-quality upper bounds on the amount of variance explainable by any set of sparse and mutually orthogonal components, and round this relaxation to obtain a feasible set of sparse PCs in Algorithm \ref{alg:greedymethod2}. Further, we propose a Lagrangian alternating maximization scheme (Algorithm \ref{alg:alternatingmin} in Section \ref{sec:deflation}) which gives an alternative upper bound and high-quality feasible solutions. Finally, we derived a new combinatorial upper bound on sparse PCA with multiple PCs, and used this bound to design a new relax-and-round scheme (Algorithm \ref{alg:disjoint.linalg} in Section \ref{sec.gen.gershgorin}) which also gives a valid upper bound and high-quality feasible solutions. Across suites of numerical experiments (Section \ref{sec:numres}), we demonstrate that the best solution obtained by any of our methods substantially outperforms methods from the literature in terms of obtaining sparse and mutually orthogonal PCs that explain most of the variance in a dataset, and nearly matches the tightest upper bound obtained by any of our methods. All in all, we solve sparse PCA with multiple components to near-optimality for $p=1000$s and $r>1$ for the first time. 

\subsection*{Acknowledgments}
We are very grateful to the area editor, Samuel Burer, for several valuable comments that helped us reorganize and clarify the manuscript's contributions. We are also grateful to the associate editor and the three anonymous referees for their valuable suggestions, which improved the manuscript.

\theendnotes
{
\scriptsize
\bibliographystyle{informs2014}

\begin{thebibliography}{70}
    \providecommand{\natexlab}[1]{#1}
    \providecommand{\url}[1]{\texttt{#1}}
    \providecommand{\urlprefix}{URL }
    
    \bibitem[{Ahmadi et~al.(2017)Ahmadi, Dash, \protect\BIBand{}
      Hall}]{ahmadi2017optimization}
    Ahmadi AA, Dash S, Hall G (2017) Optimization over structured subsets of
      positive semidefinite matrices via column generation. \emph{Discrete
      Optimization} 24:129--151.
    
    \bibitem[{Alizadeh \protect\BIBand{} Goldfarb(2003)}]{alizadeh2003second}
    Alizadeh F, Goldfarb D (2003) Second-order cone programming. \emph{Mathematical
      Programming} 95(1):3--51.
    
    \bibitem[{Amini \protect\BIBand{} Wainwright(2008)}]{amini2008high}
    Amini AA, Wainwright MJ (2008) High-dimensional analysis of semidefinite
      relaxations for sparse principal components. \emph{2008 IEEE International
      Symposium on Information Theory}, 2454--2458 (IEEE).
    
    \bibitem[{Asteris et~al.(2015)Asteris, Papailiopoulos, Kyrillidis,
      \protect\BIBand{} Dimakis}]{asteris2015sparse}
    Asteris M, Papailiopoulos D, Kyrillidis A, Dimakis AG (2015) Sparse {PCA} via
      bipartite matchings. \emph{Advances in Neural Information Processing Systems}
      28.
    
    \bibitem[{Atamt{\"u}rk \protect\BIBand{} Gomez(2025)}]{atamturk2019rank}
    Atamt{\"u}rk A, Gomez A (2025) Rank-one convexification for sparse regression.
      \emph{Journal of Machine Learning Research} 26(35):1--50.
    
    \bibitem[{Avellaneda \protect\BIBand{} Lee(2010)}]{avellaneda2010statistical}
    Avellaneda M, Lee JH (2010) Statistical arbitrage in the {US} equities market.
      \emph{Quantitative Finance} 10(7):761--782.
    
    \bibitem[{Barker \protect\BIBand{} Carlson(1975)}]{barker1975cones}
    Barker G, Carlson D (1975) Cones of diagonally dominant matrices. \emph{Pacific
      Journal of Mathematics} 57(1):15--32.
    
    \bibitem[{Behdin \protect\BIBand{} Mazumder(2025)}]{behdin2021sparse}
    Behdin K, Mazumder R (2025) Sparse {PCA}: A new scalable estimator based on
      integer programming. \emph{Annals of Statistics} .
    
    \bibitem[{Benidis et~al.(2016)Benidis, Sun, Babu, \protect\BIBand{}
      Palomar}]{benidis2016orthogonal}
    Benidis K, Sun Y, Babu P, Palomar DP (2016) Orthogonal sparse {PCA} and
      covariance estimation via procrustes reformulation. \emph{IEEE Transactions
      on Signal Processing} 64(23):6211--6226.
    
    \bibitem[{Berk \protect\BIBand{} Bertsimas(2019)}]{berk2019certifiably}
    Berk L, Bertsimas D (2019) Certifiably optimal sparse principal component
      analysis. \emph{Mathematical Programming Computation} 11(3):381--420.
    
    \bibitem[{Berthet \protect\BIBand{} Rigollet(2013)}]{berthet2013optimal}
    Berthet Q, Rigollet P (2013) Optimal detection of sparse principal components
      in high dimension. \emph{The Annals of Statistics} 41(4):1780--1815.
    
    \bibitem[{Bertsekas(1996)}]{bertsekas1996constrained}
    Bertsekas DP (1996) \emph{Constrained Optimization and Lagrange Multiplier
      Methods} (Athena Scientific).
    
    \bibitem[{Bertsimas \protect\BIBand{}
      Cory-Wright(2020)}]{bertsimas2020polyhedral}
    Bertsimas D, Cory-Wright R (2020) On polyhedral and second-order cone
      decompositions of semidefinite optimization problems. \emph{Operations
      Research Letters} 48(1):78--85.
    
    \bibitem[{Bertsimas et~al.(2021)Bertsimas, Cory-Wright, \protect\BIBand{}
      Pauphilet}]{bertsimas2019unified}
    Bertsimas D, Cory-Wright R, Pauphilet J (2021) A unified approach to
      mixed-integer optimization problems with logical constraints. \emph{SIAM
      Journal on Optimization} 31(3):2340--2367.
    
    \bibitem[{Bertsimas et~al.(2022{\natexlab{a}})Bertsimas, Cory-Wright,
      \protect\BIBand{} Pauphilet}]{bertsimas2020mixed}
    Bertsimas D, Cory-Wright R, Pauphilet J (2022{\natexlab{a}}) Mixed-projection
      conic optimization: A new paradigm for modeling rank constraints.
      \emph{Operations Research} 70(6):3321--3344.
    
    \bibitem[{Bertsimas et~al.(2022{\natexlab{b}})Bertsimas, Cory-Wright,
      \protect\BIBand{} Pauphilet}]{bertsimas2020solving}
    Bertsimas D, Cory-Wright R, Pauphilet J (2022{\natexlab{b}}) Solving
      large-scale sparse {PCA} to certifiable (near) optimality. \emph{Journal of
      Machine Learning Research} 23(13):1--35.
    
    \bibitem[{Bertsimas \protect\BIBand{} Kitane(2023)}]{bertsimassparse2022}
    Bertsimas D, Kitane DL (2023) Sparse {PCA}: A geometric approach. \emph{Journal
      of Machine Learning Research} 24:32--1.
    
    \bibitem[{Bienstock et~al.(2023)Bienstock, Pia, \protect\BIBand{}
      Hildebrand}]{bienstock2023complexity}
    Bienstock D, Pia AD, Hildebrand R (2023) Complexity, exactness, and rationality
      in polynomial optimization. \emph{Mathematical Programming} 197(2):661--692.
    
    \bibitem[{Boutsidis et~al.(2011)Boutsidis, Drineas, \protect\BIBand{}
      Magdon-Ismail}]{boutsidis2011sparse}
    Boutsidis C, Drineas P, Magdon-Ismail M (2011) Sparse features for {PCA}-like
      linear regression. \emph{Advances in Neural Information Processing Systems}
      24.
    
    \bibitem[{Bresler et~al.(2018)Bresler, Park, \protect\BIBand{}
      Persu}]{bresler2018sparse}
    Bresler G, Park SM, Persu M (2018) Sparse {PCA} from sparse linear regression.
      \emph{Advances in Neural Information Processing Systems} 31.
    
    \bibitem[{B{\"u}hler(2014)}]{buhler2014flexible}
    B{\"u}hler T (2014) \emph{A flexible framework for solving constrained ratio
      problems in machine learning}. Ph.D. thesis, Saarland University.
    
    \bibitem[{d'Aspremont et~al.(2008)d'Aspremont, Bach, \protect\BIBand{}
      El~Ghaoui}]{d2008optimal}
    d'Aspremont A, Bach F, El~Ghaoui L (2008) Optimal solutions for sparse
      principal component analysis. \emph{Journal of Machine Learning Research}
      9(7).
    
    \bibitem[{d'Aspremont et~al.(2007)d'Aspremont, El~Ghaoui, Jordan,
      \protect\BIBand{} Lanckriet}]{d2007direct}
    d'Aspremont A, El~Ghaoui L, Jordan MI, Lanckriet GR (2007) A direct formulation
      for sparse {PCA} using semidefinite programming. \emph{SIAM Review}
      49(3):434--448.
    
    \bibitem[{Del~Pia(2023)}]{del2022sparse}
    Del~Pia A (2023) Sparse {PCA} on fixed-rank matrices. \emph{Mathematical
      Programming} 198(1):139--157.
    
    \bibitem[{Deshpande \protect\BIBand{}
      Montanari(2014{\natexlab{a}})}]{deshpande2014information}
    Deshpande Y, Montanari A (2014{\natexlab{a}}) Information-theoretically optimal
      sparse {PCA}. \emph{2014 IEEE International Symposium on Information Theory},
      2197--2201 (IEEE).
    
    \bibitem[{Deshpande \protect\BIBand{}
      Montanari(2014{\natexlab{b}})}]{deshpande2014sparse}
    Deshpande Y, Montanari A (2014{\natexlab{b}}) Sparse {PCA} via covariance
      thresholding. \emph{Advances in Neural Information Processing Systems} 27.
    
    \bibitem[{Dey et~al.(2022{\natexlab{a}})Dey, Mazumder, \protect\BIBand{}
      Wang}]{dey2021using}
    Dey SS, Mazumder R, Wang G (2022{\natexlab{a}}) Using $\ell_1$-relaxation and
      integer programming to obtain dual bounds for sparse {PCA}. \emph{Operations
      Research} 70(3):1914--1932.
    
    \bibitem[{Dey et~al.(2022{\natexlab{b}})Dey, Molinaro, \protect\BIBand{}
      Wang}]{dey2020solving}
    Dey SS, Molinaro M, Wang G (2022{\natexlab{b}}) Solving sparse principal
      component analysis with global support. \emph{Mathematical Programming}
      199:421--459.
    
    \bibitem[{Ding et~al.(2024)Ding, Kunisky, Wein, \protect\BIBand{}
      Bandeira}]{ding2023subexponential}
    Ding Y, Kunisky D, Wein AS, Bandeira AS (2024) Subexponential-time algorithms
      for sparse {PCA}. \emph{Foundations of Computational Mathematics}
      24:865--914.
    
    \bibitem[{Eckart \protect\BIBand{} Young(1936)}]{eckart1936approximation}
    Eckart C, Young G (1936) The approximation of one matrix by another of lower
      rank. \emph{Psychometrika} 1(3):211--218.
    
    \bibitem[{Fan et~al.(2016)Fan, Liao, \protect\BIBand{} Wang}]{fan2016projected}
    Fan J, Liao Y, Wang W (2016) Projected principal component analysis in factor
      models. \emph{Annals of Statistics} 44(1):219.
    
    \bibitem[{Fisher(1981)}]{fisher1981lagrangian}
    Fisher ML (1981) The lagrangian relaxation method for solving integer
      programming problems. \emph{Management Science} 27(1):1--18.
    
    \bibitem[{Gally \protect\BIBand{} Pfetsch(2016)}]{gally2016computing}
    Gally T, Pfetsch ME (2016) Computing restricted isometry constants via
      mixed-integer semidefinite programming. \emph{Optimization Online} .
    
    \bibitem[{Garey \protect\BIBand{} Johnson(1979)}]{garey1979computers}
    Garey MR, Johnson DS (1979) \emph{Computers and Intractability: A Guide to the
      Theory of NP-completeness} (W.H.Freeman \& Co Ltd).
    
    \bibitem[{G{\"u}nl{\"u}k \protect\BIBand{}
      Linderoth(2010)}]{gunluk2010perspective}
    G{\"u}nl{\"u}k O, Linderoth J (2010) Perspective reformulations of mixed
      integer nonlinear programs with indicator variables. \emph{Mathematical
      Programming} 124(1):183--205.
    
    \bibitem[{Gupta et~al.(2023)Gupta, Van~Parys, \protect\BIBand{}
      Ryu}]{gupta2022branch}
    Gupta SD, Van~Parys BP, Ryu EK (2023) Branch-and-bound performance estimation
      programming: A unified methodology for constructing optimal optimization
      methods. \emph{Mathematical Programming} .
    
    \bibitem[{Hein \protect\BIBand{} B{\"u}hler(2010)}]{hein2010inverse}
    Hein M, B{\"u}hler T (2010) An inverse power method for nonlinear eigenproblems
      with applications in 1-spectral clustering and sparse {PCA}. \emph{Advances
      in Neural Information Processing Systems} 23.
    
    \bibitem[{Hopcroft \protect\BIBand{} Karp(1973)}]{hopcroft1973n}
    Hopcroft JE, Karp RM (1973) An n\^{}5/2 algorithm for maximum matchings in
      bipartite graphs. \emph{SIAM Journal on Computing} 2(4):225--231.
    
    \bibitem[{Horn \protect\BIBand{} Johnson(1985)}]{johnson1985matrix}
    Horn RA, Johnson CR (1985) \emph{Matrix analysis} (Cambridge {U}niversity
      {P}ress, {N}ew {Y}ork).
    
    \bibitem[{Hotelling(1933)}]{hotelling1933analysis}
    Hotelling H (1933) Analysis of a complex of statistical variables into
      principal components. \emph{Journal of Educational Psychology} 24(6):417.
    
    \bibitem[{Jeffers(1967)}]{jeffers1967two}
    Jeffers JN (1967) Two case studies in the application of principal component
      analysis. \emph{Journal of the Royal Statistical Society: Series C (Applied
      Statistics)} 16(3):225--236.
    
    \bibitem[{Johnstone \protect\BIBand{} Lu(2009)}]{johnstone2009consistency}
    Johnstone IM, Lu AY (2009) On consistency and sparsity for principal components
      analysis in high dimensions. \emph{Journal of the American Statistical
      Association} 104(486):682--693.
    
    \bibitem[{Jolliffe et~al.(2003)Jolliffe, Trendafilov, \protect\BIBand{}
      Uddin}]{jolliffe2003modified}
    Jolliffe IT, Trendafilov NT, Uddin M (2003) A modified principal component
      technique based on the {LASSO}. \emph{Journal of Computational and Graphical
      Statistics} 12(3):531--547.
    
    \bibitem[{Journ{\'e}e et~al.(2010)Journ{\'e}e, Nesterov, Richt{\'a}rik,
      \protect\BIBand{} Sepulchre}]{journee2010generalized}
    Journ{\'e}e M, Nesterov Y, Richt{\'a}rik P, Sepulchre R (2010) Generalized
      power method for sparse principal component analysis. \emph{Journal of
      Machine Learning Research} 11(2).
    
    \bibitem[{Kim et~al.(2022)Kim, Tawarmalani, \protect\BIBand{}
      Richard}]{kim2021convexification}
    Kim J, Tawarmalani M, Richard JPP (2022) Convexification of
      permutation-invariant sets and an application to sparse principal component
      analysis. \emph{Mathematics of Operations Research} 47(4):2547--2584.
    
    \bibitem[{Krauthgamer et~al.(2015)Krauthgamer, Nadler, \protect\BIBand{}
      Vilenchik}]{krauthgamer2015semidefinite}
    Krauthgamer R, Nadler B, Vilenchik D (2015) Do semidefinite relaxations solve
      sparse pca up to the information limit? \emph{The Annals of Statistics}
      43(3):1300--1322.
    
    \bibitem[{Li \protect\BIBand{} Xie(2024)}]{li2021beyond}
    Li Y, Xie W (2024) Beyond symmetry: Best submatrix selection for the sparse
      truncated {SVD}. \emph{Mathematical Programming} 208:1--50.
    
    \bibitem[{Li \protect\BIBand{} Xie(2025)}]{li2020exact}
    Li Y, Xie W (2025) Exact and approximation algorithms for sparse {PCA}.
      \emph{INFORMS Journal on Computing} 37(3):582--602.
    
    \bibitem[{Lu \protect\BIBand{} Zhang(2012)}]{lu2012augmented}
    Lu Z, Zhang Y (2012) An augmented {L}agrangian approach for sparse principal
      component analysis. \emph{Mathematical Programming} 135(1):149--193.
    
    \bibitem[{Lubin et~al.(2022)Lubin, Vielma, \protect\BIBand{}
      Zadik}]{lubin2022mixed}
    Lubin M, Vielma JP, Zadik I (2022) Mixed-integer convex representability.
      \emph{Mathematics of Operations Research} 47(1):720--749.
    
    \bibitem[{Mackey(2008)}]{mackey2008deflation}
    Mackey L (2008) Deflation methods for sparse {PCA}. \emph{Advances in Neural
      Information Processing Systems} 21.
    
    \bibitem[{Marshall \protect\BIBand{} Olkin(1979)}]{marshall1979inequalities}
    Marshall AW, Olkin I (1979) \emph{Inequalities: Theory of Majorization and its
      Applications} (Academic, New York).
    
    \bibitem[{Naikal et~al.(2011)Naikal, Yang, \protect\BIBand{}
      Sastry}]{naikal2011informative}
    Naikal N, Yang AY, Sastry SS (2011) Informative feature selection for object
      recognition via sparse {PCA}. \emph{2011 International Conference on Computer
      Vision}, 818--825 (IEEE).
    
    \bibitem[{Nocedal \protect\BIBand{} Wright(2006)}]{nocedal2006numerical}
    Nocedal J, Wright SJ (2006) \emph{Numerical Optimization} (Springer).
    
    \bibitem[{Overton \protect\BIBand{} Womersley(1992)}]{overton1992sum}
    Overton ML, Womersley RS (1992) On the sum of the largest eigenvalues of a
      symmetric matrix. \emph{SIAM Journal on Matrix Analysis and Applications}
      13(1):41--45.
    
    \bibitem[{Pearson(1901)}]{pearson1901liii}
    Pearson K (1901) Liii. on lines and planes of closest fit to systems of points
      in space. \emph{The London, Edinburgh, and Dublin Philosophical Magazine and
      Journal of Science} 2(11):559--572.
    
    \bibitem[{Probel \protect\BIBand{} Tropp(2011)}]{probel2011large}
    Probel CJ, Tropp JA (2011) Large-scale {PCA} with sparsity constraints.
      Technical report, California Institute of Technology.
    
    \bibitem[{Ramana(1997)}]{ramana1997exact}
    Ramana MV (1997) An exact duality theory for semidefinite programming and its
      complexity implications. \emph{Mathematical Programming} 77(1):129--162.
    
    \bibitem[{Reuther et~al.(2018)Reuther, Kepner, Byun, Samsi, Arcand, Bestor,
      Bergeron, Gadepally, Houle, Hubbell, Jones, Klein, Milechin, Mullen, Prout,
      Rosa, Yee, \protect\BIBand{} Michaleas}]{reuther2018interactive}
    Reuther A, Kepner J, Byun C, Samsi S, Arcand W, Bestor D, Bergeron B, Gadepally
      V, Houle M, Hubbell M, Jones M, Klein A, Milechin L, Mullen J, Prout A, Rosa
      A, Yee C, Michaleas P (2018) Interactive supercomputing on 40,000 cores for
      machine learning and data analysis. \emph{2018 IEEE High Performance extreme
      Computing Conference (HPEC)}, 1--6 (IEEE).
    
    \bibitem[{Rudin et~al.(2022)Rudin, Chen, Chen, Huang, Semenova,
      \protect\BIBand{} Zhong}]{rudin2022interpretable}
    Rudin C, Chen C, Chen Z, Huang H, Semenova L, Zhong C (2022) Interpretable
      machine learning: Fundamental principles and 10 grand challenges.
      \emph{Statistics Surveys} 16:1--85.
    
    \bibitem[{Tan et~al.(2014)Tan, Petersen, \protect\BIBand{}
      Witten}]{tan2014classification}
    Tan KM, Petersen A, Witten D (2014) Classification of {RNA}-seq data.
      \emph{Statistical Analysis of Next Generation Sequencing Data}, 219--246
      (Springer).
    
    \bibitem[{Tropp et~al.(2017)Tropp, Yurtsever, Udell, \protect\BIBand{}
      Cevher}]{tropp2017practical}
    Tropp JA, Yurtsever A, Udell M, Cevher V (2017) Practical sketching algorithms
      for low-rank matrix approximation. \emph{SIAM Journal on Matrix Analysis and
      Applications} 38(4):1454--1485.
    
    \bibitem[{Udell et~al.(2016)Udell, Horn, Zadeh, \protect\BIBand{}
      Boyd}]{udell2016generalized}
    Udell M, Horn C, Zadeh R, Boyd S (2016) Generalized low rank models.
      \emph{Foundations and Trends{\textregistered} in Machine Learning}
      9(1):1--118.
    
    \bibitem[{Vu et~al.(2013)Vu, Cho, Lei, \protect\BIBand{} Rohe}]{vu2013fantope}
    Vu VQ, Cho J, Lei J, Rohe K (2013) Fantope projection and selection: A
      near-optimal convex relaxation of sparse pca. \emph{Advances in Neural
      Information Processing Systems} 26.
    
    \bibitem[{Vu \protect\BIBand{} Lei(2013)}]{vu2013minimax}
    Vu VQ, Lei J (2013) Minimax sparse principal subspace estimation in high
      dimensions. \emph{The Annals of Statistics} 41(6):2905--2947.
    
    \bibitem[{Wang et~al.(2023)Wang, Dey, \protect\BIBand{} Xie}]{wang2023variable}
    Wang J, Dey SS, Xie Y (2023) Variable selection for kernel two-sample tests.
      \emph{arXiv preprint arXiv:2302.07415} .
    
    \bibitem[{Wei et~al.(2022)Wei, G{\'o}mez, \protect\BIBand{}
      K{\"u}{\c{c}}{\"u}kyavuz}]{wei2022ideal}
    Wei L, G{\'o}mez A, K{\"u}{\c{c}}{\"u}kyavuz S (2022) Ideal formulations for
      constrained convex optimization problems with indicator variables.
      \emph{Mathematical Programming} 192(1):57--88.
    
    \bibitem[{Witten et~al.(2009)Witten, Tibshirani, \protect\BIBand{}
      Hastie}]{witten2009penalized}
    Witten DM, Tibshirani R, Hastie T (2009) A penalized matrix decomposition, with
      applications to sparse principal components and canonical correlation
      analysis. \emph{Biostatistics} 10(3):515--534.
    
    \bibitem[{Yuan \protect\BIBand{} Zhang(2013)}]{yuan2013truncated}
    Yuan XT, Zhang T (2013) Truncated power method for sparse eigenvalue problems.
      \emph{Journal of Machine Learning Research} 14(4).
    
    \bibitem[{Zou et~al.(2006)Zou, Hastie, \protect\BIBand{}
      Tibshirani}]{zou2006sparse}
    Zou H, Hastie T, Tibshirani R (2006) Sparse principal component analysis.
      \emph{Journal of Computational and Graphical Statistics} 15(2):265--286.
    
    \end{thebibliography}

}

\newpage
\ECSwitch
\ECHead{Supplementary Material}

\section{Proof of Proposition \ref{prop:orth.to.rank}}\label{ec.orth.to.rank}
\proof{Proof of Proposition \ref{prop:orth.to.rank}} We decompose each matrix $\bm{Y}^t$ into $\bm{Y}^t = \bm{u}_t \bm{u}_t^\top$ with $\| \bm{u}_t \|^2 = \operatorname{tr}(\bm{Y}^t) = 1$. Hence, for any pair $(t,t')$, $\langle \bm{Y}^t, \bm{Y}^{t'} \rangle = \bm{u}_t^\top \bm{Y}^{t'} \bm{u}_t = (\bm{u}_t^\top \bm{u}_{t'})^2 \geq 0$.

$(\Rightarrow)$ If $\bm{Y} := \sum_{t' \in [r]} \bm{Y}^{t'} \preceq \mathbb{I}$, then, for any $t\in [r]$, $\bm{u}_t^\top \bm{Y} \bm{u}_t \leq \| \bm{u}_t \|^2 = 1$. However, $\bm{u}_t^\top \bm{Y} \bm{u}_t = 1 + \sum_{t' \neq t} \langle \bm{Y}^t, \bm{Y}^{t'}\rangle$. Hence, for all $t' \neq t$, we must have $\langle \bm{Y}^t, \bm{Y}^{t'}\rangle = 0$.

$(\Leftarrow)$ If $\langle \bm{Y}^t, \bm{Y}^{t'}\rangle = 0$ for all $t' \neq t$, then $\{ \bm{u}_t \}_{t \in [r]}$ is an orthonormal family that can be completed to form an orthonormal basis $\{ \bm{u}_t \}_{t \in [p]}$. For any $t \in [{\color{black}r}], t' \in [r]$, $\bm{u}_t^\top \bm{Y}^{t'} \bm{u}_t = 1$ if $t=t'$, $0$ otherwise so for any $t\in [{\color{black}r}]$, $\bm{u}_t^\top \bm{Y} \bm{u}_t \leq \| \bm{u}_t \|^2$ and $\bm{Y} = \sum_{t' \in [r]}  \bm{Y}^{t'} \preceq \mathbb{I}$.
\hfill\Halmos
\endproof

\section{Proof of Theorem \ref{thm:validinequalities}}\label{sec:ec.validineq}
In this section, we derive the valid inequalities introduced in Theorem \ref{thm:validinequalities} from first principles. We proceed in two ways. First, 
we derive {\color{black}inequalities} which hold for each $\bm{Y}^t$ separately. Second, we observe that these inequalities can be generalized to also apply for $\bm{Y}=\sum_{t \in [r]}\bm{Y}^t$, and also derive new inequalities which reflect the interaction of the sparsity and rank constraints. 

We repeatedly reference two results on the convex hulls of convex quadratic functions under logical constraints that were established by \cite{wei2022ideal}, building upon the work of \cite{gunluk2010perspective, atamturk2019rank}:

\begin{lemma}{\citep[Theorem 3 of][]{wei2022ideal}}\label{lemma:persp}
The convex closure of the set
\begin{align*}
    \mathcal{S}=&\left\{(\bm{x}, \bm{z}, t) \in \mathbb{R}^n \times \{0, 1\}^n \times \mathbb{R}: \quad t \geq \sum_{i=1}^n x_i^2,\ \bm{e}^\top \bm{z} \leq k, \ x_i=0 \ \text{if} \ z_i=0, \ \forall i \in [n]\right\} \\
\text{is given by}&\\
    \mathcal{S}^c=&\left\{(\bm{x}, \bm{z}, \bm{\theta}, t) \in \mathbb{R}^n \times [0, 1]^n \times \mathbb{R}^n \times \mathbb{R}: \quad t \geq \sum_{i=1}^n \theta_i,\ \bm{e}^\top \bm{z} \leq k,\ \theta_i z_i \geq x_i^2 \ \forall i \in [n]\right\}.
\end{align*}
\end{lemma}
The above result is sometimes known as a perspective reformulation, since we strengthen the quadratic constraint $t \geq \sum_i x_i^2$ by replacing $x_i^2$ with its perspective $z_i (x_i/z_i)^2$. 
\begin{lemma}{\citep[Proposition 4 of][]{wei2022ideal}}\label{lemma:rankone}
Let $k \geq 2$. Then, the convex closure of the set
\begin{align*}
    \mathcal{T}=&\left\{(\bm{x}, \bm{z}, t) \in \mathbb{R}^n \times \{0, 1\}^n \times \mathbb{R}: \quad t \geq (\bm{e}^\top \bm{x})^2,\ \bm{e}^\top \bm{z} \leq k, \ x_i=0 \ \text{if} \ z_i=0, \ \forall i \in [n]\right\},\\
\text{is given by} &\\
    \mathcal{T}^c=&\left\{(\bm{x}, \bm{z}, t) \in \mathbb{R}^n \times [0, 1]^n \times \mathbb{R}: \quad t \cdot \min(1, \bm{e}^\top \bm{z}) \geq (\bm{e}^\top \bm{x})^2,\ \bm{e}^\top \bm{z} \leq k\right\}. 
\end{align*}
\end{lemma}
Lemma \ref{lemma:rankone} 
is extremely useful when a single continuous variable depends upon multiple indicator variables, as occurs in certain substructures of our reformulations of Problem \eqref{prob:spcaorth_compact}.

\paragraph{Rank-One Valid Inequalities}
First, inspired by \cite{bertsimas2020polyhedral}, we observe that in a feasible solution to Problem \eqref{prob:disjunctiverelax} the $2 \times 2$ minors of $\bm{Y}^t$ are certainly non-negative, i.e.,
\begin{align*}
    (Y^t_{i,j})^2\leq Y^t_{i,i}Y^t_{j,j}, \ \forall i,j \in [p].
\end{align*}
These constraints are implied by $\bm{Y}^t \succeq \bm{0}$ and hence redundant in and of themselves. However, we can sum over all such constraints $i \in [p]$ and use $\mathrm{tr}(\bm{Y}^t)=1$, to obtain the constraint
\begin{align*}
    \sum_{i \in [p]}{(Y^t_{i,j})}^2 \leq Y^t_{j,j}, \ \forall j \in [p].
\end{align*}
This constraint is a sum of redundant constraints and hence redundant. However, we can strengthen it, by noting that it is a separable convex quadratic inequality under logical constraints. Indeed, by Lemma \ref{lemma:persp}, its convex closure under the logical constraints $Y_{i,j}^t=0$ if $Z_{i,t}=0$ is given by:
\begin{align}
    \sum_{j=1}^p {(Y^t_{i,j})}^2\leq Y^t_{i,i}Z_{i,t}, \ \forall i \in [p], t \in [{\color{black}r}].
\end{align}

\paragraph{Rank-$r$ Valid Inequalities}

{In the same spirit as in the rank one case, we can obtain strong valid inequalities by summing the $2 \times 2$ minors of $Y_{i,i}=\sum_{t=1}^r Y_{i,i}^t$. Indeed, since $\bm{Y}$ is positive semidefinite, summing its $2 \times 2$ minors implies that:
\begin{align*}
    \sum_{j=1}^p {(Y_{i,j})}^2\leq r Y_{i,i}.
\end{align*}
}Moreover, since $Y_{i,j}=\sum_{t=1}^r \bm{Y}^t_{i,j}$ is a rank-one quadratic under logical constraints $Y_{i,j}^t=0$ if $Z_{i,t}=0$, invoking Lemma \ref{lemma:rankone} reveals that the convex closure of this quadratic constraint under these logical constraints is given by the strengthened inequality:
\begin{align}
    \sum_{j=1}^p {(Y_{i,j})}^2\leq r Y_{i,i}\min\left(1, \sum_{t=1}^r  Z_{i,t}\right), \ \forall i \in [p].
\end{align}

Second, in any feasible solution we have:
\begin{align*}
    \vert Y_{i,j}\vert \leq \sum_{t=1}^r \vert Y^t_{i,j}\vert =\sum_{t=1}^r \vert U_{i,t}\vert \vert U_{j,t}\vert.
\end{align*}
Let us denote by $k_t$ the sparsity of the $t$th column of $\bm{U}$. Then, it is well known that $\Vert \bm{U}_t\Vert_1 \leq \sqrt{k_t}$. Therefore:
\begin{align*}
    \sum_{j=1}^p \vert Y_{i,j}\vert \leq\sum_{j=1}^p \left(\sum_{t=1}^r \vert U_{i,t}\vert U_{j,t}\vert \right)\leq \sum_{t=1}^r \sqrt{k_t}\vert U_{i,t}\vert.
\end{align*}

Next, squaring both sides and invoking the Cauchy-Schwarz inequality reveals that
\begin{align*}
    \left(\sum_{j=1}^p \vert Y_{i,j}\vert\right)^2 \leq \left(\sum_{t=1}^r U_{i,t}^2 \right) \left(\sum_{t=1}^r k_t \right)=k Y_{i,i}.
\end{align*}
Finally, noting that the expression $\left(\sum_{j=1}^p \vert Y_{i,j}\vert\right)^2 \leq k Y_{i,i}$ is a convex quadratic under logical constraints $ Y_{i,j}^t=0$ if $Z_{i,t}=0$ and invoking Lemma \ref{lemma:rankone} to obtain its convex closure yields the strengthened second-order cone inequality
\begin{align}
    \left(\sum_{j=1}^p \vert Y_{i,j}\vert\right)^2 \leq k Y_{i,i} \min\left(1, \sum_{t \in [r]} Z_{i,t}\right), \ \forall i \in [p], t \in [{\color{black}r}].
\end{align}

Third, in the same spirit, the $2 \times 2$ minors of $\bm{Y} \preceq \mathrm{Diag}\left(\min\left(\bm{e}, \sum_{t \in [r]} \bm{Z}_t\right)\right)$ are
\begin{align*}
    \left(\min\left(1, \sum_{t \in [r]} Z_{i,t}\right)-Y_{i,i}\right) \left( \min\left(1,\sum_{t \in [r]}Z_{j,t}\right)-Y_{j,j}\right)\geq {\left(\min\left(1,\sum_{t \in [r]}Z_{i,t}\right)\delta_{i,j}-Y_{i,j}\right)}^2,
\end{align*}
where $\delta_{i,j}=\bm{1}\{i=j\}$ is an indicator denoting whether $i=j$. Summing these constraints over all indices $i \neq j$ and using $k-r+1$ as an upper bound on $ \sum_{i \in [p]: i \neq j}\sum_{t \in [r]} Z_{i,t}-Y_{i,i}$ then yields 
\begin{align*}
    (k-r+1) \left( \min\left(1,\sum_{t \in [r]}Z_{j,t}\right)-Y_{j,j}\right)\geq \sum_{i \in [p]: i \neq j}Y_{i,j}^2, \ \forall j \in [p].
\end{align*}
Finally, we recognize the right hand side as a sum of rank-one quadratic terms $(\sum_{t=1}^r Y_{i,j}^t)^2$ under logical constraints $Y_{i,j}^t=0$ if $Z_{j,t}=0$ and invoke Lemma \ref{lemma:rankone} to obtain the convex closure, giving:
\begin{align}
    (k-r+1) \min\left(1,\sum_{t \in [r]}Z_{j,t}\right)( \min\left(1,\sum_{t \in [r]}Z_{j,t}\right)-Y_{j,j})\geq \sum_{i \in [p]: i \neq j}Y_{i,j}^2 \ \forall j \in [p].
\end{align}

The result then follows by introducing a vector $\bm{w}$ such that $w_i$ models $\min(1, \sum_{t \in [r]}Z_{i,t})$ via $\bm{w} 
\in [0, 1]^p$, $\bm{w} \leq \bm{Z}\bm{e}$, and noting that we can replace $\bm{Y} \preceq \mathrm{Diag}(\min(\bm{e}, \bm{Z}\bm{e})$ with $\bm{Y} \preceq \mathrm{Diag}(\bm{w})$.

\section{Proof of Proposition \ref{prop:validineqskt}}\label{ssec:append.proof_prop_validineqskt}

\proof{Proof of Proposition \ref{prop:validineqskt}}
First, let us observe that if $\sum_{i \in [p]} Z_{i,t} \leq k$ and $\bm{Y}^t=\bm{U}_t\bm{U}_t^\top$ is a rank-one matrix such that $\Vert \bm{U}\Vert_2=1$ then we have $\Vert \bm{U}\Vert_1 \leq \sqrt{k_t}$ by norm equivalence. Therefore
\begin{align*}
    \sum_{j=1}^p \vert Y_{i,j}^t \vert \leq \sum_{j=1}^p \vert U_{i,t} \vert \vert U_{j,t}\vert \leq \sqrt{k_t} \vert U_{i,t}\vert.
\end{align*}
Squaring both sides of this inequality then yields 
\begin{align*}
    \left(\sum_{j=1}^p \vert Y_{i,j}^t \vert\right)^2 \leq k_t Y_{i,i}^t,
\end{align*}
and combining Lemma \ref{lemma:persp} with this inequality yields \eqref{ell1_split}.

Second, in the same spirit, since $\bm{U}_t \bm{U}_t^\top$ is only supported on indices where $\bm{Z}_t$ is non-zero, we have that $\bm{Y}^t \preceq \mathrm{Diag}(\bm{Z}_t)$. This constraint implies the following $2 \times 2$ minors are non-negative
\begin{align*}
    (Z_{i,t}-Y_{i,i}^t) ( Z_{j,t}-Y_{j,j}^t)\geq {(\delta_{i,j}-Y_{i,j}^t)}^2, \ \forall i,j \in [p],
\end{align*}
where $\delta_{i,j}=1$ if $i=j$ and $0$ otherwise. Summing these inequalities over indices $i \neq j$ and setting $k_t-1$ as a valid upper bound on $\sum_{i \in [p]: i \neq j}Z_{i,t}-Y_{i,i}^t$ whenever $Z_{j,t}=1$ (as $Y_{i,j}^t=0$ if $Z_{j,t}=0$) gives
\begin{align*}
    (k_t-1) (Z_{j,t}-Y_{j,j}^t)\geq \sum_{i \in [p]: i \neq j}{Y_{i,j}^t}^2, \ \forall j \in [p].
\end{align*}
Finally, using Lemma \ref{lemma:persp} to take the convex closure of this inequality under the logical constraints $Y^t_{i,j}=0$ if $Z_{j,t}=0$ gives Equation \eqref{eqref:soc2}.
\hfill \Halmos
\endproof

\section{Complete Formulations for the Semidefinite and Second-Order Cone Relaxations}\label{sec:ec.formulations}
In this section, we provide the complete formulation for our SDP relaxation \eqref{prob:disjunctiverelax_permutationinvariant} as well as its second-order cone approximation.

\subsection{Semidefinite Relaxation} \label{ssec:ec.sdp}
In Section \ref{sec:validineqk}, we proposed a semidefinite relaxation, \eqref{prob:disjunctiverelax_permutationinvariant}, in the case where a sparsity budget for each PC, $k_t$, is provided. In particular, this relaxation involves valid inequalities that \cite{kim2021convexification} have derived in the single-PC case.
To the best of our knowledge, their formulation 
leads to the strongest known relaxation for sparse PCA with $r=1$ which can be solved in polynomial time. 
We note however that invoking a fixed but sufficiently large level of the sum-of-squares hierarchy may give tighter relaxations, although we do not write these relaxations down as they involve very large semidefinite constraints and are therefore intractable in practice (see also \citet{dey2021using} for an NP-hard relaxation that uses the $\ell_1$ norm).

In \eqref{prob:disjunctiverelax_permutationinvariant}, we concisely denoted $(\bm{Y}^t, {\bm{Z}_t}) \in \mathcal{T}(k_t)$ the set of valid inequalities involved in \cite{kim2021convexification}'s ``T-relaxation''. We now elicit the constraints involved in the set $\mathcal{T}(k_t)$ and provide the complete formulation of the SDP relaxation \eqref{prob:disjunctiverelax_permutationinvariant}. For each $t \in [r]$, we introduce an additional variable $\bm{F}^t$ to capture the entry-wise absolute value of $\bm{Y}^t$, and an additional matrix $\bm{G}^t$ which contains a sorted version of $\bm{F}^t$. We obtain:
\begin{align}
    \max_{\substack{\bm{Z} \in [0, 1]^{p \times r}:\\ \langle \bm{E}, \bm{Z}\rangle \leq k, \\ \bm{w} \in [0, 1]^p}
    }\max_{\substack{\bm{Y} \in \mathcal{S}^p_+, \bm{Y}^t, \bm{F}^t, \bm{G}^t \in \mathcal{S}^p_+,\\ \bm{T}^t  \in \mathbb{R}^{p \times p}_+,\\ \bm{r}^{t,D} \in \mathbb{R}^{p-1}, \bm{t}^{t,D} \in \mathbb{R}^{p \times p-1}_+}} \quad & \langle \bm{Y}, \bm{\Sigma}\rangle \label{prob:disjunctiverelax_permutationinvariant.full}\\
    \text{s.t.} \quad & \bm{Y} \preceq \mathrm{Diag}(\bm{w}),\ \bm{Y}=\sum_{t=1}^{\color{black}r} \bm{Y}^t,\ \bm{w} \leq \bm{Z}\bm{e}, \nonumber \\ 
    & \sum_{j=1}^p Y_{i,j}^2\leq r Y_{i,i}w_i & \forall i \in [p],\nonumber\\
   & \left(\sum_{j=1}^p \vert Y_{i,j}\vert\right)^2 \leq k Y_{i,i}w_i & \forall i \in [p],\nonumber\\
   &  \sum_{i \in [p]: i \neq j}Y_{i,j}^2 \leq (k-r+1) w_j ( w_j-Y_{j,j}) \ & \forall j \in [p],\nonumber\\
   & \pm \bm{Y}^t\leq \bm{F}^t \ & \forall t \in [r],\nonumber\\
   & G^{t}_{i,1} \geq G^{t}_{i,2} \geq \ldots \geq G^{t}_{i,k_t} \ & \forall i \in [k_t], \ \forall t \in [r],\nonumber\\
    & G^{t}_{i,j}=0 \ & \forall i > k_t \ \text{or} \ j > k_t, \ \forall t \in [r],\nonumber\\
   & \mathrm{tr}(\bm{Y}^t)=\mathrm{tr}(\bm{G}^t) =\mathrm{tr}(\bm{F}^t)=1 \ & \forall t \in [r],\nonumber\\
    & \langle \bm{E}, \bm{G}^t \rangle = \langle \bm{E}, \bm{F}^t\rangle \ & \forall t \in [r],\nonumber\\
    & \sum_{i=1}^j G^t_{i,i}\geq j r^{D,t}_j+\sum_{j=1}^n t_{i,j}^{t,D} \ & \forall j \in [p-1], \ t \in [r],\nonumber\\
    & Y^t_{i,i} \leq r_{j}^{t,D} + t_{i,j}^{t,D} \ & \forall i \in [p], j \in [p-1], t \in [r],\nonumber\\
    & {(F_{i,j}^{t})}^2 \leq T^t_{i,j}T^t_{j,i},\ T^t_{i,i}=F_{i,i}^t \ & \forall i \in [p], j \in [i-1], \ t \in [r],\nonumber\\
    & \sum_{j \in [p]} T^t_{i,j}=Z_{i,t}, \ \sum_{i \in [p]} T^t_{i,j}=k_t F^t_{j,j}\ & \forall i \in [p], \forall j \in [p], t \in [r],\nonumber\\
   & 0 \leq T^t_{i,j}\leq F^t_{i,j} \ & \forall i,j \in [p], t \in [r]\nonumber.
\end{align}
The additional variables $\bm{r}^{t,D}$ and $\bm{t}^{t,D}$ are introduced to enforce coupling constraints between the diagonal entries of $\bm{F}^t$ and $\bm{G}^t$ \citep[][eq. 44]{kim2021convexification}, while $\bm{T}^t$ allows to couple $\bm{F}^t$ with the binary variables $\bm{Z}$ \citep[][eq. 50]{kim2021convexification}.
In contrast to \cite{kim2021convexification}, we explicitly require that each $\bm{F}^t$ is positive semidefinite (rather than that its $2 \times 2$ minors are), in order to obtain a stronger relaxation; we consider the $2 \times 2$ minors when developing a more tractable relaxation in the next section.

\subsection{A Second-Order Cone Relaxation for Large-Scale Instances}\label{ssec:ec.soc}
Unfortunately, \eqref{prob:disjunctiverelax_permutationinvariant} cannot scale beyond $p=100$, at least with current technology, due to the presence of multiple semidefinite matrices and constraints. 

We now develop a more tractable, albeit less tight, version of the relaxation of \eqref{prob:disjunctiverelax_permutationinvariant} which scales to $p>100$ features. Namely, we replace all semidefinite constraints of the form $\bm{X} \in \mathcal{S}^p_+$ with the non-negativity of their $2 \times 2$ minors, $X_{i,i}X_{j,j}\geq X_{i,j}^2 \ \forall i,j \in [p]$, as presented by \cite{bertsimas2020polyhedral} and references therein. This gives the following second-order cone relaxation of \eqref{prob:disjunctiverelax_permutationinvariant}:
\begin{align}
    \max_{\substack{\bm{Z} \in [0, 1]^{p \times r}:\\ \langle \bm{E}, \bm{Z}\rangle \leq k, \bm{w} \in [0, 1]^p}}\max_{\substack{\bm{Y} \in \mathcal{S}^p, \bm{Y}^t, \bm{F}_t, \bm{G}_t \in \mathcal{S}^p,\\ \bm{T}_t  \in \mathbb{R}^{p \times p}_+ \ \forall t \in [{\color{black}r}],\\ \bm{r}^{t,D} \in \mathbb{R}^{p-1}, \bm{t}^{t,D} \in \mathbb{R}^{p \times p-1}_+
    }} \quad & \langle \bm{Y}, \bm{\Sigma}\rangle \label{prob:disjunctiverelax_permutationinvariant_soc}\\
    \text{s.t.} \quad &  \bm{Y}=\sum_{t=1}^{\color{black}r} \bm{Y}^t,\ \mathrm{tr}(\bm{Y}^t)=1, \bm{w} \leq \bm{Z}\bm{e} \ & \ \forall t \in [r], \nonumber\\
    & {Y^t_{i,j}}^2 \leq Y_{i,i}^t Y_{j,j}^t \ & \forall i,j \in [p], \forall t \in [r], \nonumber\\
    & (\delta_{i,j}-Y_{i,j})^2\leq (w_i-Y_{i,i})(w_j-Y_{j,j})\ & \forall i,j \in [p], \forall t \in [r], \nonumber\\
    & \sum_{j=1}^p {Y_{i,j}}^2\leq r Y_{i,i}w_i \ & \forall i \in [p]\nonumber\\
   & \left(\sum_{j=1}^p \vert Y_{i,j}\vert\right)^2 \leq k Y_{i,i}w_i, \ \pm \bm{Y}^t\leq \bm{F}_t \ & \forall i \in [p], \ \forall t \in [r],\nonumber\\
   &  \sum_{i \in [p]: i \neq j}Y_{i,j}^2 \leq (k-r+1) w_j ( w_j-Y_{j,j}) \ & \forall j \in [p],\nonumber\\
   & Y^t_{i,i} \leq t_{i,j}^{t,D}+r_{j}^{t,D} \ & \forall i \in [p], j \in [p-1], t \in [r],\nonumber\\
   & \mathrm{tr}(\bm{Y}^t)=\mathrm{tr}(\bm{G}^t) =\mathrm{tr}(\bm{F}_t)=1 \ & \forall t \in [r],\nonumber\\
   & \langle \bm{E}, \bm{G}_t-\bm{F}_t\rangle=0 \ & \forall t \in [r],\nonumber\\
   & {G^t_{i,j}}^2 \leq G_{i,i}^t G_{j,j}^t \  & \forall i,j \in [p], \forall t \in [r], \nonumber\\
    & G^{t}_{i,1} \geq G^{t}_{i,2} \geq \ldots \geq G^{t}_{i,k_t} \ & \forall i \in [k_t], \ \forall t \in [r],\nonumber\\
    & G^{t}_{i,j}=0 \ & \forall i > k_t \ \text{or} \ j > k_t \ \forall t \in [r],\nonumber\\
    & \sum_{i=1}^j G^t_{i,i}\geq j r^{D,t}_j+\sum_{j=1}^n t_{i,j}^{t,D} \ & \forall j \in [p-1], \ t \in [r],\nonumber\\
    & {F_{i,j}^{t}}^2 \leq T^t_{i,j}T^t_{j,i},\ T^t_{i,i}=F_{i,i}^t \ & \forall i \in [p], j \in [i-1], \ t \in [r],\nonumber\\
    & \sum_{j \in [p]} T^t_{i,j}=Z_{i,t}, \ \sum_{i \in [p]} T^t_{i,j}=k_t F^t_{j,j}\ & \forall i \in [p], \forall j \in [p], t \in [r],\nonumber\\
   & 0 \leq T^t_{i,j}\leq F^t_{i,j}, {F^t_{i,j}}^2 \leq F_{i,i}^t F_{j,j}^t\ & \forall i,j \in [p], \forall t \in [r]. \nonumber
\end{align}


{\blue
\section{Proof of Theorem~\ref{thm:resolve.nphard}} \label{sec:a.nphard}
We prove that Problem \eqref{prob:resolve} is NP-hard via a reduction from a variant of the exact 3-set cover problem. 

\paragraph{Restricted Exact Cover by 3-Sets Problem (X3C-R)} Consider a set $X$ with $|X|=3q$ and a family $\mathcal{C}=\{S_1,\dots,S_m\}$ of 3-element subsets of $X$ of size $m:=\vert \mathcal{C}\vert$. 
Further, assume that for any pair of sets $S,S' \in \mathcal{C}$, $|S \cap S'| \leq 1$.
Then, X3C-R corresponds to asking whether there exists an exact cover of $X$, namely a subfamily $\mathcal{C}'\subseteq \mathcal{C}$ of size $q$, whose sets are pairwise disjoint and whose union is $X$.

Problem X3C-R is NP-complete, as we prove in Proposition \ref{prop.ec.X3CisNPhard}. Leveraging that result, we now prove the following proposition, which implies Theorem \ref{thm:resolve.nphard}.
\begin{proposition}\label{prop:resolve.nphard}
    For any instance of X3C-R, we can construct integers $(p,r)$, a binary mask $\hat{\bm{Z}} \in \{0,1\}^{p \times r}$, a symmetric matrix $\bm{\Sigma} \in \mathcal{S}^p$, and a scalar $\tau$, such that X3C-R is equivalent to finding a feasible solution to \eqref{prob:resolve} with objective value at least $\tau$. Moreover, $p=\vert X\vert+\vert \mathcal{C}\vert=3q+m$, $r=\vert \mathcal{C}\vert$, $\tau=2q$ are all polynomial in the size of X3C-R.
\end{proposition}
\begin{remark}
    Given the orthogonality constraints in \eqref{prob:resolve}, we can replace the matrix $\bm{\Sigma}$ in Problem \eqref{prob:resolve} by $\bm{\Sigma} + \beta \bm{I}_p$ for any $\beta > 0$ and obtain the same optimal solution (it only shifts the objective value by $\beta r$). In particular, if $\bm{\Sigma} \in \mathcal{S}^p$ in Problem \eqref{prob:resolve}, we can assume $\bm{\Sigma} \succeq \bm{0}$ without loss of generality.
\end{remark}
    
\proof{Proof of Proposition \ref{prop:resolve.nphard}}
Consider an instance of X3C-R and let us construct an instance of~\eqref{prob:resolve}. 

Let us first define the coordinate space: We create one coordinate for each element $x\in X$, and one ``dummy'' coordinate $d_t$ for each 3-set $S_t\in\mathcal{C}$. Hence, $p = |X| + |\mathcal{C}| = 3q +m$. 

We consider one PC per 3-set. Hence, $r = |\mathcal{C}| = m$. Each column $t \in [r]$ is allowed to be nonzero only on the coordinates corresponding to $x \in S_t$ and on the dummy $d_t$. In other words,
\begin{align*}
    \hat Z_{i,t}=1, \quad & \mbox{ if } i\in S_t\cup\{d_t\}, \\ 
    \hat Z_{i,t}=0, \quad  & \text{ otherwise.}
\end{align*}
Note that $r \leq m$ and Problem \eqref{prob:resolve} is always feasible (a feasible solution is obtained by setting $\bm{u}_t$ equal to 1 on its dummy coordinate $d_t$ and 0 elsewhere).

Finally, the matrix $\bm{\Sigma}$ is non-zero on the block indexed by $X \times X$ only, i.e., 
\begin{align*}
    \bm{\Sigma} = \begin{pmatrix} \bm{A} & \bm{0} \\ \bm{0} & \bm{0} \end{pmatrix},
\quad \mbox{ with } \quad 
    A_{xy} =
\begin{cases}
1, & \text{if there exists } S\in\mathcal{C} \text{ with } x,y \in S,\\
0, & \text{otherwise.}
\end{cases}
\end{align*}
Consider one column $t \in [r]$. It is associated with one 3-set $S_t$. On this set, the submatrix of $\bm{A}$ is equal to
\begin{align*}
    \bm{A}_{S_t,S_t} = \begin{pmatrix} 0 & 1 & 1 \\ 1 & 0 & 1 \\ 1 & 1 & 0\end{pmatrix} = \bm{e}\bm{e}^\top - \bm{I}_3. 
\end{align*}
So, for any non-zero vector $\bm{z} \in \mathbb{R}^3$, we have by Cauchy-Schwarz
\begin{align}
    \bm{z}^\top \bm{A}_{S_t,S_t} \bm{z} = (\bm{e}^\top \bm{z})^2 - \| \bm{z} \|^2 \leq (\|\bm{z}\|_0 - 1) \| \bm{z} \|^2, \label{eqn:reduc.quad.obj}
\end{align}
with $\|z \|_0 = |\{ j : z_j \neq 0 \} |$.

Set $\tau = 2q$. Let us now show that there exists a feasible solution to \eqref{prob:resolve} with objective value at least $\tau$ iff there is an affirmative answer to X3C-R. 

Consider a solution to X3C-R, $\mathcal{C}'$. For each set $S_t \in \mathcal{C}'$, define $\bm{u}_t$ as equal to $\tfrac{1}{\sqrt{3}}(1,1,1)$ on the coordinates indexed by $S_t$ and 0 on the dummy coordinate. Alternatively, for each set $S_t \notin \mathcal{C}'$, set $\bm{u}_t$ equal to $\bm{0}$ on the coordinates indexed by $S_t$ and 1 on the dummy coordinate. By construction, because each element $x \in X$ is covered by exactly one 3-set of $\mathcal{C}'$, for any $t \neq t'$, the columns $\bm{u}_t$ and $\bm{u}_{t'}$ have disjoint supports and $\bm{u}_t^\top \bm{u}_{t'} = \delta_{t,t'}$. Second, 
\begin{align*}
    \sum_{t \in [r]} \bm{u}_t^\top \bm{\Sigma} \bm{u}_t = \dfrac{1}{3}  \sum_{S_t \in \mathcal{C}'} \bm{e}^\top \bm{A}_{S_t,S_t} \bm{e} = 2q = \tau.
\end{align*}

Conversely, assume there exists a solution to \eqref{prob:resolve} with objective value at least $\tau$. For any $x \in X$, there exists at most one column $\bm{u}_t$ such that $u_{t,x} \neq 0$. Indeed, by contradiction, assume there exist two distinct columns $t,t'$ such that $u_{t,x} \neq 0$ and $u_{t',x} \neq 0$. Then, $x \in S_t \cap S_{t'}$ so we must have $\{ x \} = S_t \cap S_{t'}$ by assumption. In this case, the constraint $\bm{u}_t^\top \bm{u}_{t'} = 0$ reads ${u}_{t,x} u_{t',x} = 0$ and forces one of the coefficients ${u}_{t,x}$ or $ u_{t',x}$ to be equal to 0, which is a contradiction. 
We say that a column of $\bm{U}$, $\bm{u}_t$, is $s$-heavy ($s\in \{0,1,2,3\}$) if $s$ of its coordinates indexed by $X$ are non-zeros. Let $n_s$ denote the number of $s$-heavy columns. From our earlier reasoning, each $x \in X$ is covered by at most one column of $\bm{U}$ so we must have $3 n_3 + 2n_2 \leq 3 n_3 + 2n_2 + n_1 \leq 3q$. The solution $\bm{U}$ achieves an objective value at least equal to $\tau = 2q$ so \eqref{eqn:reduc.quad.obj} implies  $2q = \tau \leq 2 n_3 + n_2$. However, 
\begin{align*}
    2 n_3 + n_2 = \dfrac{2}{3} (3 n_3 + \dfrac{3}{2} n_2) \leq \dfrac{2}{3} (3 n_3 + 2 n_2) \leq 2q,
\end{align*}
so $2q = 2 n_3 + n_2$ and the inequalities above are tight, i.e., $n_2 = 0$ and $n_3 = q$. Consequently, $n_1 = 0$ as well. In conclusion, define $\mathcal{C}' = \{ S_t \: : \: \bm{u}_t \text{ is } 3\text{-heavy} \}$. For any $x \in X$, there exists at most one column $\bm{u}_t$ such that $u_{t,x} \neq 0$, so the sets in $\mathcal{C}'$ are pairwise disjoint. Furthermore, $|\mathcal{C}'| = n_3 = q = |X|/3$ so  $\mathcal{C}'$ must cover the entire set $X$.
\hfill\Halmos
\endproof
\begin{remark}
    Observe that the proof of Proposition \ref{prop:resolve.nphard} shows that there exists a solution to X3C-R iff there exists a feasible solution to Problem \eqref{prob:resolve} with objective value $\tau - \eta$ with $\tau = 2q > 2p/3$ being the optimal value and $\eta \in [0,1)$. Hence,  it shows that it is NP-hard to approximate Problem~\eqref{prob:resolve} within a factor $1-\epsilon$ with $\epsilon = 1/\tau < 3/(2p)$.
\end{remark}

\subsection{X3C-R is NP-Complete} \label{ssec:a.x3cr}
To the best of our knowledge, there is no existing proof of the complexity of X3C-R. 
Problem X3C-R is a restriction of the classical exact cover by 3-sets problem (X3C), which is known to be NP-complete 
\citep[][chapter A3.1]{garey1979computers}.

\paragraph{Exact Cover by 3-Sets Problem (X3C)} Consider a set $X$ with $|X|=3q$ and a family $\mathcal{C}=\{S_1,\dots,S_m\}$ of 3-element subsets of $X$, of size $m:=\vert \mathcal{C}\vert$. Does there exist a subfamily $\mathcal{C}'\subseteq \mathcal{C}$ of size $q$ whose sets are pairwise disjoint and whose union is $X$?

Formally, we have the following result:
\begin{proposition}\label{prop.ec.X3CisNPhard}
Any instance of X3C can be reduced to an instance of X3C-R in polynomial time, and hence Problem X3C-R is NP-hard.
\end{proposition}

Any instance of X3C-R is definitionally an instance of X3C, so Proposition~\ref{prop.ec.X3CisNPhard} implies that X3C-R is NP-complete.

\proof{Proof of Proposition \ref{prop.ec.X3CisNPhard}} We consider an instance $(X,\mathcal{C})$ of X3C. 
Define the set $$Y = X \cup \{z_{S,t} \: : \: S \in \mathcal{C},\ t \in \{1,\dots,6\} \}.$$ 
So, $|Y| = |X| + 6 |\mathcal{C}| = 3q + 6m$.

For any 3-set $C = \{x_1,x_2,x_3\} \in \mathcal{C}$ (the indexing is arbitrary and fixed), define the following sets
\begin{align*}
\begin{array}{lcccccccccr}
    D_{S,1} = \{ x_1 & & & z_{S,1} & & & z_{S,4} & & & \},\\
    D_{S,2} = \{ & x_2 & & & z_{S,2} & & & z_{S,5} & &  \},\\
    D_{S,3} = \{ & & x_3 &  & & z_{S,3} & & & z_{S,6} & \},\\
    D_{S,4} = \{ & & & z_{S,1} & z_{S,2} & z_{S,3} & & & & \},\\
    D_{S,5} = \{ & & & & & & z_{S,4} & z_{S,5} & z_{S,6} & \},
\end{array}
\end{align*}
and define $\mathcal{D} = \{  D_{S,t} \: : \: S \in \mathcal{C},\ t \in \{1,\dots,5\} \}$. Hence, $|\mathcal{D}| = 5 |\mathcal{C}| = 5m$. By construction, for any pair of sets $D,D' \in \mathcal{D}$, we must have $|D \cap D'| \leq 1$. In other words, $(Y,\mathcal{D})$ is a valid X3C-R instance. 

Let us now show that there exists an exact cover of $(X,\mathcal{C})$ iff there exists an exact cover of $(Y,\mathcal{D})$.

Let $\mathcal{C}'$ be an exact cover of $(X,\mathcal{C})$. Define $\mathcal{D}' = \{ D_{S,1}, D_{S,2}, D_{S,3} \: : \: S \in \mathcal{C}' \} \cup \{ D_{S,4}, D_{S,5} \: : \: S \in \mathcal{C} \setminus \mathcal{C}' \}$. We have $|\mathcal{D}'| = 3|\mathcal{C}'| + 2|\mathcal{C} \setminus \mathcal{C}'| = |\mathcal{C}'| + 2|\mathcal{C}| = q + 2m = |Y| / 3$. Because $\mathcal{C}'$ is a cover, for any $x \in X$, there exists $S \in \mathcal{C}'$ such that $x \in S$. So there exists $t \in \{1,2,3\}$ such that $x \in D_{S,t}$ and $ D_{S,t} \in \mathcal{D}'$. For any $z_{S,t}, S \in \mathcal{C}, t\in \{1,\dots,6\}$, either $S \in \mathcal{C}'$ and there must exists $t' \in \{1,2,3\}$ such that $z_{S,t} \in D_{S,t'}$ and $D_{S,t'} \in \mathcal{D}'$, or $S \notin \mathcal{C}$' and there must exists $t' \in \{4,5\}$ such that $z_{S,t} \in D_{S,t'}$ and $D_{S,t'} \in \mathcal{D}'$. Altogether, we showed that any element of $Y$ belongs to an element of $\mathcal{D}'$ (i.e., $\mathcal{D}'$  is a cover). Because elements of $\mathcal{D}'$ are 3-sets and $|\mathcal{D}'| = |Y| / 3$, $\mathcal{D}'$  is an exact cover.

Let $\mathcal{D}'$ be an exact cover of $(Y,\mathcal{D})$ (in particular, $|\mathcal{D}'| = q + 2m$).  Define $\mathcal{C}' = \{ S \in \mathcal{C} \: : \: \exists t \in \{1,2,3\}, D_{S,t} \in \mathcal{D}' \}$. For any $x \in X$, $x \in Y$ and $\mathcal{D}'$ is a cover, so there exist $S \in \mathcal{C}, t \in \{1,2,3\}$ such that $x \in D_{S,t}$ and $ D_{S,t} \in \mathcal{C}'$. Hence, $\mathcal{C}'$ covers $X$. 
Let us consider a set $S \in \mathcal{C}$ such that $D_{S,1} \in \mathcal{D}'$, we must have $D_{S,2}, D_{S,3} \in \mathcal{D}'$. Indeed, we must have $D_{S,4}, D_{S,5} \notin \mathcal{D}'$, because they both overlap with $D_{S,1}$ and elements of $\mathcal{D}'$ are pairwise disjoint. The variable $z_{S,2}$ is covered by $\mathcal{D}'$ so necessarily, we must have $D_{S,2} \in \mathcal{D}'$. Similarly for $z_{S,3}$. Hence, $3 |\mathcal{C}'| + 2 |\mathcal{C} \setminus \mathcal{C}'| = |\mathcal{C}'| + 2m$. Therefore, $|\mathcal{C}'| = q = |X| /3$ and $\mathcal{C}'$ is an exact cover of $X$.
\hfill\Halmos
\endproof

}

{\blue 
\section{Asymptotic Feasibility of the Solutions Returned by Algorithm \ref{alg:alternatingmin}} \label{sec:a.lagrangian.feasibility}
In this section, we analyze the asymptotic feasibility of solutions generated by Algorithm~\ref{alg:alternatingmin} and show how they depend on the method used to solve each subproblem. First, we prove that Algorithm~\ref{alg:alternatingmin} generates a subsequence of solutions that are asymptotically feasible, provided that each subproblem in Algorithm~\ref{alg:alternatingmin} is solved to near-optimality (Section \ref{ssec:ec.bb_conv}). We acknowledge that this assumption may not hold for the truncated power method of \citet{yuan2013truncated} that we use in our implementation (Section \ref{ssec:asymptoticfeas_tpm}). Secondly, under less restrictive assumptions that are satisfied by the truncated power method, we show convergence to stationary points of the function $\bm{U} \mapsto \sum_{t' > t} (\bm{U}_t^\top \bm{U}_{t'})^2$ over the feasible region. Unfortunately, this result cannot rule out convergence to infeasible solutions, which is a known limitation of the quadratic penalty method, even when subproblems are solved with increasing accuracy \citep[see][Theorem 17.2]{nocedal2006numerical}. 

We introduce some compact notation from \citet[][chapter 17]{nocedal2006numerical}. Letting $\mathbb{S}^p_k$ denote the set of $k$-sparse vectors from the $p$-dimensional unit sphere, we rewrite Problem \eqref{prob:spcaorth_compact} as a generic-equality constrained maximization problem of the form
\begin{align*}
    \max_{\bm{U} \in \mathcal{U}} \:  f(\bm{U})\quad \text{ s.t. }\quad \bm{c}(\bm{U})=\bm{0},
\end{align*}
where the feasible set $\mathcal{U} := \{ \bm{U} \in \mathbb{R}^{p \times r} \: : \: \forall t \in [r], \ \bm{U}_t \in \mathbb{S}_{k_t}^p\}$ is non-convex, the objective function $f$ is defined as  $f(\bm{U})= \sum_{t \in [r]} \bm{U}_t^\top \bm{\Sigma} \bm{U}_t$, and the constraint function $\bm{c} : \mathbb{R}^{n \times r} \rightarrow \mathbb{R}^{r(r-1)/2}$ is defined as $c_{t,t'}(\bm{U}) = \bm{U}_t^\top \bm{U}_{t'}$ for $t > t'$. 
\subsection{Asymptotic Convergence of Algorithm \ref{alg:alternatingmin} With Branch-and-Bound for Subproblems}\label{ssec:ec.bb_conv}
Suppose that, at each iteration $\ell$ of Algorithm~\ref{alg:alternatingmin}, the vector $\bm{U}_t^{(\ell)}$ is computed by solving 
\begin{align} \label{eqn:alg2.subp}
    \max_{\bm{u} \in \mathbb{S}^{p}_{k_t}} \quad \left\langle  \bm{u}\bm{u}^\top, \bm{\Sigma} - \sum_{t' \in [r]: t' \neq t} \lambda_{t,t'}^{(\ell)} \bm{U}_{t'}^{\blue (\ell-1)}\bm{U}_{t'}^{\blue (\ell-1) \top} \right\rangle.
\end{align}
Under this assumption, we now establish the asymptotic feasibility of a subsequence of the solutions generated by Algorithm \ref{alg:alternatingmin}, provided the subproblems \eqref{eqn:alg2.subp} are solved to finite (but not necessarily improving) accuracy, and $k_t \geq r$ for each $t \in [r]$.

\begin{remark}\label{rem:jacobi} The subproblem \eqref{eqn:alg2.subp} for updating column $t$ at iteration $\ell$ corresponds to a block-Jacobi update rule where the penalty on orthogonality violation is expressed using the value of $\bm{U}$ at the previous iteration $\ell-1$. In particular, it requires storing both $\bm{U}^{(\ell)}$ and $\bm{U}^{(\ell-1)}$ and does not use the updated values of $\bm{U}^{(\ell)}_{t'}$ for $t'<t$. Instead, we can implement our algorithm (and adapt our analysis) using Gauss-Seidel updates instead, which would involve subproblems of the form
\begin{align*} 
    \max_{\bm{u} \in \mathbb{S}^{p}_{k_t}} \quad \left\langle  \bm{u}\bm{u}^\top, \bm{\Sigma} - \sum_{t' < t} \lambda_{t,t'}^{(\ell)} \bm{U}_{t'}^{ (\ell)}\bm{U}_{t'}^{(\ell) \top} - \sum_{t' > t} \lambda_{t,t'}^{(\ell)} \bm{U}_{t'}^{\blue (\ell-1)}\bm{U}_{t'}^{\blue (\ell-1) \top} \right\rangle.
\end{align*}
\end{remark}
\begin{proposition} \label{prop.ec.convsubsequence}
Assume that $k_t \geq r$, for all $t \in [r]$. Consider a sequence of penalty parameters $\{ \lambda_{t,t'}^{(\ell)} \}$ with $\mu_{\ell} := \min_{t,t' \in [r]} \lambda_{t,t'}^{(\ell)} \rightarrow + \infty$  and a sequence of iterates $\{ \bm{U}^{(\ell)} \}$ such that, for any $\ell \geq 1$, any $t \in [r]$, and any sparse vector $\bm{u} \in \mathbb{S}^p_{k_t}$, 
\begin{align}
    \bm{U}_{t}^{(\ell)\top}\bm{\Sigma}\bm{U}_{t}^{(\ell)} - \sum_{t' \neq t} \lambda_{t,t'}^{(\ell)} (\bm{U}_{t}^{(\ell)\top} \bm{U}^{(\ell-1)}_{t'})^2 \geq \bm{u}^{\top}\bm{\Sigma}\bm{u} - \sum_{t' \neq t} \lambda_{t,t'}^{(\ell)} (\bm{u}^{\top} \bm{U}^{(\ell-1)}_{t'})^2 - \epsilon_\ell,
\end{align}
with $\epsilon_\ell / \mu_{\ell} \rightarrow 0$. Then, for any convergent subsequence of $\{ \bm{U}^{\ell}\}$, its limit point $\bar{\bm U}$ is feasible. 
\end{proposition} 
\begin{remark} Because $\bm{U}^{(\ell)} \in \mathcal{U}$, which is compact, a converging subsequence exists by the Bolzano-Weierstrass theorem.
\end{remark}
When using a provably optimal method, like the branch-and-bound algorithm of \citet{berk2019certifiably}, for solving the subproblems \eqref{eqn:alg2.subp}, we can set the suboptimality gap $\epsilon_\ell$ explicitly as the termination criterion and ensure the condition of Proposition \ref{prop.ec.convsubsequence} holds (e.g., by keeping $\epsilon_\ell$ constant). However, empirically, solutions from the truncated power method (TPM) are very close to the optimal solution, but we cannot guarantee that TPM solves \eqref{eqn:alg2.subp} to $\epsilon_\ell$-optimality for some controllable $\epsilon_\ell > 0$.

\proof{Proof of Proposition \ref{prop.ec.convsubsequence}} Consider the $\ell$th iteration and fix $t \in [r]$. Since $k_t \geq r$, we can construct $\bm{u} \in \mathbb{S}^p_{k_t}$ such that $\bm{u}^\top \bm{U}_{t'}^{(\ell-1)} = 0, \forall t' \neq t$ (e.g., select $r$ coordinates, the restriction of the $r-1$ vectors $\bm{U}_{t'}^{(\ell-1)}, t' \neq t$ to these $r$ coordinates cannot span the entire space $\mathbb{R}^r$ so we can find a unit vector orthogonal to all of them, which we can view as an $r$-sparse unit vector in $\mathbb{R}^p$ by padding with zeros). For this vector, we have, by assumption
\begin{align*}
    \bm{U}_{t}^{(\ell)\top}\bm{\Sigma}\bm{U}_{t}^{(\ell)} - \sum_{t' \neq t} \lambda_{t,t'}^{(\ell)}  (\bm{U}_{t}^{(\ell)\top} \bm{U}^{(\ell-1)}_{t'})^2 \geq \bm{u}^{\top}\bm{\Sigma}\bm{u}- \epsilon_\ell,
\end{align*}
leading to
\begin{align*}
    \sum_{t' \neq t} (\bm{U}_{t}^{(\ell)\top} \bm{U}^{(\ell-1)}_{t'})^2 \leq \dfrac{1}{\mu_\ell} \left( \bm{U}_{t}^{(\ell)\top}\bm{\Sigma}\bm{U}_{t}^{(\ell)} - \bm{u}^{\top}\bm{\Sigma}\bm{u} \right) + \dfrac{\epsilon_\ell}{\mu_\ell} \leq \dfrac{\lambda_{\max}(\bm{\Sigma}) + \epsilon_\ell}{\mu_\ell} .
\end{align*}
Taking limits when $\ell \rightarrow \infty$, we get $\displaystyle \sum_{t' \neq t} (\bar{\bm{U}}_{t}^{\top} \bar{\bm{U}}_{t'})^2 = 0$, hence $\bar{\bm{U}}_{t}^{\top} \bar{\bm{U}}_{t'} = 0,\ \forall t \neq t'$.
\hfill\Halmos\endproof

\subsection{Asymptotic Stationarity of Algorithm \ref{alg:alternatingmin} With Truncated Power Method for Subproblems}\label{ssec:asymptoticfeas_tpm}
In this section, we analyze the convergence of Algorithm \ref{alg:alternatingmin} when the subproblems \eqref{eqn:alg2.subp} are solved using the truncated power method of \citet{yuan2013truncated}. In particular, we leverage the following property, which follows immediately from the description of the method in \citet{yuan2013truncated}:

\begin{lemma}\label{lemma.ec.tpm} Consider a symmetric matrix $\bm{A}\in\mathbb{R}^{p\times p}$, $k \leq p$, and the sparse quadratic problem
\begin{align} \label{eqn:sparse.eigen}
\max_{\bm{x} \in\mathbb{R}^p} \quad \bm{x}^\top \bm{A x}
\ \text{s.t.}\ \bm{x} \in \mathbb{S}^p_k,
\end{align}
The truncated power method of \citet{yuan2013truncated} applied to \eqref{eqn:sparse.eigen} converges to a sparse eigenvector of $\bm{A}$, i.e., a vector $\bar{\bm x} \in \mathbb{R}^p : \| \bm{x} \|_2  =1$, such that the support of $\bar{\bm x}$, $S := \{ j \: : \: \bar{x}_j \neq 0 \}$ is of size at most $k$ and such that there exists $\lambda \in \mathbb{R}$ such that $\bm{A}_{S,S} \bar{\bm x}_S = \lambda \bar{\bm x}_S$.
\end{lemma}

For any $\bm{x} \in \mathbb{S}^p_k$, we denote $\mathcal{T}_{\bm{x}}(\mathbb{S}^p_k)$ the tangent cone to $\mathbb{S}^p_k$ at $\bm{x}$, and $\operatorname{Proj}_{\bm{x}}$ the projection onto $\mathcal{T}_{\bm{x}}(\mathbb{S}^p_k)$. In particular, we have 
\begin{align*}
    \mathcal{T}_{\bm{x}}(\mathbb{S}^p_k) = \{ \bm{y} \in \mathbb{R}^p \: : \: \operatorname{supp}(\bm{y}) \subseteq \operatorname{supp}(\bm{x}),\ \bm{y}^\top \bm{x} = 0 \}.
\end{align*}
With these notations, the equation $\bm{A}_{S,S} \bar{\bm x}_S = \lambda \bar{\bm x}_S$ in Lemma \ref{lemma.ec.tpm} implies 
\begin{align*}
\bm{y}^\top \bm{A} \bar{\bm x} = 0, \quad \forall \bm{y} \in  \mathcal{T}_{\bm{x}}(\mathbb{S}^p_k), \quad
    \mbox{ and } \quad & \left\| \operatorname{Proj}_{\bm{x}}\left( \bm{A} \bar{\bm x} \right)\right\| = \bm{0}.
\end{align*}

At iteration $\ell$ of Algorithm~\ref{alg:alternatingmin}, the vector $\bm{U}_t^{\ell}$ is computed by solving a problem of the form \eqref{eqn:sparse.eigen} with $\bm{A} = \bm{\Sigma} - \sum_{t' \neq t}  \lambda_{t,t'}^{(\ell)} \bm{U}^{(\ell-1)}_{t'}\bm{U}_{t'}^{(\ell-1)\top}$ using the truncated power method. Hence, Lemma~\ref{lemma.ec.tpm} yields 
\begin{align*}
    \left\| \operatorname{Proj}_{\bm{U}_t^{(\ell)}} \left( \bm{\Sigma}\bm{U}_{t}^{(\ell)} - \sum_{t' \neq t} \lambda_{t,t'}^{(\ell)} (\bm{U}_{t}^{(\ell)\top} \bm{U}^{(\ell-1)}_{t'})\,\bm{U}_{t'}^{(\ell-1)} \right) \right\| = 0.
\end{align*}
Thanks to this condition, we can show that the sequence of iterates generated by Algorithm~\ref{alg:alternatingmin} converges to a stationary point of the quadratic penalty function. 
Unfortunately, stationarity does not imply feasibility and Algorithm~\ref{alg:alternatingmin} can converge to an infeasible, yet stationary, point. This behavior is a known limitation of the quadratic penalty method, which occurs even when subproblems are solved to increasing accuracy \citep[see][Theorem 17.2]{nocedal2006numerical}.

Formally, we have the following result: 
\begin{proposition}\label{prop:penalty.guarantee} Consider a sequence of homogeneous penalty parameters  $\{  \lambda_{t,t'}^{(\ell)} \} = \{ \mu_\ell \}$ with $\mu_\ell \rightarrow + \infty$  and a sequence of iterates $\{ \bm{U}^{(\ell)} \}$ such that  
\begin{align*}
    \max_{t \in [r]} \: \left\| \operatorname{Proj}_{\bm{U}_t^{(\ell)}}\left( \bm{\Sigma}\bm{U}_{t}^{(\ell)} - \mu_\ell \sum_{t' \neq t} (\bm{U}_{t}^{(\ell)\top} \bm{U}^{(\ell-1)}_{t'})\,\bm{U}_{t'}^{(\ell-1)}  \right) \right\| \leq \tau_\ell,
\end{align*}
with $\tau_\ell/\mu_\ell \rightarrow 0$. Then, for any converging subsequence of $\{ \bm{U}^{\ell}\}$ such that the sequence $\{ \operatorname{supp}(\bm{U}_t^{(\ell)})\}, t \in [r]$ converges, its limit point $\bar{\bm U}$ is a stationary point of the function $\bm{U} \mapsto \| \bm{c}(\bm{U}) \|^2$ on $\mathcal{U}$, namely 
\begin{align*}
\operatorname{Proj}_{\bar{\bm{u}}_t} \left( 
\sum_{t' \neq t} (\bar{\bm u}_t^\top \bar{\bm u}_{t'}) \, \bar{\bm u}_{t'} \right) = \bm{0}.
\end{align*}
\end{proposition} 
Proposition \ref{prop:penalty.guarantee} asserts that limit points are stationary points of $\bm{U} \mapsto \| \bm{c}(\bar{\bm{U}})\|^2$ over the feasible set, meaning that for each column $t$ of a limit point $\bar{\bm{U}}$, the directional derivative of $\| \bm{c}(\bar{\bm{U}})\|^2$ vanishes along every direction in $\mathcal{T}_{\bar{\bm{U}}_t}(\mathbb{S}^p_{k_t})$. 
To see that stationarity does not imply feasibility, observe, for example, that with $r=2$, Proposition~\ref{prop:penalty.guarantee} does not rule out convergence to a matrix $\bm{U} = [\bm{u}_1,\bm{u}_1] \in \mathbb{R}^{p \times 2}$ with $\bm{u}_1 \in \mathbb{S}^p_k$.
Provided the truncated power method is run until convergence at each iteration, Algorithm \ref{alg:alternatingmin} satisfies the assumptions of Proposition \ref{prop:penalty.guarantee} with $\tau_\ell = 0$. Nonetheless, Proposition \ref{prop:penalty.guarantee} shows convergence to a stationary point even when this condition is approximately satisfied.

\proof{Proof of Proposition \ref{prop:penalty.guarantee}} We follow the proof steps of \citet[][Theorem 17.2]{nocedal2006numerical}. Denoting $\bm{w}^{(\ell)}_t := \bm{\Sigma}\bm{U}_{t}^{(\ell)} - {\mu_\ell} \sum_{t' \neq t} (\bm{U}_{t}^{(\ell)\top} \bm{U}^{(\ell-1)}_{t'})\,\bm{U}_{t'}^{(\ell-1)}$, we have 
\begin{align*}
 \operatorname{Proj}_{\bm{U}_t^{(\ell)}}  \left( \sum_{t' \neq t} (\bm{U}_t^{(\ell)\top}\bm{U}_{t'}^{(\ell-1)}) \bm{U}_{t'}^{(\ell-1)} \right) = \dfrac{1}{\mu_\ell} \,  \operatorname{Proj}_{\bm{U}_t^{(\ell)}} \left(  \bm{\Sigma} \bm{U}_t^{(\ell)} - \bm{w}^{(\ell)}_t \right)
\end{align*}
and our assumption leads to
\begin{align*}
\left\| \operatorname{Proj}_{\bm{U}_t^{(\ell)}}  \left( \sum_{t' \neq t} (\bm{U}_t^{(\ell)\top}\bm{U}_{t'}^{(\ell-1)}) \bm{U}_{t'}^{(\ell-1)} \right)  \right\| &\leq \dfrac{1}{\mu_\ell} \left( \|  \operatorname{Proj}_{\bm{U}_t^{(\ell)}}  \left( \bm{\Sigma} \bm{U}_t^{(\ell)}\right) \| + \|  \operatorname{Proj}_{\bm{U}_t^{(\ell)}}  \left( \bm{w}^{(\ell)}_t \right) \| \right) \\
&\leq \dfrac{1}{\mu_\ell} \left( \lambda_{\max}(\bm{\Sigma}) + \tau_{\ell} \right)
\end{align*}

Let us consider a converging subsequence of $\{ \bm{U}^{\ell}\}$ such that the sequences $\{ \operatorname{supp}(\bm{U}_t^{(\ell)})\}$ for $ t \in [r]$ converge. Since the sequences $\{ \operatorname{supp}(\bm{U}_t^{(\ell)}) \}$ take discrete values and converge, there exists $\ell_0$ such that, for any $\ell \geq \ell_0$, $t \in [r]$, $\operatorname{supp}(\bm{U}_t^{(\ell)}) = \operatorname{supp}(\bar{\bm{U}}_t)$. Hence, for any $\ell \geq \ell_0$, $\operatorname{Proj}_{\bm{U}_t^{(\ell)}}  \left( \cdot \right) \rightarrow \operatorname{Proj}_{\bar{\bm{U}}_t}  \left( \cdot \right)$ (for $\ell \geq \ell_0$, the projection operators only differ in the linear constraint defining the tangent cone) and  
taking  the limit in the inequality above, for each $t\in[r]$, 
\begin{align*}
\operatorname{Proj}_{\bar{\bm{U}}_t} \left( 
\sum_{t' \neq t} (\bar{\bm U}_t^\top \bar{\bm U}_{t'}) \, \bar{\bm U}_{t'} \right) = \bm{0}. \hfill\Halmos
\end{align*}
\endproof

\begin{remark} At the expense of more complicated notation, Proposition \ref{prop:penalty.guarantee} can be extended to non-homogeneous penalty sequences $\{ \lambda_{t,t'}^{(\ell)} \}$ such that $\mu_\ell := \min_{t,t'} \lambda_{t,t'}^{(\ell)} \rightarrow + \infty$ and such that, for each pair $(t,t')$, there exists $\lambda_{t,t'}>0$ s.t. $\lambda_{t,t'}^{(\ell)} /\mu_\ell \rightarrow \lambda_{t,t'}$ (i.e., all $\lambda_{t,t'}^{(\ell)}$s grow at the same rate), in which case we have convergence to a stationary point of the function  $\displaystyle \bm{U} \mapsto \sum_{t > t'} \lambda_{t,t'}(\bm{U}_t^\top \bm{U}_{t'})^2$.    
\end{remark}

}

\section{Alternative Proof of Theorem \ref{thm.MILPbound_multiPCs}} \label{sec:ec.combinatorial.linalgproof}
For a fixed $\bm{Z}$, the objective value in 
\begin{align*}
    \max_{\bm{U} \in \mathbb{R}^{p \times r}} \quad \langle \bm{U}\bm{U}^\top, \bm{\Sigma}\rangle
    \quad \text{s.t.} \quad \bm{U}^\top \bm{U}=\mathbb{I}, \,  U_{i,t}=0 \ \text{if} \ Z_{i,t}=0, \ \forall i \in [p], \ \forall t \in [r],
\end{align*}
is equal to
\begin{align*}
\langle \bm{U}\bm{U}^\top, \bm{\Sigma}\rangle = \sum_{t \in [r]} \langle \bm{U}_t\bm{U}_t^\top, \bm{\Sigma}\rangle = \sum_{t \in [r]} \langle \bm{U}_t\bm{U}_t^\top, \mathrm{Diag}(\bm{Z}_t) \bm{\Sigma} \mathrm{Diag}(\bm{Z}_t)\rangle = \operatorname{vec}(\bm{U})^\top \bm{S} \operatorname{vec}(\bm{U}),
\end{align*}
with 
\begin{align*}
    \bm{S} := \begin{pmatrix}
    \mathrm{Diag}(\bm{Z}_1) \bm{\Sigma} \mathrm{Diag}(\bm{Z}_1) & & \\ 
    & \ddots & \\ & & \mathrm{Diag}(\bm{Z}_r) \bm{\Sigma} \mathrm{Diag}(\bm{Z}_r)
    \end{pmatrix},
\end{align*}
and $\operatorname{vec}(\bm{U}) \in \mathbb{R}^{pr}$ the vector obtained by concatenating the columns of $\bm{U}$, $\bm{u}_t$, vertically.
As in the proof of \eqref{eqn:bound.topr}, we use the fact that $\bm{S} \preceq \operatorname{Diag}(\bm{s})$ with $\bm{s} \in \mathbb{R}^{pr}$ the vector of absolute row-sums of $\bm{S}$, i.e., $s_{(t-1) \times p + i} = \sum_{j \in [p]} Z_{i,t} Z_{j,t} |\Sigma_{i,j}|$ for $i \in [p]$, $t \in [r]$. This leads to
\begin{align*}
\langle \bm{U}\bm{U}^\top, \bm{\Sigma} \rangle \leq 
\operatorname{vec}(\bm{U})^\top \operatorname{Diag}(\bm{s}) \operatorname{vec}(\bm{U}) 
= \sum_{t \in [r]} \sum_{i \in [p]} U_{i,t}^2 \left( \sum_{j \in [p]} Z_{i,t} Z_{j,t} |\Sigma_{i,j}| \right).
\end{align*}
Hence, we can bound the objective value by 
\begin{align*}
    \max_{\bm{U} \in \mathbb{R}^{p \times r}} \quad \sum_{t \in [r]} \sum_{i \in [p]} U_{i,t}^2 \left( \sum_{j \in [p]} Z_{i,t} Z_{j,t} |\Sigma_{i,j}| \right)
    \quad \text{s.t.} \quad \bm{U}^\top \bm{U}=\mathbb{I}.
\end{align*}
On one side, the orthogonality constraint implies $\| \bm{U}_t \|_2^2 = \sum_{i \in [p]} U_{i,t}^2 = 1$, for any $t \in [r]$. On the other side, $\bm{U}^\top \bm{U}=\mathbb{I} \implies \bm{U} \bm{U}^\top \preceq \mathbb{I} \implies \sum_{t \in [r]} U_{i,t}^2 \leq 1, \forall i \in [p]$. Hence, optimizing for $\mu_{i,t} := U_{i,t}^2$, we get a bound of the form
\begin{align*}
    \max_{\bm{\mu} \geq 0} \quad \sum_{t \in [r]} \sum_{i \in [p]} \mu_{i,t} \left( \sum_{j \in [p]} Z_{i,t} Z_{j,t} |\Sigma_{i,j}| \right)
    \quad \text{s.t.} & \sum_{i \in [p]} \mu_{i,t} = 1, \forall t \in [r] \\ & \sum_{t \in [r]} \mu_{i,t} \leq 1, \forall i \in [p].
\end{align*}
Finally, we recognize that some optimal solution $\bm{\mu}$ to the above problem is an extreme point (i.e., binary) by the linearity of the objective, giving the overall result.

\section{Algorithmic Benchmark: Non-Convex QCQP Solvers}
\label{sec:ec.branchandbound}
The sparse PCA problem with multiple PCs, either in its original formulation \eqref{prob:spcaorth_compact} or its equivalent reformulation \eqref{prob:spca_extended}, can be seen as a non-convex mixed-integer quadratically constrained problem---for \eqref{prob:spca_extended}, the rank constraints can be encoded as non-convex quadratic constraints $(\bm{Y}^t)^2 = \bm{Y}^t$ or $\bm{Y}^t = \bm{U}_t \bm{U}_t^\top$. However, current non-convex MIQCP solvers cannot handle SDP variables and constraints. Accordingly, we now discuss how to solve our sparse PCA problem exactly using commercial global optimization solvers via formulation \eqref{prob:spcaorth_compact}, and evaluate this option numerically. 

\subsection{Solving \eqref{prob:spcaorth_compact} With an Off-the-Shelf MIQCQP Solver}
Since feasible solutions are extremely challenging for spatial branch-and-bound solvers to recover when quadratic equality constraints are imposed exactly, especially when these solutions are irrational or of exponential size \citep[cf.][]{ramana1997exact, bienstock2023complexity},
we relax the constraint $\bm{U}^\top \bm{U}=\mathbb{I}$ to require that it is satisfied to within an elementwise tolerance of $\epsilon$.
This gives:
\begin{align}\label{prob:spcaorth_compact2}
    \max_{\substack{\bm{Z} \in \{0, 1\}^{p \times r}:\\ \langle \bm{E}, \bm{Z}\rangle \leq k}}\max_{\bm{U} \in \mathbb{R}^{p \times r}} \quad & \langle \bm{U}\bm{U}^\top, \bm{\Sigma}\rangle \\
    \text{s.t.} \quad & \Vert \bm{U}^\top \bm{U}-\mathbb{I}\Vert_\infty \leq \epsilon,\nonumber\\
    & U_{i,t}=0 \ \text{if} \ Z_{i,t}=0, \quad & \forall i \in [p], \ t \in [r]. \nonumber
\end{align}
We set $\epsilon=10^{-4}/r^2$ so that the total constraint violation does not exceed $10^{-4}$. 
Problem \eqref{prob:spcaorth_compact2} is a non-convex quadratically constrained mixed-integer problem with $pr$ continuous variables, $pr$ binaries, and $r^2$ quadratic constraints. 

In addition, we strengthen Problem \eqref{prob:spcaorth_compact2} with valid inequalities derived from the $\ell_1$ relaxation of sparse PCA, as explored by \citet{dey2020solving,dey2021using}. Indeed, if the sparsity of each PC, $k_t$, is specified a priori, we have the valid inequalities
\begin{align}\label{eqn:ell1relax1}
    \Vert \bm{U}_t\Vert_1 \leq \sqrt{k_t}, \ \forall t \in [r].
\end{align}
Moreover, if $k$ is specified but $k_t$ is not, we instead impose the second-order cone inequalities
\begin{align}\label{eqn:ell1relax2}
    \Vert \bm{U}_t\Vert_1^2 \leq \sum_{i \in [p]}Z_{i,t}, \ \forall t \in [r],
\end{align}
which allows us to model $k_t=\sum_{j=1}^p Z_{j,t}$ in a tractable fashion. 

In practice, this approach allows MIQCP solvers to solve Problem \eqref{prob:spca_extended} to optimality for $pr < 100$ (Section \ref{ssec:feasiblemethods}) and obtain high-quality solutions at larger problem sizes. Note that we avoid mixing both sets of inequalities, as we observed in some preliminary numerical experiments that this sometimes induces numerical instability. 

To further improve branch-and-bound, we consider two acceleration strategies:
\begin{itemize}
    \item {Using} the solution generated by Algorithm \ref{alg:greedymethod2} as a warm-start;
    \item {Inform branching decisions by} the combinatorial upper bound derived in Section \ref{sec.gen.gershgorin}. 
\end{itemize}

As spatial branch-and-bound technology improves over time, we believe that it should be possible to solve Problem \eqref{prob:spcaorth_compact2} exactly at larger problem sizes. Indeed, recent works, e.g. \citet{gupta2022branch}, solve some quadratically constrained problems with up to $50$ variables to optimality using custom branch-and-bound solvers, and \citet{gupta2022branch} reports that \verb|Gurobi|'s off-the-shelf QCQP solver has achieved a machine-independent speedup factor of $67.5$ in less than two years, which suggests that larger instances of \eqref{prob:resolve}--\eqref{prob:spcaorth_compact2} may soon be in reach.

\subsection{Numerical Performance on Pitprops}
In this section, we {\color{black}numerically evaluate the quality of our approaches} on the \texttt{pitprops} dataset (Table \ref{tab:comparison_feasible_pitprops}) and an approach based on the currently available non-convex MIQCQP technology.

\begin{table}[h]
\centering\footnotesize
\begin{tabular}{@{}l l r r r r r r r r r r r@{}} \toprule
$r$ &  $k_t$ & \multicolumn{4}{c@{\hspace{0mm}}}{Alg. \ref{alg:greedymethod2}} & \multicolumn{3}{c@{\hspace{0mm}}}{Alg. \ref{alg:alternatingmin}} & \multicolumn{3}{c@{\hspace{0mm}}}{Alg. \ref{alg:disjoint.linalg}} \\
\cmidrule(l){3-6} \cmidrule(l){7-9} \cmidrule(l){10-12} &  & UB & Obj. & Viol. & T(s) & Obj. & Viol. & T(s) & Obj. & Viol. & T(s) \\\midrule
2 &	2 &	0.295 & \textbf{0.295} &	0 &	20.64 & 
\textbf{0.295} &	0 &	7.54 & \textbf{0.295} &	0 & 0.02 \\
2 &	4 &	0.408 & 0.378 &	0 &	20.80 & {0.400} &	0 &	4.94 & \textbf{0.404} &	0 &	0.52 \\
2 &	6 &	0.477& 0.437 &	0 &	20.62 &	{0.452} &	0 &	6.77 & 0.443 &	0 &	0.62 \\
2 &	8 &	0.501 & 0.375 &	0 &	20.93 &	\textbf{0.476} &	0 &	7.83 & 0.446 &	0 &	0.64 \\
2 &	10 & 0.507 & 0.463 &	0 &	21.01 &	\textbf{0.500} &	0 &	6.36 & 0.464 &	0 &	0.72 \\\midrule
 3 &	2 &	0.435 & \textbf{0.435} &	0 &	20.81 &	{0.424} &	0 &	9.98 & \textbf{0.435} & 0 & 0.03 \\
 3 &	4 &	0.572 & {0.525} &	0 &	21.03 & {0.551} &	0 &	7.24 & \textbf{0.555} &	0 &	1.01 \\
3 &	6 &	0.641 & 0.463 &	0 &	23.03 &	\textbf{0.608} &	0 &	10.21 & 0.569 &	0 &	0.92 \\
3 &	8 &	0.652 & 0.580 &	0 &	20.94 &	\textbf{0.638} &	0 &	8.88 & 0.569 &	0 &	1.06 \\
3 &	10 &	0.652 & 0.392 &	0 &	20.81 &	\textbf{0.650} &	0 &	11.43 & 0.569 &	0 &	1.39 \\\midrule
4 &	2 &	0.554 & \textbf{0.554} &	0 &	21.22 &	\textbf{0.554} &	0 &	11.05 & \textbf{0.554} &	0 &	0.97 \\
4 &	4 &	0.704 & 0.470 &	0 &	21.07 &	 {0.657} &	0.003 &	13.81 & \textbf{0.657} &	0 &	3.25 \\
4 &	6 &	0.737 & 0.537 &	0 &	22.22 &	\textbf{0.697} &	0.002 &	12.14 & 0.644 &	0 &	2.07 \\
4 &	8 &	0.737 & 0.553 &	0 &	22.71 &	\textbf{0.720} &	0 &	11.84 & 0.644 &	0 &	3.48 &	\\
4 &	10 & 0.737 &	0.508 &	0 &	21.03 &	\textbf{0.736} &	0 &	11.22 & 0.644 &	0 &	3.16 \\\midrule
5 &	2 &	0.657 & 0.455 &	0 &	20.87 &	{0.647} &	0 &	12.23 & \textbf{0.648} &	0 &	1.57 \\
5 &	4 &	0.795 & 0.586 &	0 &	21.27 &	\textbf{0.743} &	0 &	15.60 & 0.709 &	0 &	8.12 \\
5 &	6 &	0.807 & 0.538 &	0 &	23.47 &	\textbf{0.779} &	0.016 &	13.95 & 0.713 &	0 &	16.54 \\
5 &	8 &	0.807 & 0.563 &	0 &	21.02 &	\textbf{0.800} &	0.004 &	17.94 & 0.713 &	0 &	24.78 \\
5 &	10 &	0.807 & 0.525 &	0 &	21.97 &	\textbf{0.807} &	0.001 &	15.01 & 0.713 &	0 &	7.39 \\\midrule
6 &	2 &	0.749 & {0.581} &	0 &	21.45 &	0.746 & 0 &	12.61 & \textbf{0.749} &	0 &	4.87 \\
6 &	4 &	0.866 & 0.576 &	0 &	21.23 &	{0.807} &	0.035 &	23.08 & \textbf{0.780} &	0 &	54.2\\
6 &	6 &	0.870 & 0.617 &	0 &	23.45 &	\textbf{0.839} &	0.044 &	15.63 & 0.780 &	0 &	12.54    \\
6 &	8 &	0.870 & 0.628 &	0 &	21.69 &	\textbf{0.849} &	0.011 &	17.20 & 0.780 &	0 &	8.45  \\
6 &	10 &	0.870 & 0.664 &	0 &	21.12 &	\textbf{0.866} &	0.006 &	25.91 & 0.780 &	0 &	86.88 \\\midrule
Avg & & 0.668 & 0.508 &	0 &	21.46 &	\textbf{0.649} & 0.005 &	12.42 & 0.610 &	0 &	9.85 &	\\
\bottomrule
\end{tabular}
\caption{Performance of Algorithms \ref{alg:greedymethod2}--\ref{alg:disjoint.linalg} 
on the pitprops dataset ($p=13$) using the experimental setup laid out in Section \ref{ssec:feasiblemethods}. We denote the best-performing solution (in terms of the proportion of variance explained minus the total orthogonality constraint violation) in bold. 
We report the semidefinite upper bound obtained from solving Problem \eqref{prob:disjunctiverelax_permutationinvariant} as part of our analysis of Algorithm $1$ since the semidefinite upper bound is the tightest bound proposed in the paper, but do not report the upper bound from the other methods to avoid redundancy. Note that $k_t$ denotes the sparsity of each individual component, meaning a set of $r$ PCs have a collective sparsity budget of $k_t r$, and that all objective values are reported in terms of the proportion of variance explained by dividing by $p$, the number of features. Note that the relative optimality gap from Gurobi was less than $10^{-4}$ at termination for all results in this table.
} 
\label{tab:comparison_feasible_pitprops}
\end{table}

Tables \ref{tab:comparison_gurobi_pitprops1}-\ref{tab:comparison_gurobi_pitprops2} compare commercial spatial branch-and-bound with/without the combinatorial upper bound \eqref{prob.linearalgebraicbound}, and with/without a warmstart from Algorithm \ref{alg:greedymethod2}.

In comparison with the performance of Algorithms \ref{alg:alternatingmin}--\ref{alg:disjoint.linalg}, we observe that (i) Gurobi's upper bounds are uniformly worse than the SDP relaxation for $r\geq 4$; 
(ii) on instances where  $\sum_t k_t \leq p$ it explains a comparable amount of variance to Algorithms \ref{alg:alternatingmin}--\ref{alg:disjoint.linalg}, although it explains significantly less variance than Algorithms \ref{alg:alternatingmin} on instances where $\sum_t k_t > p$; (iii) It requires 3--4 orders of magnitude more time than any of our algorithms.

\begin{table}[h]
\centering\footnotesize
\begin{tabular}{@{}l l r r r r r r r r r r r r@{}} \toprule
$r$ &  $k_t$ & \multicolumn{6}{c@{\hspace{0mm}}}{Branch-and-Bound} & \multicolumn{6}{c@{\hspace{0mm}}}{Branch-and-Bound with \eqref{prob.linearalgebraicbound}}\\
\cmidrule(l){3-8} \cmidrule(l){9-14} &  & UB & Obj. & Viol. & Nodes & Gap (\%) & T(s) & UB & Obj. & Viol. & Nodes & Gap (\%) & T(s) \\\midrule
2 &	2 &	0.295 & \textbf{0.295} &	0 &	5100 &	0.00 & 10.9 &	0.295 & \textbf{0.295} &	0 &	3557 &	0.00 & 19.42\\
2 &	4 &	0.404 & \textbf{0.404} &	0 &	99800 &	0.01 & 59.49 &	0.404 & \textbf{0.404} &	0 &	109485 &	0.01 & 238.48\\
2 &	6 &	0.514 & \textbf{0.456} &	0 &	1405600 &12.72	 & $>600$ &	0.521 & 0.452 &	0 &	1038059 &15.20	 & $>600$\\
2 &	8 &	0.595 & 0.467 &	0 &	1266600 & 27.28 &	$>600$ &	0.604 & 0.465 &	0 &	1413722 &29.70 &	$>600$\\
2 &	10 &	0.633 & {0.486} &	0 &	640900 &	30.13 & $>600$ &	0.635 & 0.488 &	0 &	848164 & 30.03 &	$>600$\\\midrule
3 &	2 &	0.435 & \textbf{0.435} &	0 &	22400 &	0.01 & 17.92 &	0.435 & \textbf{0.435} &	0 &	17924 &	0.00 & 176.79\\
3 &	4 &0.717 & 	0.530 &	0 &	542700 &	35.43 & $>600$ &	0.753 & 0.524 &	0 &	222894 &	43.68 & $>600$\\
3 &	6 & 0.846 & 0.551 &	0 &	518100 &	53.62 & $>600$ &	0.879 & 0.560 &	0 &	387403 & 56.96 &	$>600$\\
3 &	8 &	0.933 & 0.585 &	0 &	564200 &	59.45 & $>600$ &	0.935 & 0.570 &	0 &	612382 &	64.06 & $>600$\\
3 &	10 & 0.962 & 0.627 &	0 &	383100 &	53.51 & $>600$ &	0.963 & 0.595 &	0 &	455231 & 61.95 &	$>600$\\\midrule
4 &	2 &	0.566 & \textbf{0.554} &	0 &	783600 &	2.02 & $>600$ &	0.774 & \textbf{0.554} &	0 &	36000 & 39.64 &	$>600$\\
4 &	4 &	1.089 & {0.636} &	0 &	269900 &	71.24 & $>600$ &	1.114 & 0.610 &	0 &	203231 & 82.70 & 	$>600$\\
4 &	6 &	1.213 & {0.642} &	0 &	268300 &	88.92 & $>600$ &	1.226 & 0.617 &	0 &	201271 &	98.73 & $>600$\\
4 &	8 & 1.267 & 0.648 &	0 &	251800 &	95.45 & $>600$ &	1.266 & 0.699 &	0 &	308274 &	81.18 & $>600$\\
4 &	10 & 1.289 & {0.713} &	0 &	163200 &	80.67 & $>600$ &	1.289 & 0.714 &	0 &	266105 &	80.62 & $>600$\\\midrule
5 &	2 &	0.946 & {0.641} &	0 &	702700 &	47.60 & $>600$ &	1.073 & 0.603 &	0 &	44795 &	77.98 & $>600$\\
5 &	4 & 1.430 & 0.697 &	0 &	215200 &	105.20 & $>600$ &	1.468 & 0.652 &	0 &	113951 &	125.31 & $>600$\\
5 &	6 &	1.555 & 0.692 &	0 &	199400 &	124.65 & $>600$ &	1.555 & 0.686 &	0 &	152744 &	126.63 & $>600$\\
5 &	8 &	1.592 & 0.761 &	0 &	191100 &	109.14 & $>600$ &	1.610 & 0.714 &	0 &	156767 &	125.47 & $>600$\\
5 &	10 & 1.616 & 0.801 &	0 &	147500 &	101.68 & $>600$ &	1.620 & 0.804 &	0 &	108311 &	101.51 & $>600$\\\midrule
6 &	2 &	1.323 & 0.702 &	0 &	323700 &	88.34 & $>600$ &	1.430 & 0.698 &	0 &	43880 &	104.84 & $>600$\\
6 &	4 &	1.822 & 0.761 &	0 &	139900 &	139.43 & $>600$ &	1.807 & 0.711 &	0 &	93572 &	154.35 & $>600$\\
6 &	6 &	1.903 & {0.771} &	0 &	114000 &	146.82 & $>600$ &	1.907 & 0.776 &	0 &	43754 &	145.87 & $>600$\\
6 &	8 & 1.932 & \textbf{0.846} &	0 &	141700 &	128.42 & $>600$ &	1.946 & 0.833 &	0 &	62657 &	133.63 & $>600$\\
6 &	10 & 1.947 & \textbf{0.868} &	0 &	31200 &	124.29 & $>600$ &	1.947 & 0.864 &	0 &	68429 &	125.36 & $>600$\\\midrule
Avg & & 1.113 & {0.623} &	0 &	375700 &	69.04 & 533.28 &	1.138 & 0.613 &	0 &	280500 &	76.22 & 547.62\\
\bottomrule
\end{tabular}
\caption{Performance of branch-and-bound without warmstart on the pitprops dataset ($p=13$) using the experimental setup laid out in Section \ref{ssec:feasiblemethods}, except we use a time limit of $600$s for branch-and-bound. 
We report the performance of branch-and-bound with and without the upper bound developed in Section \ref{ssec:bounds} separately. 
The column ``UB'' reports the upper bound obtained by the branch-and-bound scheme at the time limit. 
We use $>600$ to denote an instance where branch-and-bound terminates at the $600$s time limit. The column gap denotes the relative optimality gap reported by Gurobi at termination (in \%).
We denote the best-performing solution (in terms of the proportion of variance explained minus the total orthogonality constraint violation) in bold (cont.). 
} 
    \label{tab:comparison_gurobi_pitprops1}
\end{table}

\begin{table}[h]
\centering\footnotesize
\begin{tabular}{@{}l l r r r r r r r r r r r r r@{}} \toprule
$r$ &  $k_t$ & \multicolumn{6}{c@{\hspace{0mm}}}{Branch-and-Bound (warm-start)} & \multicolumn{6}{c@{\hspace{0mm}}}{Branch-and-Bound with \eqref{prob.linearalgebraicbound} (warm-start)}\\
\cmidrule(l){3-8} \cmidrule(l){9-14} &  & UB & Obj. & Viol. & Nodes & Gap (\%) & T(s) & UB & Obj. & Viol. & Nodes & Gap (\%) & T(s) \\\midrule
2 &	2 &	0.295 & \textbf{0.295} &	0 &	6400 &	0.00 & 10.90 &	0.295 & \textbf{0.295} &	0 &	143 &	0.00 & 19.42\\
&	4 &	0.404 & \textbf{0.404} &	0 &	75200 &	0.01 & 59.49 &	0.404 & \textbf{0.404} &	0 &	88735 &	0.01 & 238.5\\
& 6 & 0.521 & 0.453 &	0 &	1051700 &	15.10 & $>600$ &	0.524 & 0.445 &	0 &	744000 &17.89 &	$>600$\\
& 8 & 0.598 & 0.463 &	0 &	826700 &	29.31 & $>600$ &	0.604 & 0.462 &	0 &	1447055 &	30.73 & $>600$\\
& 10 &0.633 & 0.489 &	0 &	686100 &	29.55 & $>600$ &	0.636 & 0.487 &	0 &	687944 &	30.63 & $>600$\\\midrule
3 & 2 & 0.435 & \textbf{0.435} &	0 &	48600 &	0.00 & 17.92 &	0.435 & \textbf{0.435} &	0 &	2490 &	0.00 & 176.8\\
& 4 & 0.716 & 0.536 &	0 &	491800 &	33.59 & $>600$ &	0.752 & 0.524 &	0 &	240207 &	43.42 & $>600$\\
& 6 & 0.863 & 0.560 &	0 &	470900 &	54.16 & $>600$ &	0.855 & 0.567 &	0 &	598767 &	50.86 & $>600$\\
& 8 & 0.935 & 0.578 &	0 &	383000 &	61.89 & $>600$ &	0.937 & 0.566 &	0 &	239381 &	65.68 & $>600$\\
& 10 & 0.963 & 0.603 &	0 &	297600 &	59.74 & $>600$ &	0.963 & 0.595 &	0 &	355350 &	61.91 & $>600$\\\midrule
4 & 2 & 0.554 & \textbf{0.554} &	0 &	604900 &	0.01 & $>600$ & 0.554 & 	\textbf{0.544} &	0 &	29467 &	38.09 & $>600$\\
& 4 & 1.107 & 0.627 &	0 &	264300 &	76.52 & $>600$ &	1.126 & 0.646 &	0 &	136767 &	74.26 & $>600$\\
& 6 & 1.208 & 0.633 &	0 &	250300 &	90.89 & $>600$ &	1.220 & 0.625 &	0 &	247504 &	95.30& $>600$\\
& 8 & 1.273 & 0.676 &	0 &	225600 &	88.26 & $>600$ &	1.271 & 0.651 &	0 &	216736 &	95.07& $>600$\\
& 10 & 1.289 & 0.702 &	0 &	210000 &	83.73 & $>600$ &	1.291 & 0.712 &	0 &	165362 &	81.32 & $>600$\\\midrule
5 & 2 & 0.907 & \textbf{0.656} &	0 &	279800 &	38.27 & $>600$ &	1.132 & 0.616 &	0 &	28514 &	83.71 & $>600$\\
& 4 & 1.419 & 0.718 &	0 &	251000 &	97.50 & $>600$ &	1.475 & 0.680 &	0 &	141932 &	116.99 & $>600$\\
& 6 & 1.549 & 0.699 &	0 &	225900 &	121.73 & $>600$ &	1.555 & 0.677 &	0 &	229664 &	129.75 & $>600$\\
& 8 & 1.609 & 0.743 &	0 &	223700 &	116.69 & $>600$ &	1.603 & 0.743 &	0 &	200306 &	115.90 & $>600$\\
& 10 & 1.619 & 0.800 &	0 &	126300 &	102.49 & $>600$ &	1.621 & 0.794 &	0 &	135142 &	103.99 & $>600$\\\midrule
6 & 2 & 1.285 & \textbf{0.749} &	0 &	253200 &	71.53 & $>600$ &	1.647& \textbf{0.749} &	0 &	11818 &	119.79 & $>600$\\
& 4 & 1.821 & 0.778 &	0 &	89900 &	133.89 & $>600$ &	1.859 & 0.737 &	0 &	141476 &	152.33 & $>600$\\
& 6 & 1.891 & 0.780 &	0 &	185600 &	142.52 & $>600$ &	1.906 & 0.768 &	0 &	126387 &	148.36 & $>600$\\
& 8 & 1.927 & 0.833 &	0 &	124300 &	131.24 & $>600$ &	1.942 & 0.856 &	0 &	91093 &	126.86 & $>600$\\
& 10 & 1.948 & 0.867 &	0 &	43200 &	124.74 & $>600$ & 1.948 &	0.866 &	0 &	82759 &	124.83 & $>600$\\\midrule
Avg & & 1.111 & 0.625 & 0 & 307800 & 68.13 & $>600$ & 1.150 & 0.618 & 0 & 255600 & 76.30 & $>600$\\
\bottomrule
\end{tabular}
\caption{Performance of branch-and-bound with warmstart 
on the pitprops dataset ($p=13$) using the experimental setup laid out in Section \ref{ssec:feasiblemethods}, except we use a time limit of $600$s for branch-and-bound. 
We report the performance of branch-and-bound with and without the upper bound developed in Section \ref{ssec:bounds} separately.
We use $>600$ to denote an instance where branch-and-bound terminates at the $600$s time limit. The column gap denotes the relative optimality gap reported by Gurobi at termination (in \%).
We denote the best-performing solution (in terms of the proportion of variance explained minus the orthogonality violation) in bold (cont.). 
} 
\label{tab:comparison_gurobi_pitprops2}
\end{table}

\FloatBarrier
We observe that including the combinatorial upper bound developed in Section \ref{sec.gen.gershgorin} within the branch-and-bound scheme does more harm than good, 
and that using Algorithm \ref{alg:greedymethod2} as a warmstart marginally improves the performance. 

Furthermore, we observe that the upper bound returned by branch-and-bound outperforms the semidefinite upper bound from Problem \eqref{prob:disjunctiverelax_permutationinvariant} for the smallest combinations of $r$ and $k$, but rapidly becomes worse as $k$ and $r$ increases, to the extent that it is unable to provide an upper bound better than the trivial bound of $1$ for the largest combinations of $r$ and $k$. 
This suggests that the upper bound from branch-and-bound is not practically useful for larger problem instances.

\FloatBarrier \newpage
\section{Supplementary Numerical Results on UCI Datasets}\label{sec:ec.numres}
This section provides supplementary results supporting the numerical experiments performed in Section \ref{sec:numres} on UCI datasets.

\subsection{Description of the Experimental Setup}\label{sec:ec.supp_expr_setup}
All experiments were performed on MIT's supercloud cluster \citep{reuther2018interactive}, which hosts Intel Xeon Platinum 8260 processors and Intel Xeon Gold 6248 processors. {\color{black}For experiments where $p < 100$, we use Platinum processors with $32$ GB RAM, for experiments where $p \in [100, 250]$, Platinum processors with $100$ GB RAM; and for $p>250$, Gold processors with $370$ GB RAM.}

We also implement some existing algorithmic strategies from the literature, to provide a baseline for the performance of our methods. To abide by software licensing restrictions, all existing strategies from the literature were benchmarked using a MacBook Pro laptop with a $2.9$GHz $6$-Core Intel i9 CPU, using $16$ GB DDR4 RAM. {Therefore, runtimes are not directly comparable across strategies.}

\subsection{Description of the Data Sources} \label{sec:ec.datadescr}
We perform experiments on eleven datasets from the frequently used UCI database in Sections \ref{ssec:bounds}-\ref{ssec:feasiblemethods} and \ref{ssec:symm}. 
Of the eleven datasets, six datasets are overdetermined (meaning $n >p$), while five datasets are underdetermined (meaning {\color{black}$p>n$}). Moreover, many existing works on sparse PCA report results on similar datasets. For instance, the pitprops dataset was also considered by \cite{jolliffe2003modified,zou2006sparse, journee2010generalized} among others, and three of the datasets studied by \cite{berk2019certifiably} are included within our suite of datasets. Thus, our experimental setup is broadly representative of both {\color{black}the underdetermined and the overdetermined regimes, as well as of the literature}. For completeness, we summarize the datasets we benchmark on and their dimensionality in Table \ref{tab:dataset_summary}.

\begin{table}[h!]
    \centering\footnotesize
    \begin{tabular}{l r r l}
    \toprule
    Dataset & $p$ & $n$ \\\midrule
    Pitprops & $13$ & $180$ \\
    Wine & $13$ & $178$ \\
    Ionosphere & $34$ & $351$\\
    Lung (Lung cancer) & $54$ & $32$ \\
    Geographical (Geographical Origin of Music) & $68$ & $1059$\\
    Communities (Communities and Crime) & $101$ & $1994$ \\
    Arrhythmia & $274$ & $452$ \\
    Voice (LSVT Voice Rehabilitation) & $310$ & $126$\\
    Gait (Gait Classification) & $320$ & $48$\\
    Gastro (Gastrointestinal Lesions in Regular Colonoscopy) & $466$ & $152$\\
    Micromass & $1300$ & $931$\\
    \bottomrule
    \end{tabular}
    \caption{
    Summary of the 11 datasets in our library, where $n$ denotes the number of observations and $p$ the number of features. For conciseness, the names of certain datasets are abbreviated throughout. For these datasets, we first state the abbreviation used, followed by their full names in brackets. Further, we report the dimensionality of each dataset after preprocessing, {\color{black}removing} all features with missing values. All datasets can be found in the UCI database, except the pitprops dataset, which is due to \citet{jeffers1967two} and distributed via the R package ElasticNet.}
    \label{tab:dataset_summary}
\end{table}

\FloatBarrier 
\subsection{Preliminary Experiments With Pitprops Dataset}\label{ssec:ec.pitprops}

We now provide instance-wise results for different variants of our methods on the \verb|pitprops| dataset. In particular, we consider invoking the valid inequalities \eqref{prob.gen.gershgorin.formulation} derived in Section \ref{sec.gen.gershgorin} to improve branch-and-bound further. When we do so, we also invoke a branching callback each time we expand a node to determine whether the subtree rooted at this node can improve upon the incumbent solution. {\color{black}This is justified by the fact that at each node, some variables $Z_{i,t}$ are fixed to $0$, some to $1$, and some are not fixed.} Accordingly, we can compute an upper bound on any solution with the same fixed variables by relaxing the orthogonality constraint and applying the Gershgorin circle theorem to each component separately; see \citet[Section 2.4]{bertsimas2020solving} for a discussion of this callback in the rank-one case. In particular, if the Gershgorin bound for a given subtree is weaker than an incumbent solution, then this subtree does not contain any optimal solutions, and we can prune it from our search tree. 

\FloatBarrier
\subsection{Instance-Wise Results on Larger UCI Datasets} \label{ssec:largeruci}
Next, we provide an instance-by-instance account of the results summarized in Table \ref{tab:ucilargersummarized}--Table \ref{tab:ucilargersummarized2}, in Tables \ref{tab:comparison_feasiblemethods_uci}-\ref{tab:comparison_feasiblemethods_uci3_part4}.

\begin{table}[h]
\centering\footnotesize
\begin{tabular}{@{}l l l l r r r r r r r r r r r@{}} \toprule
Dataset & $p$ & $r$ &  $k_t$ & \multicolumn{4}{c@{\hspace{0mm}}}{Alg. \ref{alg:greedymethod2}} & \multicolumn{4}{c@{\hspace{0mm}}}{Alg. \ref{alg:alternatingmin} } & \multicolumn{3}{c@{\hspace{0mm}}}{Branch-and-bound}  \\
\cmidrule(l){5-8} \cmidrule(l){9-12} \cmidrule(l){13-15} &  &  & & UB & Obj. & Viol. & T(s) & UB & Obj. & Viol. & T(s) & Obj. & Viol. & T(s)  \\\midrule
Pitprops &	13 &	2 &	5 &	0.449 &	0.429 &	0 &	0.99 &	0.524 &{0.433} &	0 &	5.01 &	\textbf{0.439} &	0 &	$>7200$ \\
&	&	2 &	10 &	0.507 &	0.380 &	0 &	0.34 &	0.642 & \textbf{0.500} &	0 &	4.10 &	{0.498} &	0 &	$>7200$ \\
&	&	3 &	5 &	0.616 &	0.541 &	0 &	0.44 &	0.786& \textbf{0.582} &	0 &	6.48 &	{0.555} &	0 &	$>7200$ \\
&	&	3 &	10 &	0.652 &	0.511 &	0 &	0.46 &	0.963 & \textbf{0.650} &	0 &	6.41 &	{0.618} &	0 &	$>7200$ \\\midrule
Wine &	13 &	2 &	5 &	0.458 &	0.401 &	0 &	0.35 &	0.529 & {0.446} &	0 &	2.19 &	\textbf{0.448} &	0 &	1073 \\
&	&	2 &	10 &	0.554 &	0.508 &	0 &	0.31 &	0.707 & \textbf{0.544} &	0 &	4.03 &	{0.537} &	0 &	$>7200$ \\
&	&	3 &	5 &	0.632 &	0.446 &	0 &	0.48 &	0.794 & \textbf{0.613} &	0 &	3.56 &	{0.576} &	0 &	$>7200$ \\
&	&	3 &	10 &	0.665 &	0.528 &	0 &	0.46 &	1.060 & \textbf{0.660} &	0 &	6.47 &	{0.640} &	0 &	$>7200$ \\\midrule
Ionosphere &	34 &	2 &	5 &	0.209 &	0.203 &	0 &	5.97 &	0.221 & \textbf{0.204} &	0 &	3.18 &	{0.202} &	0 &	$>7200$ \\
&	&	2 &	10 &	0.305 &	0.265 &	0 &	8.72 &	0.361 & \textbf{0.285} &	0 &	3.95 &	\textbf{0.285} &	0 &	$>7200$ \\
&	&	2 &	20 &	0.378 &	0.286 &	0 &	23.31 &	0.500 & \textbf{0.360} &	0 &	7.64 &	{0.344} &	0 &	$>7200$ \\
&	&	3 &	5 &	0.297 &	0.287 &	0 &	31.04 &	0.331 & {0.279} &	0 &	6.72 &	\textbf{0.289} &	0 &	$>7200$ \\
&	&	3 &	10 &	0.411 &	0.305 &	0 &	39.00 &	0.542 & {0.397} &	0 &	11.89 &	{0.375} &	0 &	$>7200$ \\
&	&	3 &	20 &	0.464 &	0.390 &	0 &	10.26 &	0.749 & \textbf{0.458} &	0 &	12.17 &	{0.383} &	0 &	$>7200$ \\\midrule
Lung &	54 &	2 &	5 &	0.119 &	\textbf{0.119} &	0 &	17.84 &	0.124 & {0.110} &	0 &	3.81 &	{0.113} &	0 &	$>7200$ \\
&	&	2 &	10 &	0.176 &	\textbf{0.175} &	0 &	30.81 &	0.178 & {0.171} &	0 &	2.06 &	{0.168} &	0 &	$>7200$ \\
&	&	2 &	20 &	0.234 &	{0.170} &	0 &	30.82 &	0.262 & {0.217} &	0 &	4.12 &	{0.185} &	0 &	$>7200$ \\
&	&	3 &	5 &	0.173 &	0.165 &	0 &	40.83 &	0.185 & {0.160} &	0 &	4.45 &	{0.169} &	0 &	$>7200$ \\
&	&	3 &	10 &	0.249 &	0.194 &	0 &	90.20 &	0.267 & {0.240} &	0 &	4.31 &	{0.188} &	0 &	$>7200$ \\
&	&	3 &	20 &	0.324 &	0.308 &	0 &	34.50 &	0.393 & {0.303} &	0 &	9.51 &	{0.213} &	0 &	$>7200$ \\\midrule
Geography &	68 &	2 &	5 &	0.147 &	0.145 &	0 &	99.45 &	0.147 & \textbf{0.147} &	0 &	2.99 &	{0.145} &	0 &	$>7200$ \\
&	&	2 &	10 &	0.294 &	0.290 &	0 &	107.27 &	0.294 & \textbf{0.294} &	0 &	1.81 &	{0.292} &	0 &	$>7200$ \\
&	&	2 &	20 &	0.433 &	0.393 &	0 &	1213.3 &	0.564 & {0.376} &	0 &	6.32 &	{0.327} &	0 &	$>7200$ \\
&	&	3 &	5 &	0.221 &	0.213 &	0 &	119.3 &	0.221 & \textbf{0.221} &	0 &	2.58 &	{0.215} &	0 &	$>7200$ \\
&	&	3 &	10 &	0.410 &	0.342 &	0 &	1453.19 &	0.441 & {0.348} &	0 &	5.39 &	{0.355} &	0 &	$>7200$ \\
&	&	3 &	20 &	0.529 &	0.457 &	0 &	1571.86 &	0.846 & {0.345} &	0 &	11.51 &	{0.352} &	0 &	$>7200$ \\\midrule
Communities &	101 &	2 &	5 &	0.095 &	\textbf{0.095} &	0 &	484.0 &	0.096 & {0.078} &	0 &	3.73 &	\textbf{0.095} &	0 &	$>7200$ \\
&	&	2 &	10 &	0.169 &	\textbf{0.169} &	0 &	1327 &	0.175 & {0.159} &	0 &	4.81 &	{0.160} &	0 &	$>7200$ \\
&	&	2 &	20 &	0.268 &	0.219 &	0 &	2438 &	0.284 & {0.244} &	0 &	6.26 &	{0.198} &	0 &	$>7200$ \\
&	&	3 &	5 &	0.141 &	\textbf{0.141} &	0 &	979.7 &	0.144 & {0.119} &	0 &	7.63 &	\textbf{0.141} &	0 &	$>7200$ \\
&	&	3 &	10 &	0.246 &	0.242 &	0 &	3553 &	0.262 & {0.243} &	0 &	8.92 &	{0.205} &	0 &	$>7200$ \\
&	&	3 &	20 &	0.385 &	0.267 &	0 &	3231 &	0.425 & \textbf{0.370} &	0 &	8.28 &	{0.300} &	0 &	$>7200$ \\\midrule
Arrhythmia &	274 &	2 &	5 &	0.031 &	0.021 &	0 &	583.7 &	0.031 & {0.027} &	0 &	16.57 &	{0.027} &	0 &	$>7200$ \\
&	&	2 &	10 &	0.055 &	0.035 &	0 &	555.4 &	0.055 & {0.047} &	0 &	24.72 &	{0.044} &	0 &	$>7200$ \\
&	&	2 &	20 &	0.086 &	0.067 &	0 &	622.8 &	0.084 & {0.071} &	0.002 &	49.41 &	{0.059} &	0 &	$>7200$ \\
&	&	3 &	5 &	0.047 &	0.031 &	0 &	1423.0 &	0.046 & {0.039} &	0 &	27.99 &	{0.044} &	0 &	$>7200$ \\
&	&	3 &	10 &	0.083 &	0.044 &	0 &	1085.8 &	0.083 & {0.067} &	0 &	38.16 &	{0.065} &	0 &	$>7200$ \\
&	&	3 &	20 &	0.129 &	0.083 &	0 &	1059.7 &	0.126 & {0.105} &	0 &	28.96 &	{0.068} &	0 &	$>7200$ \\
\bottomrule
\end{tabular}
\caption{Performance of Algorithms \ref{alg:greedymethod2}--\ref{alg:alternatingmin} and branch-and-bound on UCI datasets. $k_t$ denotes the sparsity of each individual component, meaning a set of $r$ PCs have a collective sparsity budget of $k_t r$. Note that all objective values are reported in terms of the proportion of correlation explained by dividing by $p$, the number of features.
} 
\label{tab:comparison_feasiblemethods_uci}
\end{table}

\begin{table}[h]
\centering\footnotesize
\begin{tabular}{@{}l l l l r r r r r r r r r r r@{}} \toprule
Dataset & $p$ & $r$ &  $k_t$ & \multicolumn{4}{c@{\hspace{0mm}}}{Alg. \ref{alg:greedymethod2}} & \multicolumn{4}{c@{\hspace{0mm}}}{Alg. \ref{alg:alternatingmin} } & \multicolumn{3}{c@{\hspace{0mm}}}{Branch-and-bound}  \\
\cmidrule(l){5-8} \cmidrule(l){9-12} \cmidrule(l){13-15} &  &  & & UB & Obj. & Viol. & T(s) & UB & Obj. & Viol. & T(s) & Obj. & Viol. & T(s)  \\\midrule
Voice &	310 &	2 &	5 &	0.032 &	0.024 &	0 &	375.0 &	0.032 & \textbf{0.032} &	0 &	21.14 &	\textbf{0.032} &	0 &	$>7200$ \\
&	&	2 &	10 &	0.064 &	\textbf{0.064} &	0 &	741.2 &	0.064 & {0.063} &	0 &	21.99 &	\textbf{0.064} &	0 &	$>7200$ \\
&	&	2 &	20 &	0.127 &	\textbf{0.127} &	0 &	630.6 &	0.127 & {0.124} &	0 &	20.97 &	{0.109} &	0 &	$>7200$ \\
&	&	3 &	5 &	0.048 &	\textbf{0.048} &	0 &	758.3 &	0.048 & {0.047} &	0 &	29.94 &	\textbf{0.048} &	0 &	$>7200$ \\
&	&	3 &	10 &	0.096 &	0.079 &	0 &	797.1 &	0.096 & {0.093} &	0 &	34.52 &	\textbf{0.096} &	0 &	$>7200$ \\
&	&	3 &	20 &	0.191 &	0.125 &	0 &	721.5 &	0.191 & {0.183} &	0 &	31.06 &	{0.155} &	0 &	$>7200$ \\\midrule
Gait &	320 &	2 &	5 &	0.031 &	0.018 &	0 &	399.2 &	0.031 & {0.028} &	0 &	19.40 &	{0.028} &	0 &	$>7200$ \\
&	&	2 &	10 &	0.057 &	0.036 &	0 &	450.7 &	0.057 & {0.050} &	0 &	24.64 &	{0.047} &	0 &	$>7200$ \\
&	&	2 &	20 &	0.103 &	0.062 &	0 &	392.3 &	0.103 & {0.081} &	0 &	24.51 &	{0.067} &	0 &	$>7200$ \\
&	&	3 &	5 &	0.046 &	0.036 &	0 &	933.7 &	0.046 & {0.041} &	0 &	27.80 &	\textbf{0.045} &	0 &	$>7200$ \\
&	&	3 &	10 &	0.085 &	0.049 &	0 &	829.3 &	0.085 & {0.077} &	0 &	34.99 &	{0.060} &	0 &	$>7200$ \\
&	&	3 &	20 &	0.154 &	0.060 &	0 &	792.1 &	0.154 & {0.121} &	0 &	34.11 &	{0.111} &	0 &	$>7200$ \\\midrule
Gastro &	466 &	2 &	5 &	{0.021} &	\textbf{0.021} &	0 &	1633 &	0.021 & \textbf{0.021} &	0 &	452.9 &	\textbf{0.021} &	0 &	$>7200$ \\
&	&	2 &	10 &	0.043 &	\textbf{0.043} &	0 &	1753 &	0.043 & \textbf{0.043} &	0 &	38.95 &	\textbf{0.043} &	0 &	$>7200$ \\
&	&	2 &	20 &	0.086 &	0.085 &	0 &	2166 &	0.086 & {0.085} &	0 &	535.4 &	\textbf{0.086} &	0 &	$>7200$ \\
&	&	3 &	5 &	0.032 &	\textbf{0.032} &	0 &	2682 &	0.032 & \textbf{0.032} &	0 &	59.97 &	\textbf{0.032} &	0 &	$>7200$ \\
&	&	3 &	10 &	0.064 &	\textbf{0.064} &	0 &	4307 &	0.064 & \textbf{0.064} &	0 &	73.39 &	\textbf{0.064} &	0 &	$>7200$ \\
&	&	3 &	20 &	0.129 &	0.126 &	0 &	5544 &	0.128 & {0.121} &	0 &	529.8 &	\textbf{0.128} &	0 &	$>7200$ \\\midrule
Micromass &	1300 &	2 &	5 &	0.008 &	0.005 &	0 &	1089 &	0.008 & {0.006} &	0 &	170.5 &	{0.004} &	0 &	$>7200$ \\
&	&	2 &	10 &	0.015 &	0.008 &	0 &	13620 &	0.014 & {0.011} &	0 &	158.6 &	{0.011} &	0 &	6100 \\
&	&	2 &	20 &	0.027 &	0.018 &	0 &	9213 &	0.023 & 0.019 & 0 &	440.3 &	{0.018} &	0 &	6826 \\
&	&	3 &	5 &	0.012 &	0.008 &	0 &	7953 &	0.011 & {0.010} &	0 &	238.0 &	{0.006} &	0 &	$>7200$ \\
&	&	3 &	10 &	0.023 &	0.009 &	0 &	19640 &	0.021 & {0.018} &	0 &	208.5 &	{0.009} &	0 &	$>7200$ \\
&	&	3 &	20 &	0.043 &	0.029 &	0 &	18630 &	0.034 & {0.030} &	0 &	205.9 &	{0.010} &	0 &	$>7200$ \\\midrule
Avg & & & & 0.213 & 0.176 & 0.000 & 1899 & 0.257 & 0.199 & 0.000 & 61.38 & 0.187 & 0 & $>7200$ \\
\bottomrule
\end{tabular}
\caption{Performance of Algorithms \ref{alg:greedymethod2}--\ref{alg:alternatingmin} and branch-and-bound on UCI datasets (cont). $k_t$ denotes the sparsity of each individual component, meaning a set of $r$ PCs {\color{black}has} a collective sparsity budget of $k_t r$. Note that all objective values are reported in terms of the proportion of correlation explained by dividing by $p$, the number of features.
} 
\label{tab:comparison_feasiblemethods_uci_part2}
\end{table}

\begin{table}[h!]
\centering\footnotesize
\begin{tabular}{@{}l l l l r r r r r r r r r r r @{}} \toprule
Dataset & $p$ & $r$ &  $k_t$ & \multicolumn{3}{c@{\hspace{0mm}}}{\citet{berk2019certifiably}} & \multicolumn{3}{c@{\hspace{0mm}}}{\citet{hein2010inverse}} & \multicolumn{3}{c@{\hspace{0mm}}}{\citet{zou2006sparse}}  \\
\cmidrule(l){5-7} \cmidrule(l){8-10} \cmidrule(l){11-13}  &  &  & & Obj. & Viol. & T(s) & Obj. & Viol. & T(s) & Obj. & Viol. & T(s) \\\midrule
Pitprops &	13 &	2 &	5 &	0.421 &	0.168 &	1.67 &	0.418 &	0 &	0.11 &	0.177 &	1.341 &	0.12 \\
&	&	2 &	10 &	0.502 &	0.008 &	0.14 &	0.502 &	0.008 &	0.01 &	0.139 &	1.827 &	0.22 \\
&	&	3 &	5 &	0.592 &	0.675 &	0.08 &	0.575 &	0.166 &	0.02 &	0.169 &	3.462 &	0.04 \\
&	&	3 &	10 &	0.648 &	0.073 &	0.07 &	0.647 &	0.084 &	0 &	0.181 &	3.771 &	0.36 \\\midrule
Wine &	13 &	2 &	5 &	\textbf{0.448} &	0 &	0.04 &	0.422 &	0.004 &	0.01 &	0.127 &	0.315 &	0.04 \\
&	&	2 &	10 &	0.545 &	0.020 &	0.04 &	0.545 &	0.02 &	0 &	0.068 &	0.731 &	0.06 \\
&	&	3 &	5 &	0.610 &	0.019 &	0.06 &	0.559 &	0.092 &	0 &	0.225 &	2.830 &	0.05 \\
&	&	3 &	10 &	0.654 &	0.059 &	0.06 &	0.655 &	0.093 &	0 &	0.232 &	2.771 &	0.32 \\\midrule
Ionosphere &	34 &	2 &	5 &	\textbf{0.205} &	0 &	0.08 &	0.153 &	0 &	0.08 &	0.078 &	0 &	0.02 \\
&	&	2 &	10 &	{0.289} &	0 &	0.30 &	0.288 &	0 &	0.01 &	0.106 &	0 &	0.04 \\
&	&	2 &	20 &	0.369 &	0.058 &	4.45 &	\textbf{0.370} &	0.010 &	0.17 &	0.147 &	0.305 &	0.12 \\
&	&	3 &	5 &{0.291} &	0 &	0.14 &	0.227 &	0 &	0.02 &	0.097 &	1.666 &	0.07 \\
&	&	3 &	10 &	0.392 &	0.109 &	0.38 &	0.365 &	0.255 &	0.01 &	0.100 &	1.909 &	0.12 \\
&	&	3 &	20 &	0.449 &	0.183 &	0.27 &	0.451 &	0.037 &	0.03 &	0.111 &	2.111 &	4.51 \\\midrule
Lung &	54 &	2 &	5 &	\textbf{0.119} &	0 &	0.43 &	0.107 &	0 &	0.34 &	0.040 &	0.587 &	0.04 \\
&	&	2 &	10 &	\textbf{0.176} &	0 &	0.10 &	0.170 &	0 &	0.03 &	0.044 &	0.639 &	0.12 \\
&	&	2 &	20 &	0.220 &	0.008 &	0.44 &	0.184 &	0 &	0.05 &	0.044 &	0.908 &	0.63 \\
&	&	3 &	5 &	\textbf{0.172} &	0 &	0.10 &	0.149 &	0 &	0.03 &	0.061 &	2.755 &	0.15 \\
&	&	3 &	10 &	\textbf{0.243} &	0 &	0.11 &	0.234 &	0.113 &	0.05 &	0.054 &	1.593 &	0.20 \\
&	&	3 &	20 &	0.300 &	0.219 &	0.16 &	0.261 &	0.081 &	0.04 &	0.044 &	1.703 &	1.20 \\\midrule
Geography &	68 &	2 &	5 &	\textbf{0.147} &	0 &	0.09 &	0.097 &	0 &	0.01 &	0.034 &	1.793 &	0.41 \\
&	&	2 &	10 &	\textbf{0.294} &	0 &	0.08 &	0.164 &	0 &	0 &	0.068 &	1.939 &	0.66 \\
&	&	2 &	20 &	{0.395} &	0 &	5.95 &	0.316 &	0.135 &	0.04 &	0.062 &	1.754 &	0.53 \\
&	&	3 &	5 &	\textbf{0.221} &	0 &	0.13 &	0.122 &	0 &	0.01 &	0.061 &	2.720 &	0.57 \\
&	&	3 &	10 &	\textbf{0.389} &	0 &	0.18 &	0.192 &	0 &	0.01 &	0.054 &	4.021 &	0.90 \\
&	&	3 &	20 &	0.484 &	0.273 &	23.66 &	0.387 &	0.261 &	0.06 &	0.090 &	5.009 &	1.40 \\\midrule
Communities &	101 &	2 &	5 &	\textbf{0.095} &	0 &	0.73 &	0.093 &	0 &	0 &	0.032 &	0.576 &	0.05 \\
&	&	2 &	10 &	\textbf{0.169} &	0 &	1.68 &	0.154 &	0 &	0 &	0.029 &	0.605 &	0.18 \\
&	&	2 &	20 &	\textbf{0.258} &	0 &	120 &	\textbf{0.258} &	0 &	0.07 &	0.027 &	0.090 &	1.49 \\
&	&	3 &	5 &	\textbf{0.141} &	0 &	1.18 &	0.129 &	0 &	0.01 &	0.050 &	1.854 &	0.29 \\
&	&	3 &	10 &	\textbf{0.245} &	0 &	2.85 &	0.181 &	0 &	0.02 &	0.044 &	1.504 &	1.76 \\
&	&	3 &	20 &	{0.361} &	0.058 &	180.1 &	0.350 &	0.064 &	0.02 &	0.043 &	1.869 &	5.27 \\\midrule
Arrhythmia &	274 &	2 &	5 &	\textbf{0.031} &	0 &	2.81 &	0.012 &	0 &	0.02 &	0.007 &	1.799 &	0.71 \\
&	&	2 &	10 &	\textbf{0.052} &	0 &	61.25 &	0.011 &	0 &	0.03 &	0.007 &	1.143 &	1.08 \\
&	&	2 &	20 &	\textbf{0.077} &	0 &	120.0 &	0.043 &	0.005 &	0.06 &	0.006 &	1.140 &	4.62 \\
&	&	3 &	5 &	\textbf{0.046} &	0 &	5.35 &	0.016 &	0 &	0.02 &	0.012 &	1.076 &	0.53 \\
&	&	3 &	10 &	\textbf{0.074} &	0 &	121.6 &	0.018 &	0 &	0.05 &	0.012 &	0.876 &	3.82 \\
&	&	3 &	20 &	\textbf{0.109} &	0 &	180.0 &	0.074 &	0.005 &	0.07 &	0.012 &	0.694 &	10.65 \\
\bottomrule
\end{tabular}
\caption{Performance of the methods of {\citet{berk2019certifiably}}, {\citet{hein2010inverse}}, and {\citet{zou2006sparse}} on UCI datasets.
} \label{tab:comparison_feasiblemethods_uci2_part3}
\end{table}

\begin{table}[h!]
\centering\footnotesize
\begin{tabular}{@{}l l l l r r r r r r r r r r r @{}} \toprule
Dataset & $p$ & $r$ &  $k_t$ & \multicolumn{3}{c@{\hspace{0mm}}}{\citet{berk2019certifiably}} & \multicolumn{3}{c@{\hspace{0mm}}}{\citet{hein2010inverse}} & \multicolumn{3}{c@{\hspace{0mm}}}{\citet{zou2006sparse}}  \\
\cmidrule(l){5-7} \cmidrule(l){8-10} \cmidrule(l){11-13}  &  &  & & Obj. & Viol. & T(s) & Obj. & Viol. & T(s) & Obj. & Viol. & T(s) \\\midrule
Voice &	310 &	2 &	5 &	\textbf{0.032} &	0 &	1.04 &	\textbf{0.032} &	0 &	0.05 &	0.006 &	0.874 &	0.71 \\
&	&	2 &	10 &	\textbf{0.064} &	0 &	1.09 &	\textbf{0.064} &	0 &	0.04 &	0.006 &	0.907 &	2.68 \\
&	&	2 &	20 &	\textbf{0.127} &	0 &	0.66 &	\textbf{0.127} &	0 &	0.02 &	0.006 &	1.017 &	16.02 \\
&	&	3 &	5 &	\textbf{0.048} &	0 &	1.72 &	0.039 &	0 &	0.16 &	0.009 &	1.834 &	0.94 \\
&	&	3 &	10 &	\textbf{0.096} &	0 &	1.56 &	0.069 &	0 &	0.08 &	0.012 &	0.242 &	12.8 \\
&	&	3 &	20 &	\textbf{0.190} &	0 &	1.22 &	0.187 &	0 &	0.12 &	0.021 &	2.121 &	26.18 \\\midrule
Gait &	320 &	2 &	5 &	\textbf{0.030} &	0 &	0.93 &	0.027 &	0 &	0.02 &	0.006 &	1.071 &	1.75 \\
&	&	2 &	10 &	\textbf{0.055} &	0 &	0.61 &	0.051 &	0 &	0.05 &	0.004 &	1.054 &	1.27 \\
&	&	2 &	20 &	0.094 &	0 &	1.36 &	0.080 &	0 &	0.080 &	0.005 &	0.852 &	3.88 \\
&	&	3 &	5 &	\textbf{0.045} &	0 &	1.02 &	0.041 &	0 &	0.06 &	0.01 &	1.16 &	4.81 \\
&	&	3 &	10 &	\textbf{0.082} &	0 &	1.64 &	0.070 &	0 &	0.06 &	0.009 &	1.821 &	5 \\
&	&	3 &	20 &	\textbf{0.135} &	0 &	1.94 &	0.095 &	0 &	0.15 &	0.008 &	1.477 &	9.49 \\\midrule
Gastro &	466 &	2 &	5 &	\textbf{0.021} &	0 &	2.71 &	0.020 &	0 &	0.05 &	0.007 &	1.154 &	1.14 \\
&	&	2 &	10 &	\textbf{0.043} &	0 &	1.42 &	0.039 &	0 &	0.08 &	0.007 &	0.178 &	1.62 \\
&	&	2 &	20 &	\textbf{0.086} &	0 &	1.95 &	0.076 &	0 &	0.07 &	0.005 &	0.456 &	3.61 \\
&	&	3 &	5 &	\textbf{0.032} &	0 &	2.79 &	0.029 &	0 &	0.05 &	0.006 &	2.303 &	1.59 \\
&	&	3 &	10 &	\textbf{0.064} &	0 &	3.51 &	0.053 &	0 &	0.12 &	0.007 &	1.826 &	3.33 \\
&	&	3 &	20 &	\textbf{0.128} &	0 &	2.81 &	0.085 &	0 &	0.24 &	0.008 &	1.139 &	46.05 \\
\midrule
Micromass &	1300 &	2 &	5 &	\textbf{0.008} &	0 &	45.4 &	0.004 &	0 &	0.77 &	0.002 &	0.014 &	18.05 \\
&	&	2 &	10 &	\textbf{0.014} &	0 &	120.2 &	0.007 &	0 &	1.09 &	0.002 &	0.323 &	41.97 \\
&	&	2 &	20 &	\textbf{0.023} &	0 &	120.2 &	0.012 &	0 &	1.34 &	0.002 &	0.361 &	2.64 \\
&	&	3 &	5 &	\textbf{0.011} &	0 &	71.97 &	0.005 &	0 &	1.54 &	0.002 &	3.004 &	24.30 \\
&	&	3 &	10 &	\textbf{0.020} &	0 &	180.3 &	0.008 &	0 &	1.80 &	0.002 &	2.301 &	31.89 \\
&	&	3 &	20 &	\textbf{0.034} &	0 &	180.3 &	0.013 &	0 &	2.64 &	0.002 &	1.252 &	2.02 \\\midrule
Avg & & &  & 0.205 & 0.031 & 25.57 & 0.180 & 0.023 & 0.19 & 0.049 & 1.458 & 4.95\\
\bottomrule
\end{tabular}
\caption{Performance of the methods of {\citet{berk2019certifiably}}, {\citet{hein2010inverse}}, and {\citet{zou2006sparse}} on UCI datasets (cont.).} 
\label{tab:comparison_feasiblemethods_uci2_part4}
\end{table}

\begin{table}[h!]
\centering\footnotesize
\begin{tabular}{@{}l l l l r r r r r r r@{}} \toprule
Dataset & $p$ & $r$ &  $k_t$ & \multicolumn{3}{c@{\hspace{0mm}}}{\citet{deshpande2014sparse}} & \multicolumn{4}{c@{\hspace{0mm}}}{Algorithm \ref{alg:disjoint.linalg}}  \\
\cmidrule(l){5-7} \cmidrule(l){8-11}&  &  & & Obj. & Viol. & T(s) & UB & Obj. & Viol. & T(s) \\\midrule
Pitprops &	13 &	2 &	5 &	0.422 &	0.226 &	1.00 & 0.559 & 0.422 & 0 & 0.53\\
&	&	2 &	10 &	0.501 &	0.104 &	0.00  & 0.803 & 0.456 & 0 & 0.05\\
&	&	3 &	5 &	0.592 &	0.661 &	0.00  & 0.827 & 0.568 & 0 & 0.13\\
&	&	3 &	10 &	0.644 &	0.214 &	0.00  & 1.198 & 0.569 & 0 & 0.13\\\midrule
Wine &	13 &	2 &	5 &	0.434 &	0.211 &	0.02  & 0.579 & 0.447 & 0 & 0.02\\
&	&	2 &	10 &	0.545 &	0.021 &	0  & 0.876 & 0.508 & 0 & 0.06\\
&	&	3 &	5 &	0.546 &	0.510 &	0  & 0.853 & 0.577 & 0 & 0.21\\
&	&	3 &	10 &	0.656 &	0.228 &	0  & 1.296 & 0.580 & 0 & 0.23\\\midrule
Ionosphere &	34 &	2 &	5 &	0.203 &	0 &	0.12  & 0.228 & 0.202 & 0 & 0.08\\
&	&	2 &	10 &	\textbf{0.287} &	0 &	0.02  & 0.401 & \textbf{0.290} & 0 & 0.14\\
&	&	2 &	20 &	0.368 &	0.011 &	0.02  & 0.618 & 0.357 & 0 & 0.2\\
&	&	3 &	5 &	0.276 &	0 &	0.02  & 0.340 & \textbf{0.292} & 0 & 0.28\\
&	&	3 &	10 &	0.357 &	0.143 &	0.02  & 0.597 & \textbf{0.398} & 0 & 1.51\\
&	&	3 &	20 &	0.447 &	0.166 &	0.02  & 0.920 & 0.408 & 0 & 2.32\\\midrule
Lung &	54 &	2 &	5 &	0.117 &	0 &	0.45  & 0.139 & 0.118 & 0 & 0.71\\
&	&	2 &	10 &	0.154 &	0.291 &	0.04  & 0.218 & 0.171 & 0 & 0.1\\
&	&	2 &	20 &	0.217 &	0.495 &	0.03  & 0.345 & 0.216 & 0 & 0.12\\
&	&	3 &	5 &	0.164 &	0 &	0.03  & 0.204 & 0.169 & 0 & 0.23\\
&	&	3 &	10 &	0.200 &	0.558 &	0.03  & 0.326 & 0.225 & 0 & 0.56\\
&	&	3 &	20 &	0.290 &	0.591 &	0.03  & 0.514 & 0.271 & 0 & 3.54\\\midrule
Geography &	68 &	2 &	5 &	0.130 &	0 &	0.09  & 0.147 & \textbf{0.147} & 0 & 0.64\\
&	&	2 &	10 &	0.294 &	0 &	0.07  & 0.294 & \textbf{0.294} & 0 & 0.81\\
&	&	2 &	20 &	0.395 &	0 &	0.08  & 0.568 & \textbf{0.396} & 0 & 0.61\\
&	&	3 &	5 &	0.184 &	0 &	0.08  & 0.221 & \textbf{0.221} & 0 & 2.79\\
&	&	3 &	10 &	0.378 &	0 &	0.07  & 0.441 & \textbf{0.389} & 0 & 2.04\\
&	&	3 &	20 &	0.462 &	0.073 &	0.08  & 0.852 & 0.479 & 0 & 2.56\\\midrule
Communities &	101 &	2 &	5 &	0.089 &	0 &	0.19  & 0.097 & \textbf{0.095} & 0 & 2.3\\
&	&	2 &	10 &	0.158 &	0 &	0.23  & 0.180 & 0.167 & 0 & 1.37\\
&	&	2 &	20 &	0.249 &	0.145 &	0.19  & 0.320 & \textbf{0.259} & 0 & 1.94\\
&	&	3 &	5 &	0.140 &	0 &	0.19  & 0.146 & \textbf{0.141} & 0 & 1.14\\
&	&	3 &	10 &	0.208 &	0.138 &	0.20  & 0.270 & 0.242 & 0 & 1.61\\
&	&	3 &	20 &	0.317 &	0.669 &	0.19  & 0.476 & 0.369 & 0 & 3.62\\\midrule
Arrhythmia &	274 &	2 &	5 &	0.030 &	0 &	1.79  & 0.032 & \textbf{0.031} & 0 & 13.71\\
&	&	2 &	10 &	0.051 &	0.113 &	1.84  & 0.060 & 0.051 & 0 & 11.26\\
&	&	2 &	20 &	0.074 &	0.019 &	1.88  & 0.105 & 0.075 & 0 & 8.04\\
&	&	3 &	5 &	0.039 &	0 &	1.97  & 0.049 & \textbf{0.046} & 0 & 699.5\\
&	&	3 &	10 &	0.072 &	0 &	2.35  & 0.089 & 0.072 & 0 & 223.1\\
&	&	3 &	20 &	0.103 &	0.243 &	2.23  & 0.155 & 0.107 & 0 & 105.5\\
\bottomrule
\end{tabular}
\caption{Performance of the method of \citet{deshpande2014sparse} and Algorithm \ref{alg:disjoint.linalg} on UCI datasets.} 
\label{tab:comparison_feasiblemethods_uci3}
\end{table}

\begin{table}[h!]
\centering\footnotesize
\begin{tabular}{@{}l l l l r r r r r r r@{}} \toprule
Dataset & $p$ & $r$ &  $k_t$ & \multicolumn{3}{c@{\hspace{0mm}}}{\citet{deshpande2014sparse}} & \multicolumn{4}{c@{\hspace{0mm}}}{Algorithm \ref{alg:disjoint.linalg}}  \\
\cmidrule(l){5-7} \cmidrule(l){8-11}&  &  & & Obj. & Viol. & T(s) & UB & Obj. & Viol. & T(s)  \\
\midrule
Voice &	310 &	2 &	5 &	\textbf{0.032} &	0 &	60.46  & 0.032 & \textbf{0.032} & 0 & 16.14\\
&	&	2 &	10 &	0.063 &	0 &	54.59  & 0.064 & \textbf{0.064} & 0 & 14.22\\
&	&	2 &	20 &	0.124 &	0 &	39.19  & 0.128 & \textbf{0.127} & 0 & 28.01\\
&	&	3 &	5 &	0.047 &	0 &	1.42  & 0.048 & \textbf{0.048} & 0 & 64.25\\
&	&	3 &	10 &	0.093 &	0 &	1.75  & 0.097 & \textbf{0.096} & 0 & 66.68\\
&	&	3 &	20 &	0.184 &	0 &	1.42  & 0.192 & \textbf{0.190} & 0 & 93.89\\\midrule
Gait &	320 &	2 &	5 &	\textbf{0.030} &	0 &	1.45  & 0.031 & \textbf{0.030} & 0 & 12.51\\
&	&	2 &	10 &	0.054 &	0 &	1.45  & 0.058 & \textbf{0.055} & 0 & 10.94\\
&	&	2 &	20 &	0.089 &	0 &	1.42  & 0.110 & \textbf{0.096} & 0 & 12.27\\
&	&	3 &	5 &	0.042 &	0 &	1.85  & 0.046 & \textbf{0.045} & 0 & 47.73\\
&	&	3 &	10 &	0.078 &	0 &	1.65  & 0.087 & \textbf{0.082} & 0 & 29.65\\
&	&	3 &	20 &	0.128 &	0 &	3.04  & 0.164 & \textbf{0.135} & 0 & 46.5\\\midrule
Gastro &	466 &	2 &	5 &	\textbf{0.021} &	0 &	6.59  & 0.021 & \textbf{0.021} & 0 & 35.41\\
&	&	2 &	10 &	0.042 &	0 &	4.33  & 0.043 & \textbf{0.043} & 0 & 30.51 \\
&	&	2 &	20 &	0.083 &	0 &	3.67  & 0.086 & \textbf{0.086} & 0 & 54.39\\
&	&	3 &	5 &	0.031 &	0 &	3.26  & 0.032 & \textbf{0.032} & 0 & 97.05\\
&	&	3 &	10 &	0.061 &	0 &	3.12  & 0.064 & \textbf{0.064} & 0 & 115.2\\
&	&	3 &	20 &	0.117 &	0 &	3.11  & 0.129 & \textbf{0.128} & 0 & 153.1\\\midrule
Micromass &	1300 &	2 &	5 &	0.007 &	0 &	75.38  & 0.008 & \textbf{0.008} & 0 & 982.8\\
&	&	2 &	10 &	0.012 &	0 &	147.9  & 0.014 & {0.013} & 0 & 615.4\\
&	&	2 &	20 &	0.022 &	0 &	98.34  & 0.025 & \textbf{0.022} & 0 & 481.8\\
&	&	3 &	5 &	0.010 &	0 &	71.05  & 0.012 & {0.011} & 0 & 1023.8\\
&	&	3 &	10 &	0.018 &	0 &	71.88  & 0.021 & {0.020} & 0 & 1015.1\\
&	&	3 &	20 &	0.033 &	0 &	73.00  & 0.038 & \textbf{0.033} & 0 & 593.0\\\midrule
Avg & & & & 0.197 & 0.094 & 12.04 & 0.289 & 0.199 & 0 & 108.6\\
\bottomrule
\end{tabular}
\caption{Performance of the method of \citet{deshpande2014sparse} and Algorithm \ref{alg:disjoint.linalg} on UCI datasets (cont.)} 
\label{tab:comparison_feasiblemethods_uci3_part4}
\end{table}

\FloatBarrier
\subsection{Instance-Wise Plots of Symmetry vs. Proportion of Correlation Explained} \label{ssec:symmetryplots_appended}

\newpage
\section{Non-Redundancy of Rank Constraints in Problem \eqref{prob:spca_extended}}\label{sec:nonredundancyrankconstraints}
We claimed in Remark \ref{rem:nonredundancyrankconstraints} that the rank-one constraints in Problem \eqref{prob:spca_extended} are not redundant. We now demonstrate this by example, by providing an example where, after constraining the support pattern, Problem \eqref{prob:spca_extended} attains a different optimal value than the following optimization problem:
\begin{align}\label{prob:spca_extended2}
    \max_{\substack{\bm{Z} \in \{0, 1\}^{p \times r}:\\ \langle \bm{E}, \bm{Z}\rangle \leq k}}\max_{\bm{Y} \in \mathcal{S}^p,
    \bm{Y}^t \in \mathcal{S}^p_+} \  \langle \bm{Y}, \bm{\Sigma}\rangle \ \text{s.t.} \quad & \bm{Y} \preceq \mathrm{Diag}\left(\min\left(\bm{e}, \sum_t \bm{Z}_t\right)\right), \bm{Y}=\sum_{t=1}^r \bm{Y}^t,\\
    & \mathrm{tr}(\bm{Y}^t)=1, \ \forall t \in [r], \ Y_{i,j}^t=0 \ \text{if} \ Z_{i,t}=0, \ \forall t \in [r], i,j \in [p], \nonumber
\end{align}
where we take $r=3, p=4$, and fix the support in both problems by setting
\begin{align*}
    \bm{Z}=\begin{pmatrix} 1 & 1 & 0\\ 1 & 0 & 1\\ 0 & 1 & 1\\ 0 & 1 & 1\end{pmatrix}, \ \bm{\Sigma}=\begin{pmatrix} 1 & 1 & 0 & 0 \\ 1 & 8 & 2 & 1 \\ 0 & 2 & 2 & 1\\ 0 & 1 & 1 & 2\end{pmatrix}.
\end{align*}
Solving Problem \eqref{prob:spca_extended} by letting $\bm{Y}_t=\bm{u}_t\bm{u}_t^\top$ via Gurobi with \verb|NonConvex|$=2$ gives an optimal objective value of $12.14006$. On the other hand, solving Problem \eqref{prob:spca_extended2} via Mosek gives an optimal objective value of $12.25765$. Thus, Problems \eqref{prob:spca_extended}--\eqref{prob:spca_extended2} cannot be equivalent, as they give a different optimal objective value for a fixed binary support.



\end{document}